%% file: Review_arXiv.tex
\documentclass[american,english]{article}
\usepackage{graphicx}
\usepackage{textcomp}
\usepackage{gensymb}
\usepackage{booktabs}
\usepackage{subfig}

\usepackage{geometry}
\usepackage[round, authoryear]{natbib}

\usepackage[super]{nth}
\usepackage[export]{adjustbox}
\usepackage{authblk}

\author[1]{Fabian Wein}
\author[2]{Peter D. Dunning}
\author[3]{Juli{\'a}n A. Norato}

\affil[1]{Friedrich-Alexander-Universi\"at Erlangen-N\"urnberg (FAU), Germany, {\tt fabian.wein@fau.de}}
\affil[2]{University of Aberdeen, UK, {\tt peter.dunning@abdn.ac.uk}}
\affil[3]{University of Connecticut, USA, {\tt julian.norato@uconn.edu}}

\title{A Review on Feature-Mapping Methods for Structural Optimization}

\date{}

\input{fabian_macros.tex}

\begin{document}
\sloppy 
\title{A review on feature-mapping methods for structural optimization}

\maketitle

\begin{abstract}

In this review we identify a new category of structural optimization methods that has emerged over the last 20 years, which we propose to call \emph{feature-mapping methods}. The two defining aspects of these methods are that the design is parameterized by a high-level geometric description and that features are mapped onto a fixed grid for analysis. The main motivation for using these methods is to gain better control over the geometry to, for example, facilitate imposing direct constraints on geometric features, while avoiding issues with re-meshing.
The review starts by providing some key definitions and then examines the ingredients that these methods use to map geometric features onto a fixed-grid. One of these ingredients corresponds to the mechanism for mapping the geometry of a single feature onto a fixed analysis grid, from which an ersatz material or an immersed boundary approach is used for the analysis. For the former case, which we refer to as the pseudo-density approach, a test problem is formulated to investigate aspects of the material interpolation, boundary smoothing and numerical integration. We also review other ingredients of feature-mapping techniques, including approaches for combining features (which are required to perform  topology optimization) and methods for imposing a minimum separation distance among features.
A literature review of feature-mapping methods is provided for shape optimization, combined feature/free-form optimization, and topology optimization. Finally, we discuss potential future research directions for feature-mapping methods.

\end{abstract}

\section{Introduction}
\label{sec:intro}
Structural optimization methods can be classified into size, shape and topology optimization.  Size optimization modifies dimensions of the structure such as the cross-section of truss members or the point-wise thickness of plates.  Shape optimization modifies the boundaries of the structure, but without altering its topology, i.e., without adding or removing holes. Topology optimization can simultaneously change the shape of the structure and its connectivity.  

A key aspect of these methods is the mechanism they employ to update the analysis model upon design changes. Some methods deform the analysis mesh when the design changes. Most topology optimization methods, and some shape optimization methods, use a mesh that does not conform to the boundaries of the structure. Density-based methods, which are the most prevalent topology optimization techniques, employ a pixel/voxel representation of the design, typically based on the analysis grid. Level-set methods, which can be used both for shape and topology optimization, use the zero level-set of a function to define the structural boundaries. Density-based and level-set methods endow the optimizer with substantial freedom, rendering organic, free-form designs.

In recent years, new methods have been developed that are motivated by obtaining designs that have some desired ``high-level" geometric features (which we will define later in this article) without the need to re-mesh upon design changes. 
These methods have been largely motivated by the need to embed primitive-shaped components in free-form designs, to design structures made of stock material, to control certain dimensions of the structure, and ultimately, to provide a geometric representation that is directly understood by computer aided design (CAD) systems. In addition, these methods may represent designs with a low number of variables, which may be beneficial to, for example, the use of gradient-free optimizers. These methods build on aspects of existing techniques in (density-based) topology optimization, shape optimization and level-set methods. Despite their similarities, these methods originated independently and hence do not describe themselves as part of a common category within structural optimization. Consequently there is no commonly used label for these new methods.
We propose the term \emph{feature-mapping}, which is defined as a method that uses a high-level geometric feature parameterization that is mapped on a fixed-grid for analysis, see \secref{sec:feat-map-def} for the definition.

This review article is structured as follows. \secref{sec:define} provides some key definitions. \secref{sec:geometry_mapping} reviews methods for mapping a single feature to a fixed-grid, including pseudo-density and immersed boundary methods. For pseudo-density methods we also use a test case to investigate material interpolation, boundary smoothing and numerical integration. Methods for combining features are reviewed in \secref{sec:combination}. Some feature-mapping methods employ separation constraints, which are reviewed in \secref{sec:separation}. We then give a literature review of feature-mapping methods for: shape optimization (\secref{sec:shape_opt}), hybrid methods (which combine feature-mapping with free-form topology optimization, \secref{sec:hybrid}) and topology optimization (\secref{sec:TO-feat-map-methods}). Finally, we discuss potential future research directions for feature-mapping methods in \secref{sec:discuss}.

\section{Definitions and key components}
\label{sec:define}
\subsection{High-level geometric features}
\label{sec:high-level}
In this paper, we define a geometric feature as a geometric solid with a high-level parameterization. A geometric solid is here understood as a closed regular set of points, i.e., a set that equals the closure of its interior (cf.~\cite{Shapiro:2002:Solid}). Physically, we consider the feature can either be a solid component or a hole in a solid component. By high-level parameters, we refer to those with a direct spatial dimension associated with the feature's size, position or orientation. Examples of these parameters are the radius of a fillet, the thickness of a plate, or the location of a primitive (e.g. a bar or circle).  Notably, these high-level parameters are the ones often employed to represent solids in CAD systems. 
The advantage of having these dimensions as direct design variables is that they simplify enforcing the presence of these features and to control their dimensions, as opposed to the indirect and more verbose low-level representations of solids, such as those that are pixel or voxel-based.

\subsection{Design region, fixed-grid and moving grid}
\label{sec:region}  
The design region corresponds to the sole region of space where material can be placed. By a \emph{fixed-grid} we consider a spatial partition of the design region that remains fixed throughout design updates during the optimization process, for the purpose of linking the design description to the analysis. The term ``mesh" is perhaps more adequate to describe this spatial discretization in that it does not necessarily convey that the partition is structured (as in, for example, a ``$10 \times 10$ grid"). However, the term fixed-grid is widely used in the structural optimization and computational mechanics literature to refer to the same concept and hence we adopt it here. When used alone (i.e., `grid' or `mesh'), however, we use these terms interchangeably. 

The alternative to the fixed-grid approach is to use a moving grid, re-meshing or a combination of both. In these cases, the mesh conforms to the boundaries of the structure for any given design.  Specifically, moving grid approaches are those where the nodes in the mesh are re-positioned, but the topology of the mesh remains the same.

\subsection{Explicit and implicit geometric representations}
\label{sec:exp_imp}
For reasons of clarity we first define the terms \emph{explicit} and \emph{implicit} with respect to the geometric representation of solid objects. During the course of our review, we encountered inconsistencies in the way these terms are used in the structural optimization community. In particular, the term explicit seems to be used in several instances whenever the geometric representation employs high-level parameters, regardless of the actual representation mechanism. An explicit representation is one where points on the solid (or its boundary) are generated by a rule, whereas an implicit representation is one where a rule provides a test as to whether or not a point belongs to the solid \citep{Shapiro:2002:Solid}. For example, an explicit representation of a disc of radius $R$ centered at the origin is given by $f(t, r) = \{r \, \cos(t), r \, \sin(t) \}$, where values of the parameters $0 \leq t \leq 2\,\pi$ and $0 \leq r \leq R$ generate points within the disc. On the other hand, an implicit representation of the same disc is given by $f(\vec{x}) = \{ 1 \text{ if } \| \vec{x} \|_2  \leq R, 0 \text{ otherwise} \}$.

Density-based and classical level-set methods are implicit. In the discretized representation of density-based methods, given a point $\mathbf{x}$, the element constant pseudo-density of the element that contains $\mathbf{x}$ determines if $\mathbf{x}$ is outside or inside the solid. This point classification test is not `sharp', however, since in most density-based methods the density is a relaxed continuous variable.  
Level-set methods are implicit by definition, since the value of the level-set function at $\mathbf{x}$ determines if $\mathbf{x}$ is inside or outside the structure.  In level-set methods, a sharp representation of the boundary is available for any design throughout the optimization. This is also true when a diffuse boundary is used for the analysis, as in ersatz material methods, see \citet{Sigmund:2013:Topology}. That free-form density and level-set methods use implicit representations of the geometry is not fortuitous, since a) implicit representations more easily accommodate topological changes than explicit ones \citep{Shapiro:2002:Solid}, and b) the mesh can be used to parameterize the implicit representations, which facilitates coupling with fixed-grid analysis techniques, leading to efficient, robust methods to solve the governing equation and compute the design sensitivities.  

We consider level-set methods as those that directly represent the design using an implicit function, independent of the design update approach, as discussed by \cite{VanDijk:2013:Level}. It should be noted that some feature-mapping methods reviewed in this paper also utilize implicit functions in their formulation. Thus, they could be viewed as level-set methods. However, our aim is to emphasize the ability of feature-mapping methods to control high-level geometric features, which is not a property of level-set methods in general.

\subsection{Feature-mapping}
\label{sec:feat-map-def}
We define feature-mapping methods as those that capture high-level geometric features in their design parameterization and that map those features onto a fixed grid to perform the analysis. Note that the high-level geometric description can be either explicit or implicit.

\section{Geometry mapping to fixed-grid}
\label{sec:geometry_mapping}
\begin{figure}[ht!]
  \centering
  \includegraphics[width=0.48\textwidth]{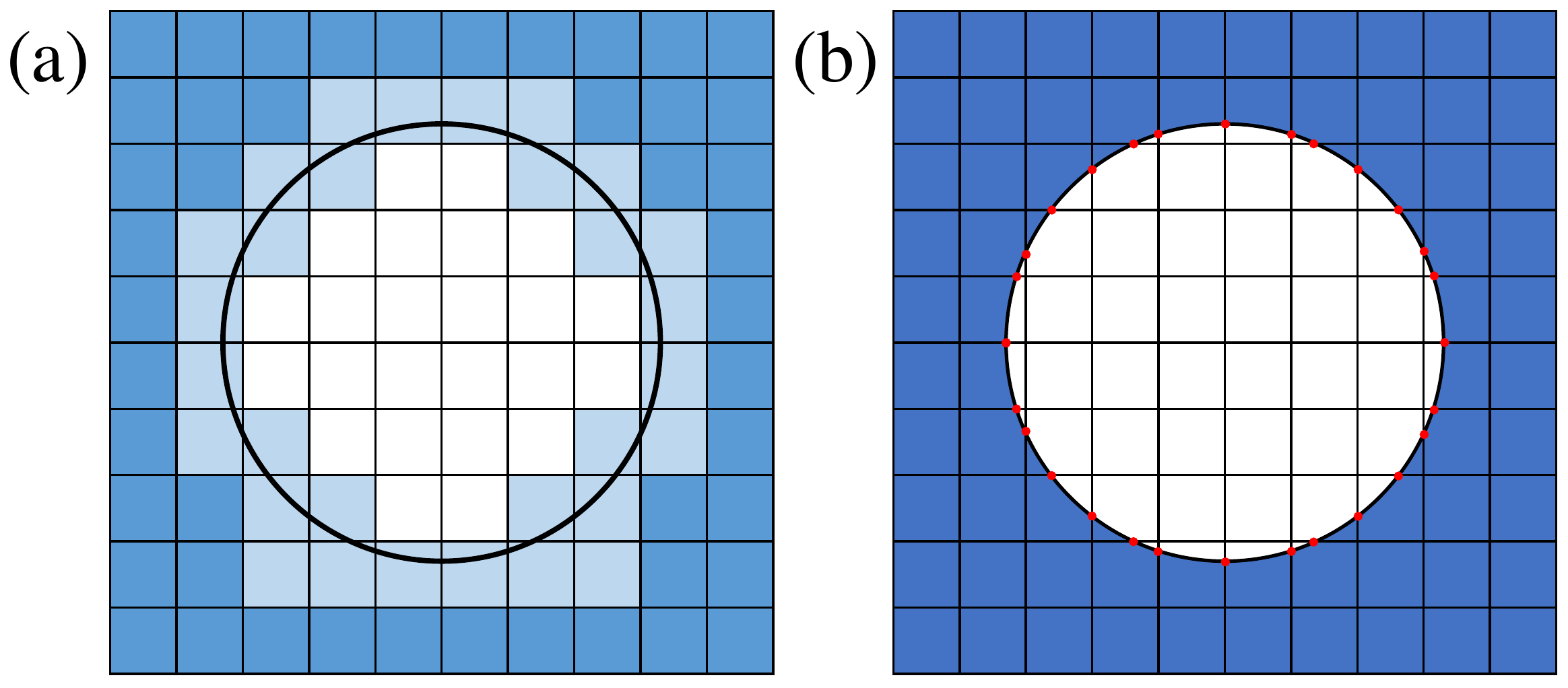}
  \caption{Fixed-grid mapping methods: a) pseudo-density, and b) immersed boundary based.}
  \label{fig:map_methods}
\end{figure}

 There are currently two main approaches to mapping high-level geometric features onto a fixed-grid for analysis: pseudo-density based mapping and immersed boundary mapping, as illustrated in \figref{fig:map_methods}. Both approaches utilize a fixed-grid for the analysis, thus circumventing the need for re-meshing during optimization.

Generally speaking, the purpose of fixed-grid analysis methods is to replace volume integrals evaluated over the structural domain $\omega$ with integrals over a domain $\Omega \supseteq \omega$ that encompasses the structure
\begin{equation}
\label{eqn:domain-integral}
    \int_{\omega} f(\Vec{x}) \,\intd \Vec{x} = \int_{\Omega} \chi_\omega(\Vec{x}) f(\Vec{x}) \, \intd \Vec{x},
\end{equation}
where $\chi_\omega$ is the \emph{characteristic} or \emph{indicator} function defined as
\begin{equation}
\label{eqn:characteristic}
    \chi_\omega(\Vec{x}) := 
    \begin{cases}
      1 &\text{if } \Vec{x} \in \omega \\
      0 &\text{otherwise}
    \end{cases}
\end{equation}
and $f$ is the domain integrand (for example, the virtual strain energy density in elasticity). 

Feature-mapping techniques that employ element-constant pseudo-densities accomplish this by replacing the corresponding element volume integral over $\Omega_e$ as
\begin{equation}
\label{eqn:element-domain-integral-density}
\int_{\Omega_e} \chi_\omega(\Vec{x}) f(\Vec{x}) \, \intd \Vec{x} \approx \rho_e \int_{\Omega_e}  f(\Vec{x}) \, \intd \Vec{x},
\end{equation}
with the element-constant pseudo-density 
\begin{equation}
\rho_e = \frac{1}{| \Omega_e |} \int_{\Omega_e} \chi_\omega(\Vec{x}).
\end{equation}
Further, $\rho_e$ might be  subject to an interpolation function $\mu$ (see \secref{sec:density_note}). Consequently, the element stiffness matrix is: $\mathbf{K}_e = \mu(\rho_e) \, \mathbf{K}_e^0$, where $\mathbf{K}_e^0$ is the `fully-solid' element matrix. This approach is also known as the \emph{ersatz material} approach. In feature-mapping methods, $\rho_e$ depends explicitly and (ideally) smoothly on the high-level parameterization of the geometric features. This mapping typically leads to elements of intermediate pseudo-density near the feature boundary, see \figref{fig:map_methods}(a) and \figref{fig:garcia_i_o_nio}. If the  mapping is differentiable, the chain rule can be readily used to obtain sensitivities with respect to the high-level geometric parameters, as we explain in \secref{sec:dens_sensitivity-analysis}. 

The immersed boundary approach employs techniques widely used in finite element analysis, such as the extended finite element method (XFEM) and isogeometric analysis, to capture sharp interfaces on a fixed-grid.  In other words, the element volume integrals are evaluated as
\begin{equation}
\label{eqn:element-domain-integral-immersed}
\int_{\Omega_e} \chi_\omega(\Vec{x}) f(\Vec{x}) \, \intd \Vec{x} = \int_{\Omega_e \cap \omega}  f(\Vec{x}) \, \intd \Vec{x}.
\end{equation}
These techniques have the advantage over pseudo-density approaches that there are no `gray regions', which require some assumption on their material properties.  For the same reason, immersed boundary methods (in principle) render more accurate analysis solutions.  These advantages come at the expense of challenges on, e.g., numerical evaluation of integrals in elements cut by the structure boundary and sensitivity calculation. These challenges and proposed solutions from the literature are discussed in \secref{sec:xfem}.

In the remainder of this section some important aspects of these approaches are discussed and examined using test case examples, in relation to their application to feature-mapping methods.

\subsection{Element-constant pseudo-density}
\label{sec:density_mapping}
Feature-mapping techniques that use pseudo-densities essentially differ in the way they compute $\rho_e$ to be used in \eqref{eqn:element-domain-integral-density}. 

An often used approach in earlier feature-mapping methods is to compute the element pseudo-density $\rho_e$ as an approximation of the volume fraction, defined as the portion of the element that intersects the feature, $|\omega \cap \Omega_e|/|\Omega_e|$.  By making simplifying assumptions about the shape of the intersected region, it is possible to use simple expressions to approximate this volume fraction.

\begin{figure}[ht!]
  \centering
  \includegraphics[width=0.48\textwidth]{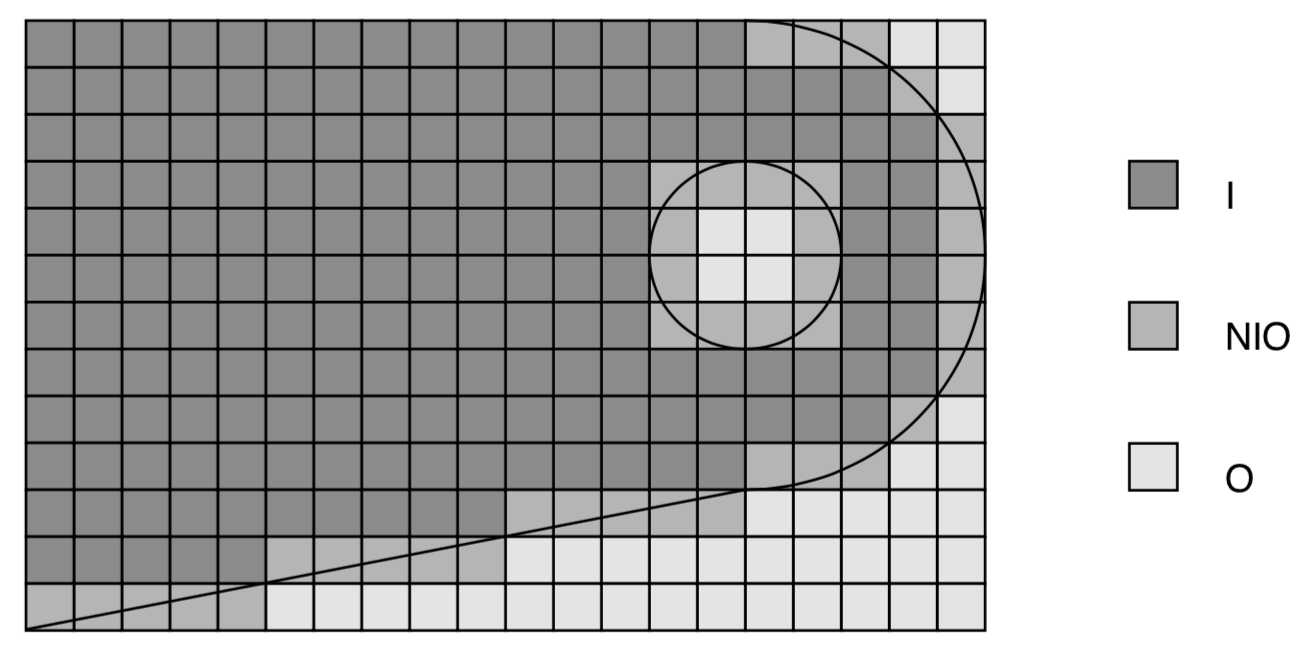}
  \caption{Fixed-grid mapping in \cite{Garcia:1999:FixedGridFE} with I (inside), O (outside) and NIO (neither inside nor outside) elements.}
  \label{fig:garcia_i_o_nio}
\end{figure}

\cite{Garcia:1999:FixedGridFE} is the earliest publication where the volume fraction approach on fixed-grids is used in the context of feature-mapping methods, see \figref{fig:garcia_i_o_nio}. Elements are classified as completely inside (I), completely outside (O), or neither inside nor outside (NIO) the structure, i.e., if the element is cut by the boundary of the geometric feature. I and O elements have pseudo-densities of 1 and $\rhomin$, respectively (with $0 < \rhomin \ll 1$ a small bound to prevent an ill-posed analysis).  The pseudo-density of NIO elements is computed as the volume fraction of the portion of the element that intersects the feature (albeit few details are provided as to how this fraction is computed), and there is no further interpolation, i.e., $\mu(\rho_e) = \rho_e$. See \secref{sec:density_note} and \secref{sec:benchmark} for discussion on the limitations of the element volume fraction approach.

As we will detail in \secref{sec:combination}, all pseudo-density techniques employ an implicit geometric representation of the feature; even when the high-level parametric representation is explicit (e.g. using a B-spline), it is first converted to an implicit representation $\phi_\omega$ satisfying the properties
\begin{equation}
\label{eqn:imp_func}
\left\{ 
  \begin{array}{l l}
   \phi_\omega (\vec{x}) > 0 , &~ \vec{x} \in \omega\\
   \phi_\omega (\vec{x}) = 0 , &~ \vec{x} \in \partial \omega\\
   \phi_\omega (\vec{x}) < 0 , &~ \vec{x} \notin \omega,
  \end{array} \right.
\end{equation}
where $\vec{x}$ is a point in the fixed-grid design domain and $\partial \omega$ is the feature boundary. Note that there is no convention in the structural optimization literature on the meaning of the implicit function sign. In \eqnref{eqn:imp_func}, we define positive values as being inside the feature. Also, the implicit function may be defined as the signed-distance function, $d(\vec{x})$, where the magnitude is computed as the shortest distance from the point $\vec{x}$ to the feature boundary $\partial \omega$ and the sign is the same as $\phi_\omega(\vec{x})$.  The distance can be computed by different approaches, which are discussed in \secref{sec:signed_distance}.

With the implicit function $\phi_\omega$ representing the feature, the Heaviside function is:
\begin{equation}
\label{eq:heaviside}
H\left(\phi_\omega(\vec{x})\right) = \left\{ 
  \begin{array}{l l}
   1  & \text{if}~ \phi_\omega (\vec{x}) \geq 0 \\
   0 & \text{if}~ \phi_\omega (\vec{x}) < 0.
  \end{array} \right.
\end{equation}

Clearly, $H$ can directly replace $\chi_\omega$ in the left-hand side \eqref{eqn:element-domain-integral-density}. The volume fraction approach directly uses $H$. However, $H$ (and $\chi_\omega$) is not differentiable, and so it is often replaced by a smooth approximation $\widetilde{H}$, see \secref{sec:boundary_smoothing} for a detailed discussion of smooth boundary modeling functions. One could consider the function values of $\widetilde{H}$ as a continuous pseudo-density field.

The element pseudo-density is then found by integrating the Heaviside function, or continuous pseudo-density field, over the element volume as
\begin{equation}
\label{eqn:rho_e_int}
   \rho_e = \frac{1}{|\Omega_e|}\int_{\Omega_e} \widetilde{H}(\phi_\omega(\Vec{x})) \, \intd x,
\end{equation}
which may be evaluated directly, or approximated by numerical integration in the form of a weighted sum
\begin{equation}
\label{eqn:rho_e_num_int}
   \rho_e = \sum_i^{N_\text{ip}} w_i \; \widetilde{H}(\phi_\omega(\Vec{x}_i)).
\end{equation}
The process of mapping conceptually consists of firstly generating the continuous pseudo-density field and then its integration. See \secref{sec:num_int} for more information on numerical integration.

Thus, feature mapping using pseudo-densities requires a choice of several key ingredients: 1) the type of material interpolation function, $\mu(\rho_e)$, 2) the form of the Heaviside (or smooth boundary) function, $\widetilde{H}$, 3) the form of the implicit function, $\phi_\omega(\vec{x})$ (signed distance or otherwise) and, 4) the integration method used to evaluate \eqnref{eqn:rho_e_int}. These ingredients are now discussed in detail, using test cases to highlight the effect of certain choices.

\subsubsection{Material interpretation of pseudo-density}
\label{sec:density_note}
Pseudo-density based feature-mapping approaches not only inherit the advantages of density-based topology optimization in terms of the in terms of the easy analysis and sensitivity computation, but they also inherit one of its challenges, which is how to interpret material properties for intermediate values of the pseudo-density.  This is dictated by the form of the function $\mu(\rho_e)$ in \eqref{eqn:element-domain-integral-density}.

\begin{figure}[ht!]
 \centering
  \includegraphics[width=.35\textwidth]{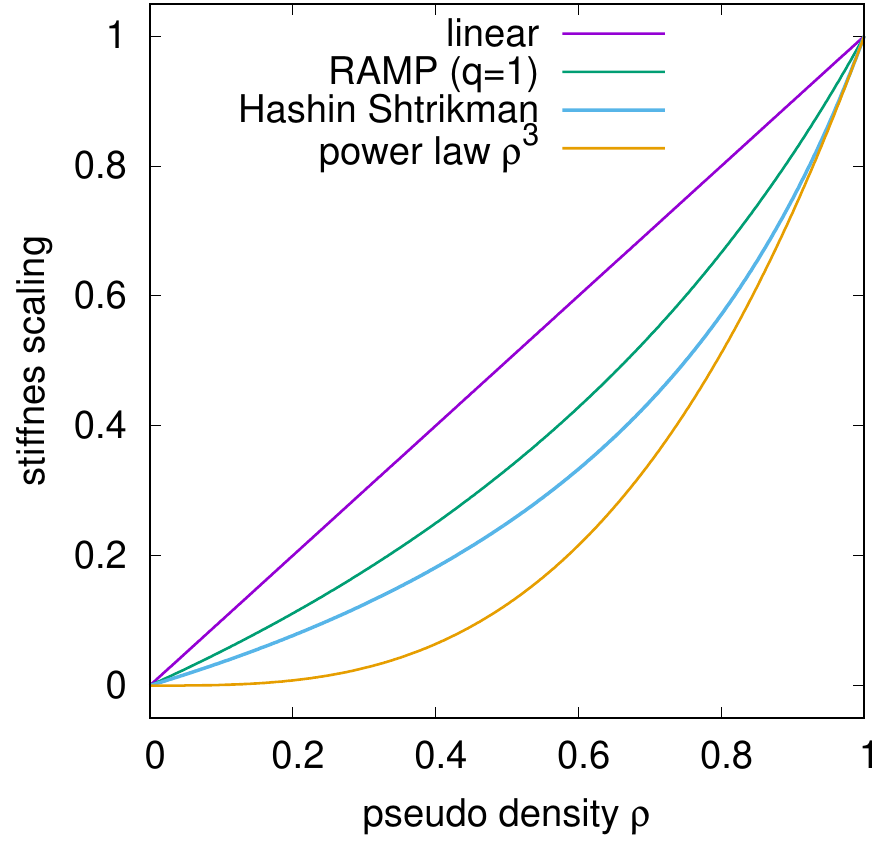}
 \caption{\label{fig:interpolation}Different material interpolation functions, see \secref{sec:density_note}. Only the Hashin-Shtrikman bound and the common SIMP power law satisfy physical limits for an isotropic material with intermediate density} 
\end{figure}

In the pioneering work on topology optimization by \cite{Bendsoe:1988:OptimalTopologies}, the stiffness properties of porous material, in between solid and void, were determined by mathematical homogenization of a periodic structure. This work showed that, due to its unfavorable stiffness-to-porosity ratio, intermediate material was barely used. This led to the famous power law approximation $\rho^p$  introduced in \cite{Bendsoe:1989:OptimalShape}, now known as the Solid Isotropic Material with Penalization (SIMP) model, to replace the homogenization step, see also \cite{Rozvany:1992:SIMP}. In \cite{Bendsoe:1999:Interpolation} it was shown that the power law 
\begin{equation}
  \mu_\text{PL}(\rho) = \rho^p    
\end{equation}
with exponent $p=3$ never overestimates the maximal physical porosity-to-stiffness relationship of isotropic material given by the upper Hashin-Shtrikman bounds, see \figref{fig:interpolation}. Homogenization and the Hashin-Shtrikman bounds show that the relative stiffness of porous isotropic material is below its volume fraction. In \figref{fig:interpolation} two further graphs are given: the linear interpolation
\begin{equation}
  \mu_\text{lin}(\rho) = \rho
\end{equation}
corresponds to the case of interpreting the volume fraction of an element covered by a geometry directly as pseudo-density, and the RAMP (Rational Approximation of Material Properties),  \cite{Stolpe:2001:RAMP}) interpolation
\begin{equation}
 \mu_\text{RAMP}(\rho) = \rho / (1 + q(1-\rho))   
 \label{eqn:ramp}
\end{equation}
is shown for $q=1$ and discussed further below. 
Both linear and RAMP (with $q=1$) interpolations are above the upper Hashin-Shtrikman bounds and therefore overestimate the stiffness of intermediate material in a non-physical way (particularly the linear interpolation). Consequently, we want to emphasize that the term \textit{penalization} in SIMP does not indicate a mathematical trick to prevent intermediate material in the optimal design, but a realistic and physical modeling of porosity to stiffness relationship.  

\subsubsection{Test problem to investigate the effect of intermediate material}
\label{sec:benchmark}

\begin{figure}[ht!]
 \centering
  \subfloat[Continuous]   
  {\includegraphics[width=.23\textwidth]{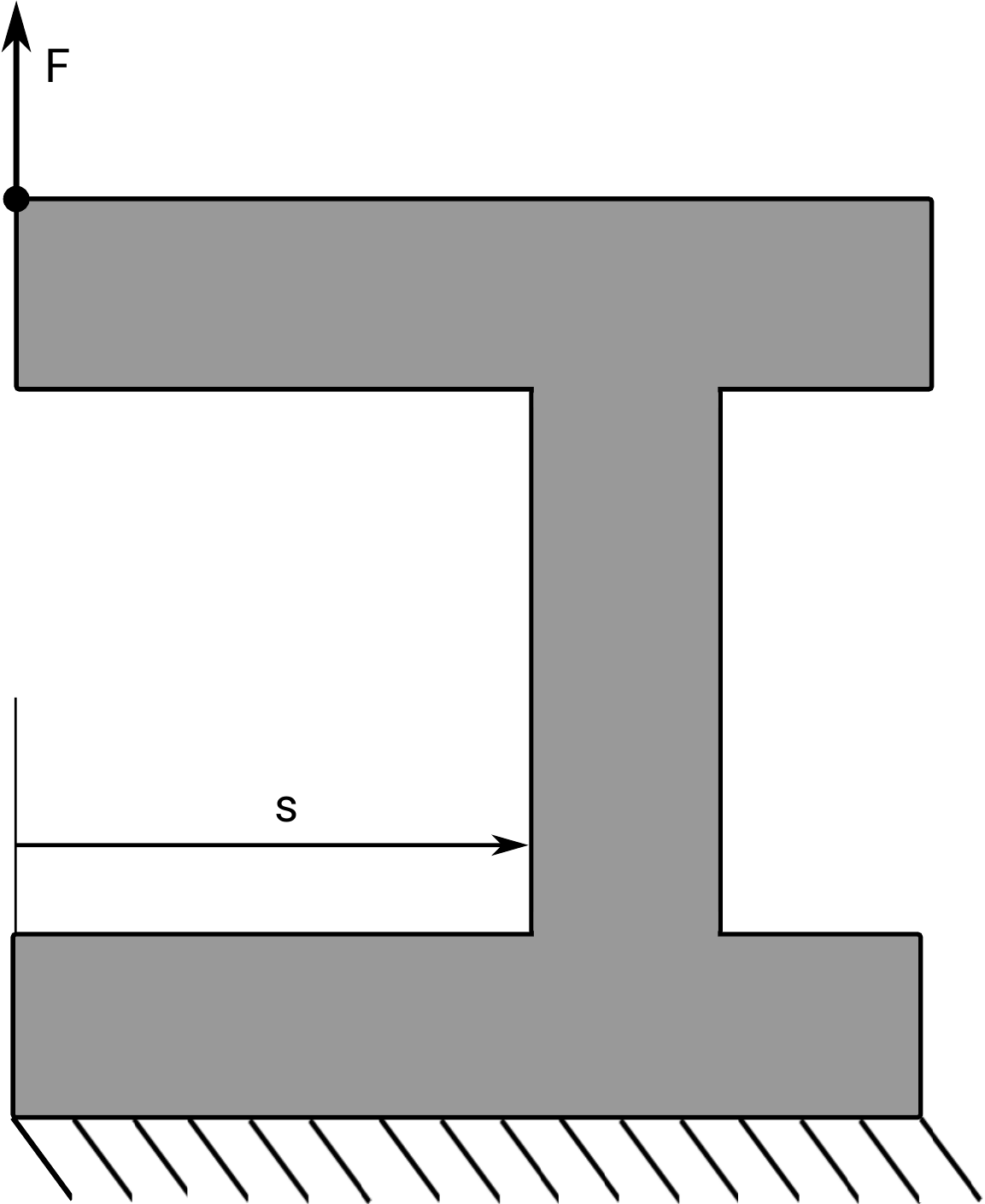}}  \,
  \subfloat[Discrete]
  {\includegraphics[width=.23\textwidth]{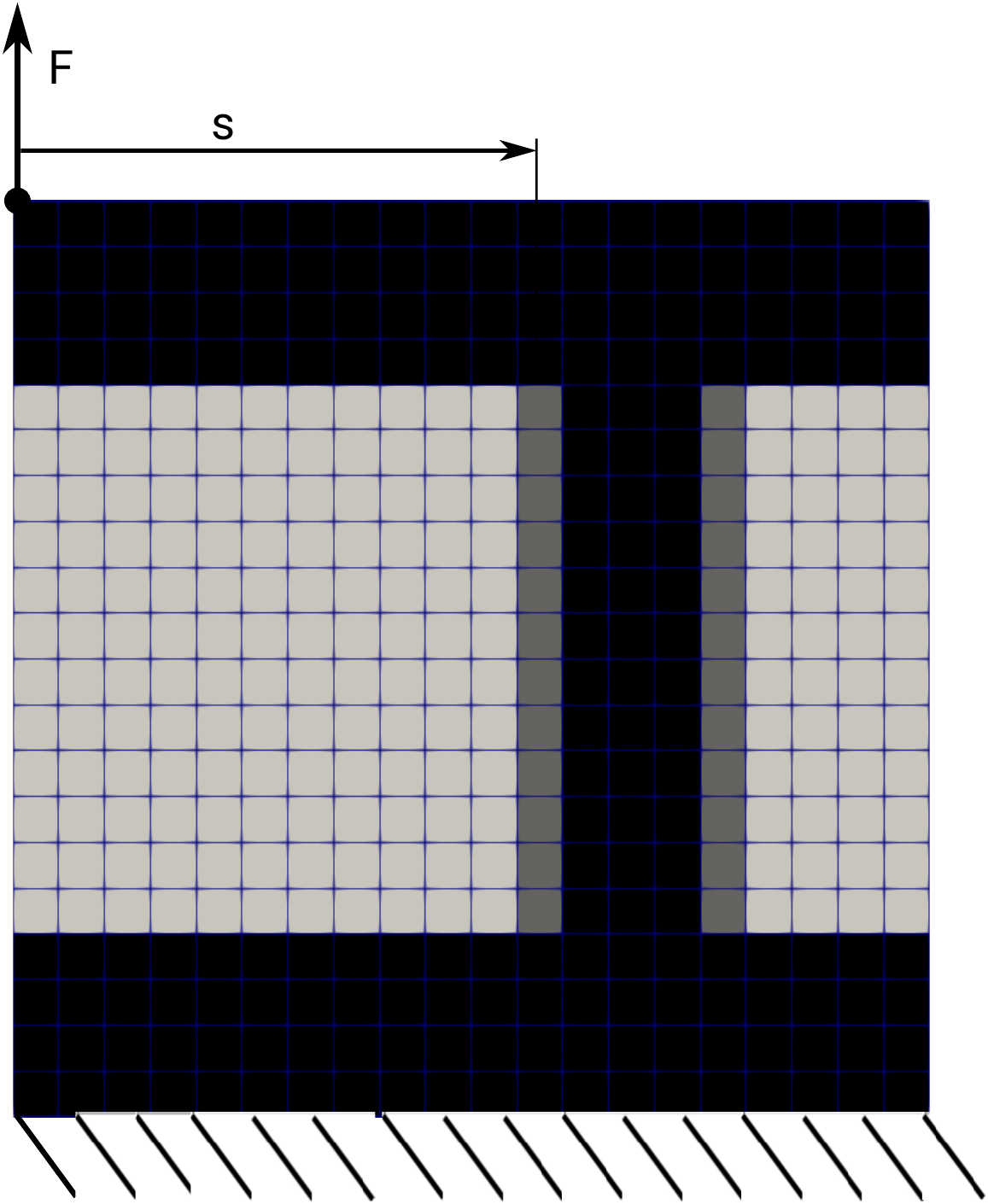}}
 \caption{\label{fig:benchmark}H-shape test problem where the horizontal position of the vertical bar is subject to the variable $s$. The force vector is applied on the upper left corner: (a) continuous setting; (b) discretized problem.}
 \end{figure}

A test problem is introduced in \figref{fig:benchmark} to investigate the effect of intermediate material modeling, particularly in the context of feature-mapping methods. The vertical bar is subject to a continuous horizontal movement with position $s$. The width of the bar is four elements. According to \figref{fig:garcia_i_o_nio} we assign a pseudo-density $\rho=1$ for elements fully contained in the bar (I), a very small value to elements fully outside the bar (O), and for the partially covered elements (NIO) a density value corresponding to the covered element volume fraction (which is the same as using the exact Heaviside in \eqref{eqn:element-domain-integral-density}). 

Note that if NIO elements along the left-hand boundary of the bar have pseudo-density $\rho$, then NIO elements along the right-hand boundary have a pseudo-density $1-\rho$, since the width of the bar is a multiple of the element size. 

The element pseudo-density $\rho_e$ is then interpolated using the functions shown in \figref{fig:interpolation}. We refer to $\mu(\rho)$ as the \textit{physical pseudo-density}, since, according to \eqref{eqn:element-domain-integral-density}, this is the element-constant material property scaling in the finite element analysis. The upper Hashin-Shtrikman bounds are given for a Poisson's ratio of 0.3 as $\mu(\rho) = \rho / (3-2\,\rho)$, see \cite{Bendsoe:1999:Interpolation}. To measure the grayness of the boundary elements we introduce
\begin{equation}
g(\mu(\rho_e)) = 4 \, \mu(\rho_e) \,(1-\mu(\rho_e)) ,
\label{eqn:grayness} 
\end{equation}  
where $\mu(\rho_e)=0.5$ results in the highest grayness value of 1.0.

\begin{figure}[ht!]
 \centering
  \subfloat[Grayness]   
  {\includegraphics[width=.4\textwidth]{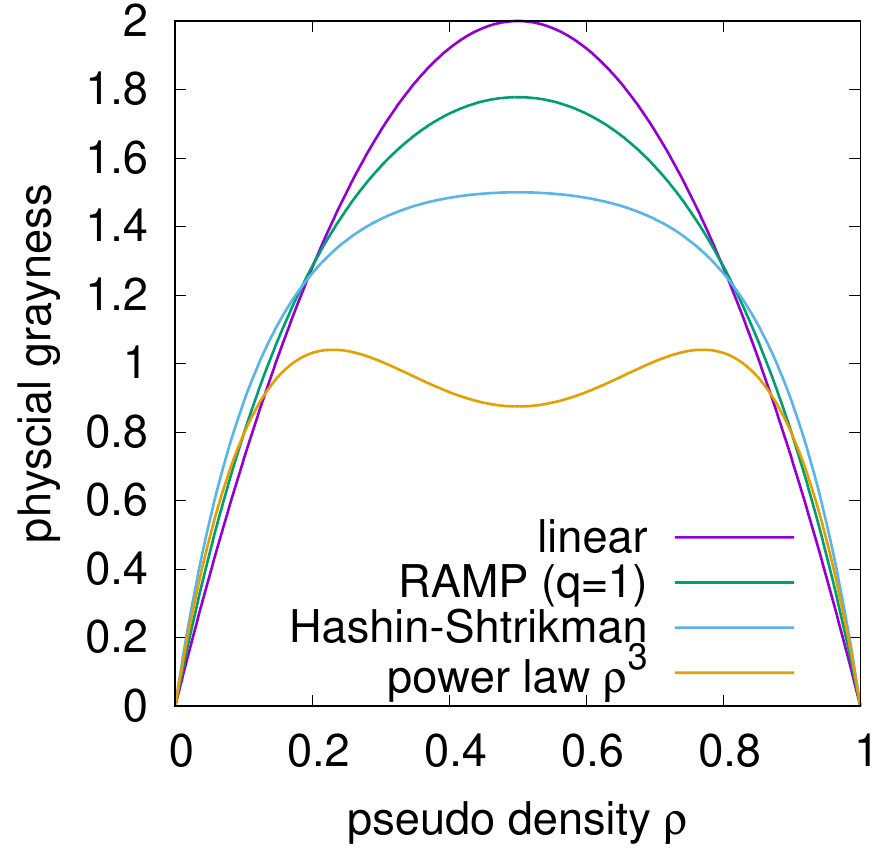}} \,
  \subfloat[Compliance]   
  {\includegraphics[width=.4\textwidth]{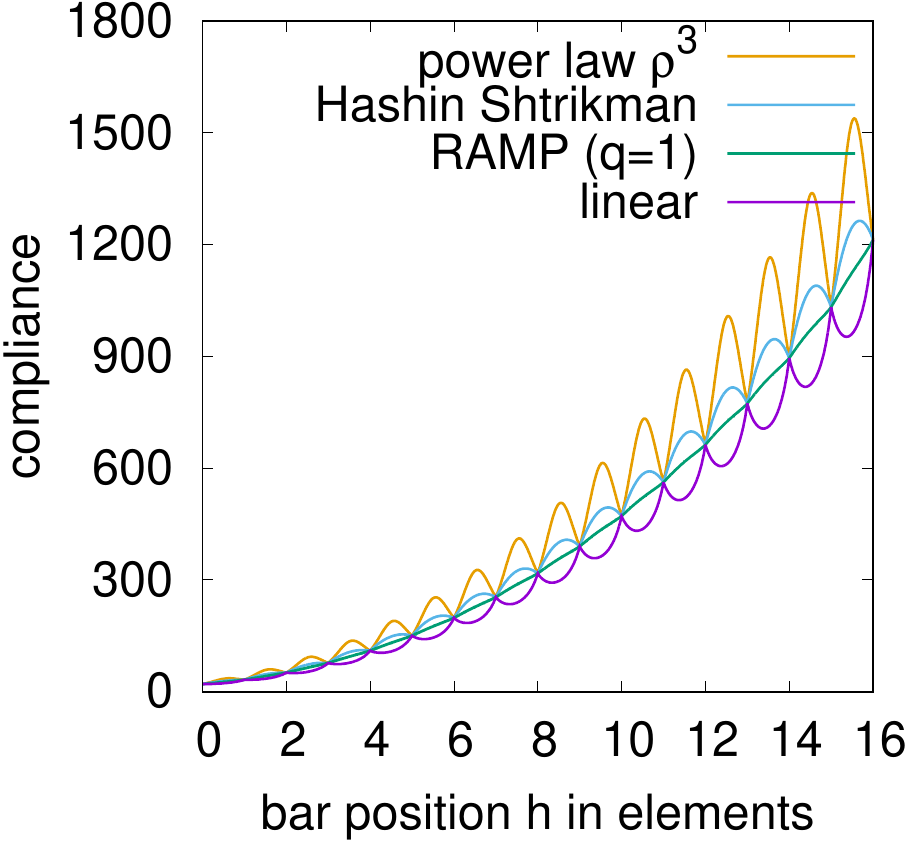}}
 \caption{\label{fig:grayness}(a) Sum of grayness values of the physical density $g(\mu(\rho))$ \eqnref{eqn:grayness} for left and right edges of the vertical bar, corresponding to moving $s$ by one element in \figref{fig:benchmark} for different interpolation functions; and (b) their corresponding structural compliance values.}
 \end{figure}

The compliance for the test problem in \figref{fig:benchmark}, is evaluated for the material interpolation functions from \figref{fig:interpolation} and the results are shown in \figref{fig:grayness}(b). The linear material interpolation shows improved (lower) compliance when the bar edge is positioned between elements, giving a high grayness value, see \figref{fig:grayness}(a). Peaks of poorer (high) compliance are seen when the bar is aligned with element edges, resulting in no gray boundary elements. Using the upper Hashin-Shtrikman bounds for material interpolation shows an increased compliance when gray elements are involved, reflecting the inefficient stiffness of porous structures in reality. This effect is amplified for the classical SIMP power law. Realistic compliance values are obtained only when the bar edges align with element boundaries, as there are no intermediate, gray densities.

These three material interpolation functions show that the compliance is non-monotonic with respect to the bar position $s$. In an optimization problem, the linear interpolation function will likely favor intermediate bar positions, while the Hashin--Shtrikman bounds and the power law will likely favor bar positions aligned to mesh elements - hence the problem becomes somewhat mesh-dependent. We note that this effect is caused only by the process of mapping the feature onto the fixed-grid using pseudo-densities, since analysis of this example with a conforming mesh would render a monotonic compliance curve. Interestingly, the RAMP interpolation function with parameter $q=1$, which is $\mu(\rho)=\frac{\rho}{2-\rho}$, exhibits an almost monotonic compliance with respect to the design change for this test problem.
The authors are only aware of one work in the literature where this interpolation function is used in feature-mapping methods \citep{zhang2017stress}. However, RAMP is not introduced with the purpose of avoiding mesh dependency, but to favor reintroduction of geometric features during optimization. 

In addition to non-monotonicity, \figref{fig:grayness}(b) reveals another aspect of the feature-mapping: the compliance is non-smooth, as it exhibits `kinks' whenever the vertical boundaries of the moving bar coincide with element boundaries. This non-smoothness is present even for the seemingly smoother RAMP interpolation.

\subsubsection{Principal boundary modeling approaches}
\label{sec:gray_boundary}

Mapping a feature to a fixed analysis grid requires modeling the boundary when the design is not exactly aligned to the mesh. However, as shown in the previous section, this may lead to non-monotonicity, non-smoothness and mesh-dependency. The results in \secref{sec:density_note} assume the boundary is modelled by an exact Heaviside function \eqnref{eq:heaviside}. 
In this section, we investigate the effect of using a smoothed Heaviside, or boundary smoothing, on the test problem. 

\begin{figure}[ht!]
 \centering
 \includegraphics[width=.75\textwidth]{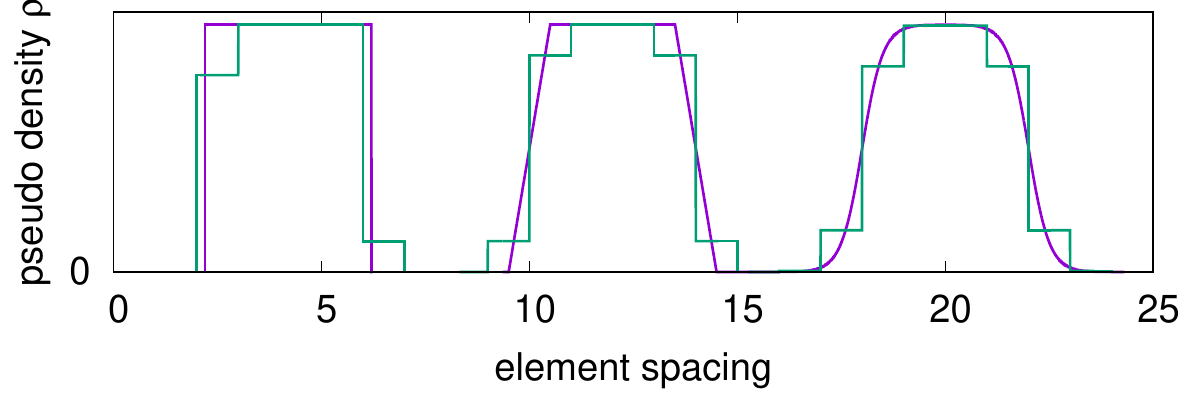} 
 \caption{\label{fig:shapes}1D boundary modeling (in magenta): exact Heaviside function (left), piecewise linear \eqnref{eq:linearHS} (center) and $\text{tanh}$ smoothing \eqnref{eqn:tanh} (right) and their element-constant pseudo-density values (in green).}
\end{figure}

First, we consider a 1D model of the bar cross-section in the test case example of \figref{fig:benchmark}. This feature is modelled by assigning a pseudo-density $\rho=1$ to any point inside the feature and a very small value $\rhomin$ to points outside the feature (similar to the characteristic function as defined in \eqnref{eqn:characteristic}). In \figref{fig:shapes}, three approaches to model the transition between material and void across the boundary are shown. The first is an exact Heaviside function, \eqnref{eq:heaviside}, as in \secref{sec:density_note}; this function is discontinuous, a fact we denote as the function being $C^{-1}$. The other functions are: a continuous, but non-differentiable piecewise linear function, which is $C^0$
\begin{equation}
\label{eq:linearHS}
\widetilde{H}_{\text{lin}}(d(\vec{x}), h) :=
\begin{cases}
\rhomin & \text{if } d(\vec{x}) < -h \\
 (1-\rhomin) \; \frac{d(\vec{x})}{2\,h} + \frac{1+\rhomin}{2} & \text{if }  |d(\vec{x})| \leq h \\
1 & \text{if } d(\vec{x}) > h,
\end{cases}
\end{equation}
where $h$ defines the size of the transition zone between material and void; and a $\text{tanh}$-like function, which is $C^\infty$, as used in \cite{wein2018combined} 
\begin{equation}
  \label{eqn:tanh}
  \widetilde{H}_\text{tanh}(d(\vec{x}), \beta) = (1-\rhomin) \left(1 - \frac{1}{e^{\beta\,d(\vec{x})} + 1} \right) + \rhomin,
\end{equation}
where $\beta$ is a parameter that controls the size of the transition zone. Note that we have chosen the signed-distance function $(d(\vec{x}))$ as the implicit function in the above equations. For further discussion on how the choice of implicit function affects boundary mapping, see \secref{sec:boundary_smoothing}.

\begin{figure}[ht!]
 \centering
\includegraphics[width=.4\textwidth]{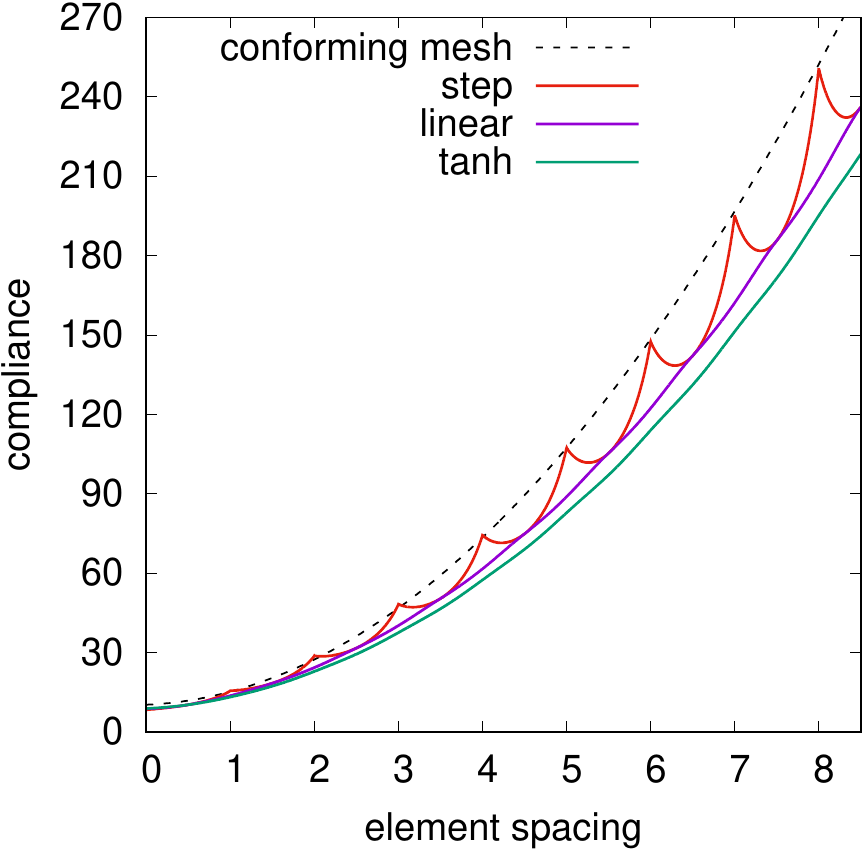} 
 \caption{\label{fig:shape_move}Compliance values for the boundary modeling functions in \figref{fig:shapes} applied to the test problem in \figref{fig:benchmark}. Due to the linear density-to-stiffness interpolation the compliance is underestimated, see \figref{fig:grayness}(b). The compliance for a conforming mesh is plotted for reference. } 
\end{figure}

To obtain the element constant pseudo-densities for the fixed analysis grid we evaluate \eqnref{eqn:rho_e_int}. 
For simplicity, we assume the linear material interpolation model from \figref{fig:interpolation}. The corresponding element-constant pseudo-densities are shown in \figref{fig:shapes}. Again, the element grayness values depends on the alignment of the boundary modeling function with respect to the mesh (see \figref{fig:grayness}(a)). However, the sum over all elements attains a constant value for the piecewise linear and $\text{tanh}$ functions, see \figref{fig:num_int_slope}.
Applying these functions to the test problem, we get the compliance values shown in \figref{fig:shape_move}. The result for the exact Heaviside function has already been given in \figref{fig:grayness}(b); the piecewise linear and $\text{tanh}$-functions result in visually smooth compliance functions that are artificially good (low). The smoothness results from the boundary smoothing, while the artificially low compliance results from the linear material interpolation.

Upon evaluation of \eqnref{eqn:rho_e_int}, the resulting pseudo-density $\rho_e(s)$ and hence the compliance become $C^0$ for the exact Heaviside function, $C^1$ for the piecewise linear function and stays $C^\infty$ for the tanh function, with respect to a change in $s$. For optimization we require $C^1$, hence the piecewise linear boundary modeling function is sufficient. For more discussion on the effect of numerical integration, see \secref{sec:num_int}.

\subsubsection{Further smooth boundary modeling approaches}
\label{sec:boundary_smoothing}

In this section, we review further smooth boundary modeling approaches used in feature-mapping methods. 
Typically they are piecewise defined with a transitioning zone controlled by the parameter $h$, and transition function $\sigma_\text{trans}$ as
\begin{equation}
\label{eq:smoothHS}
\widetilde{H}(\phi_\omega(\vec{x}), h) :=
\begin{cases}
\rhomin & \text{if } \phi_\omega(\vec{x}) < -h \\
 \sigma_\text{trans}(\phi_\omega(\vec{x}), h) & \text{if }  |\phi_\omega(\vec{x})| \leq h \\
1 & \text{if } \phi_\omega(\vec{x}) > h.
\end{cases}
\end{equation}

For example, the transition function for the piecewise linear boundary model \eqnref{eq:linearHS} is
\begin{equation}
  \label{eqn:transition_lin}
\sigma_\text{lin} = (1-\rhomin) \; \frac{\phi(\vec{x})}{2\,h} + \frac{1+\rhomin}{2}.
\end{equation}

A common choice of transition function is based on a spline representation as a cubic function with zero slope on both sides of the transition zone, which is used by e.g. \cite{zhang2016new} and \cite{Dunning2018}
\begin{equation}
  \label{eqn:transition_poly}
  \sigma_\text{poly} = \frac{3\,(1 - \rhomin)}{4}\left( \frac{\phi(\vec{x})}{h}- \frac{\phi(\vec{x})^3}{3 \, h^3} \right) + \frac{1+\rhomin}{2}.
\end{equation}
Another choice is based on a trigonometrical function:
\begin{equation}
  \label{eqn:transition_cos}
  \sigma_\text{cos} = \frac{(1-\rhomin)}{2} \, \cos \left( \left(\frac{\phi(\vec{x})}{2h}-\frac{1}{2} \right) \, \pi \right) + \frac{1+\rhomin}{2}.
\end{equation}
Note that the $\text{tanh}$-like function \eqnref{eqn:tanh} has no finite transition zone.

In this section, the more general implicit function, $\phi(\vec{x})$, is used instead of the signed-distance function, $d(\vec{x})$, as not all feature-mapping methods use a signed-distance function. If a signed-distance function is used, then the magnitude of the implicit function spatial gradient is one, $||\nabla d||=1$, and the width of the transition zone between solid and void is defined as: $w=2\,h$. In \secref{sec:num_int} we also introduce the discrete element transition zone $w_b$. If $||\nabla \phi||<1$, the transition zone will be stretched $w>2\,h$. Conversely, if $||\nabla \phi||>1$, then the transition zone will be compressed $w<2\,h$. This issue is discussed and investigated by \cite{zhou2016feature}, where it is argued that a signed-distance function should be used to avoid issues caused by a varying spatial gradient of $\phi$ around the feature boundary, as this influences the accuracy of the structural response and gradient computation. See also \secref{sec:num_int}. For further discussion on computing the signed-distance function, see \secref{sec:signed_distance}.

\secref{sec:benchmark} shows the limitations of obtaining the element pseudo-density $\rho_e$ as volume fraction of a grid cell $e$ covered by a non-smoothed feature. A particular issue demonstrated in \cite{Norato:2004:Shape} is that the volume fraction calculation becomes non-differentiable if a portion with non-zero measure of the feature boundary coincides with the element boundary---for instance, in the test problem of \figref{fig:benchmark} when the sides of the vertical bar align with element boundaries. This problem can be readily circumvented by using a circular (2D) or spherical (3D) sampling window (instead of the element itself) to compute the volume fraction, and by linearizing the boundary of the feature within the sampling window (cf.~\cite{norato2015geometry}). This leads to a closed formula for the transition function, which is given for 2D as
\begin{equation}
\label{eq:geom_proj}
\sigma_\text{circ} =
 \frac{1}{\pi} \left[ \cos^{-1}\left( \frac{-d(\vec{x})}{h} \right) +  \frac{d(\vec{x})}{h} \,\sqrt{ 1 - \left(\frac{d(\vec{x})}{h}\right)^2} \right]
\end{equation}
within the framework of \eqnref{eq:smoothHS}. However, note that numerical integration in \eqnref{eqn:rho_e_int} is not used in this case, as \eqnref{eq:geom_proj} is derived from an exact analytical integration of the volume fraction of the linearized feature boundary within the circular sampling window and thus only requires the signed-distance information from the element center to the feature boundary.

\begin{figure}[ht!]
  \centering
  \includegraphics[width=0.65\textwidth]{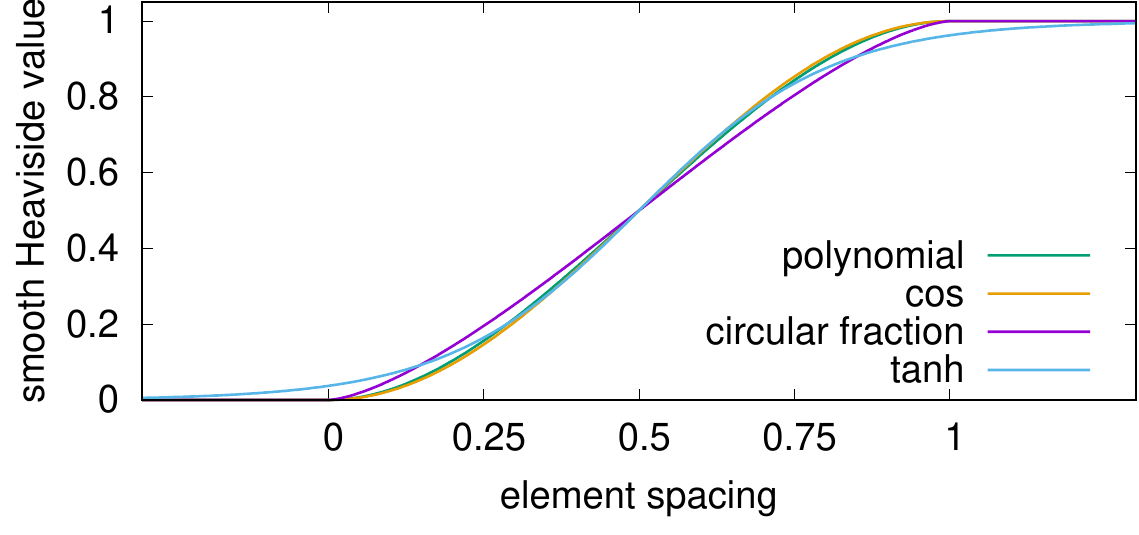}
  \caption{\label{fig:cut_smooth}Smooth boundary modelling by polynomial \eqnref{eqn:transition_poly} (almost identical to  cos \eqnref{eqn:transition_cos}) and volume fraction with circular sampling window \eqnref{eq:geom_proj} with transition zone $2\,h$ equals one element width. The $\text{tanh}$-based function \eqnref{eqn:tanh} is plotted with $\beta=6.5$. All functions are plotted assuming a signed-distance function.}
\end{figure}

A comparison of boundary smoothing functions is shown in \figref{fig:cut_smooth}.

\subsubsection{Sensitivity Analysis} 
\label{sec:dens_sensitivity-analysis}

One of the appealing features of the element pseudo-density approach in feature-mapping methods is that, as in density-based topology optimization, the computation of design sensitivities is much simpler than for approaches that must compute boundary sensitivities (as in some level-set methods). Moreover, as we will show in this section, the computation of sensitivities is closely connected to that of density-based methods.

Sensitivity analysis in density-based topology optimization is well established. It can be readily performed on the discretized algebraic system resulting from a finite element analysis for a wide range of functions, and even multiphysics problems fit one of the known generalized derivations, see \cite{Bendsoe:2003:Book}. 

We briefly review sensitivity analysis for standard density-based topology optimization and consider the easy static case, where the finite element system matrix $\vec{K}$ depends explicitly on the vector of element pseudo-densities $\brho$, the state solution $\vec{u}$ depends only implicitly on $\brho$ and the boundary conditions are assumed to be design-independent. The system of linear equations arising from the finite element discretization reads
\begin{equation}
    \vec{K}(\brho)\,\vec{u}(\brho) = \vec{f}.
\end{equation}
Using adjoint differentiation (e.g. \cite{Troeltzsch:2010:OptPDE}), the sensitivity of a function $J(\brho, \vec{u}(\brho))$ with respect to an element pseudo-density can be written as
\begin{equation}
   \total{J}{\rho_e} = \drho{J} + \blmbd(\brho)^\top \drho{K_e}\,\vec{u}(\brho).
   \label{eq:dJ-density}
\end{equation}
The partial derivative $\drho{J}$ is for many functions zero, $\drho{\vec{K}_e}$ is trivial to obtain and $\blmbd$ solves the adjoint problem
\begin{equation}
    \vec{K}(\brho)\,\blmbd(\brho) = - \left(\sens{J}{\vec{u}}\right)^\top.
    \label{eq:adjoint-problem}
\end{equation}
Notably, the adjoint solution $\blmbd$ is independent of the design parameterization, because the pseudo-load in \eqref{eq:adjoint-problem} does not depend explicitly on the design variables. Consequently, the adjoint solution needed for feature-mapping methods is the same as the one obtained for other topology optimization techniques.

Once the adjoint solution is computed, feature-mapping methods with pseudo-densities only need to compute the derivative of the boundary mapping function to obtain derivatives of pseudo-densities with respect to the high-level design parameters, $s_j$. The final derivative is obtained by the chain rule as
\begin{equation}
   \total{J}{s_j} = \sum_e^{N_\text{e}} \left\{ \left[ \sens{J}{\rho_e}  + \sens{\mu}{\rho_e}  \blmbd(\brho)^\top \vec{K}_e \,\vec{u}(\brho) \right] \sens{\rho_e}{s_j} \right\}.
   \label{eqn:dJ_ds}
\end{equation}

\begin{figure}[ht!]
  \centering
  \includegraphics[width=0.24\textwidth]{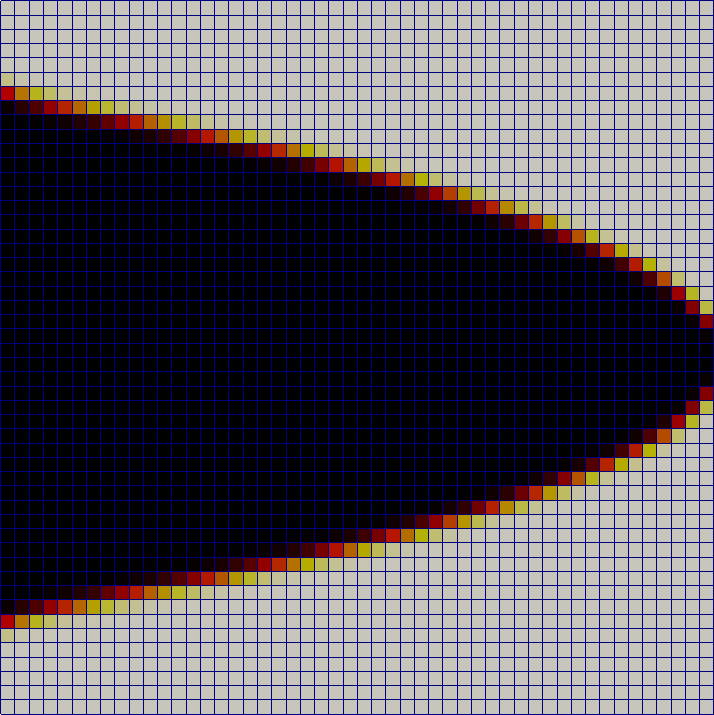} 
  \includegraphics[width=0.34\textwidth]{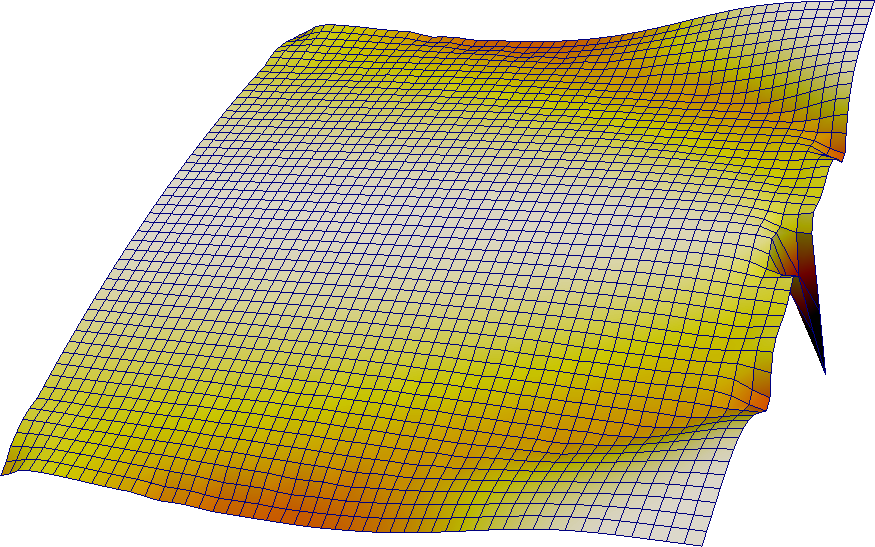} 
  \includegraphics[width=0.34\textwidth]{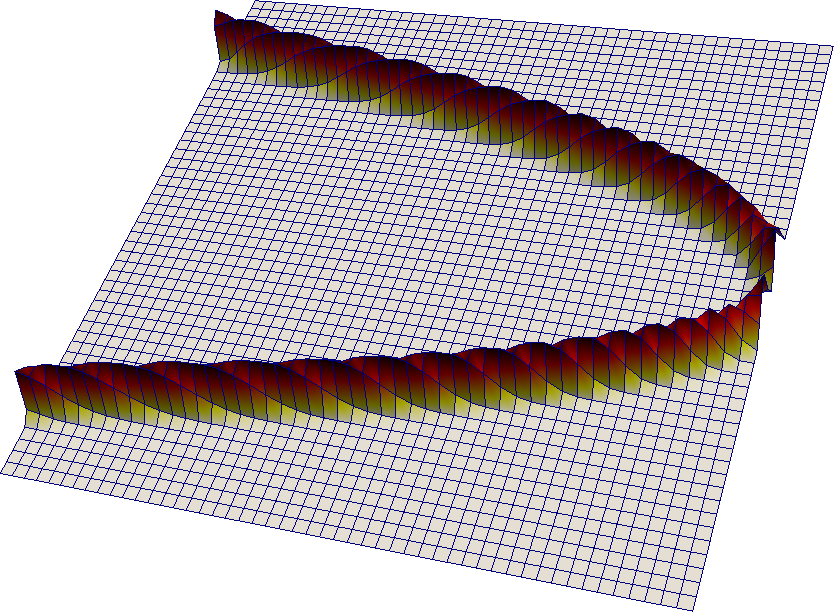} 
  \caption{For a smoothly mapped shape (left) the compliance sensitivity with respect to pseudo-density \eqnref{eq:dJ-density} is shown in the center. The element-wise summand of the sensitivity with respect to shape variables \eqnref{eqn:dJ_ds} is non-zero only where the pseudo-density has intermediate values (right). From \cite{wein2018combined}.}
  \label{fig:mapped_gradient}
\end{figure}

Note that the boundary modeling function \eqref{eq:smoothHS} is constant outside of the transition region, hence $\sens{\widetilde{H}}{s_j}=0$ in the void region (i.e.., $\widetilde{H} = \rhomin$) and the solid region (i.e.., $\widetilde{H}=1$). It thus follows from \eqref{eqn:rho_e_num_int} that $\sens{\rho_e}{s_j}$ is non-zero only in regions with intermediate pseudo-density values, namely in the gray regions around the boundaries of the structure, see \figref{fig:mapped_gradient}. The choice of width $h$ of the smoothing functions in \secref{sec:boundary_smoothing} controls the amount of information collected from \eqnref{eq:dJ-density}. 

\subsubsection{Numerical integration of the boundary mapping function}
\label{sec:num_int}
In principle, density-based feature-mapping requires the element-constant pseudo-density to be found by integrating the smoothed Heaviside, or boundary mapping function \eqnref{eqn:rho_e_int}. In the test case example from \figref{fig:benchmark}, the vertical feature is aligned with the fixed-grid. This effectively makes the volume integral of the boundary modeling function 1D. Thus, analytical integration is reasonably straight-forward and is used to generate the results above.
However, analytical integration can become involved in two and three dimensions. Therefore, many methods compute the pseudo-density by numerical integration as a weighted sum via ~\eqnref{eqn:rho_e_num_int}.

The boundary modeling function \eqnref{eq:smoothHS} influences the choice of quadrature rule and the number of sampling points. We now examine the effect of the number of sampling points when using Newton-Cotes formulae to evaluate \eqnref{eqn:rho_e_num_int}. Note that zero-degree quadrature corresponds to midpoint integration (i.e., the value of the function at the element center) and first-degree quadrature corresponds to the trapezoidal rule (i.e., the average of function values at the corner positions of the element). The number of integration points for the Newton-Cotes formula with degree $\text{deg}$ is $N_\text{ip} = (\text{deg}+1)^{\text{dim}}$. 

The investigation uses the piecewise linear \eqnref{eq:linearHS} and polynomial \eqnref{eqn:transition_poly} boundary modeling functions with transition zone $w$ of one element. The $\text{tanh}$-like function \eqnref{eqn:tanh} is also included, with a similar maximal slope. The signed-distance implicit function is used. It is clear that element-wise numerical integration of the density function does not increase the regularity with respect to the shape variables. In particular the piecewise linear boundary modeling function  \eqnref{eqn:transition_lin} stays non-continuous differentiable. Thus, this combination is not suitable for gradient based optimization, but we feel it worth including in the discussion.

\begin{figure}[ht!]
	 \centering
  \subfloat[linear \eqnref{eqn:transition_lin}]
   {\includegraphics[width=0.36\textwidth]{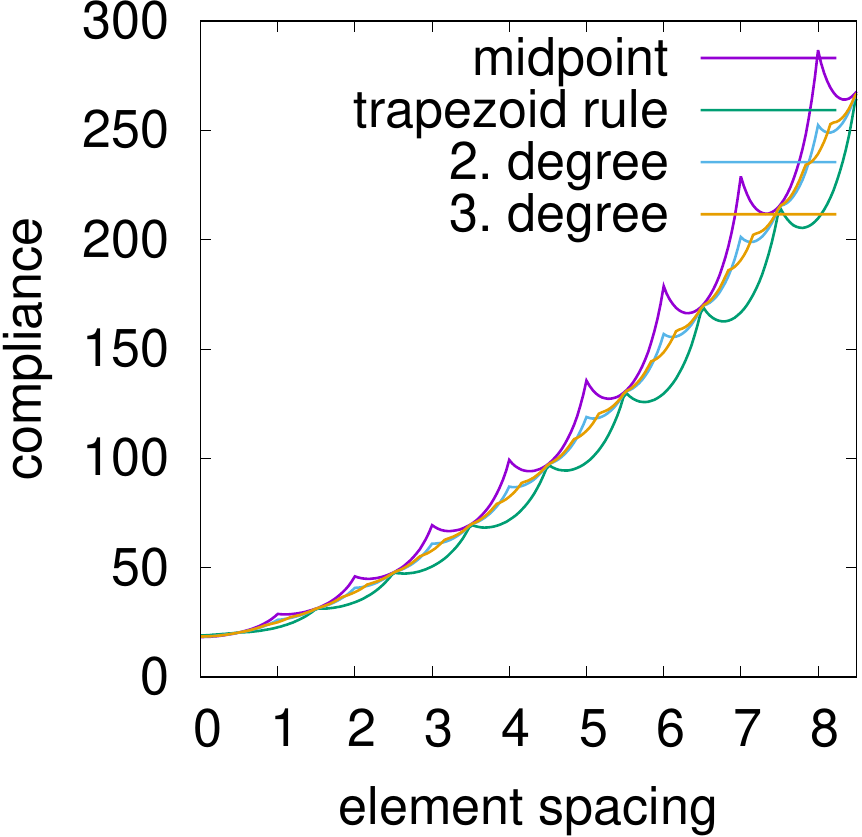}}
   \subfloat[polynomial \eqnref{eqn:transition_poly}]
   {\includegraphics[width=0.29\textwidth]{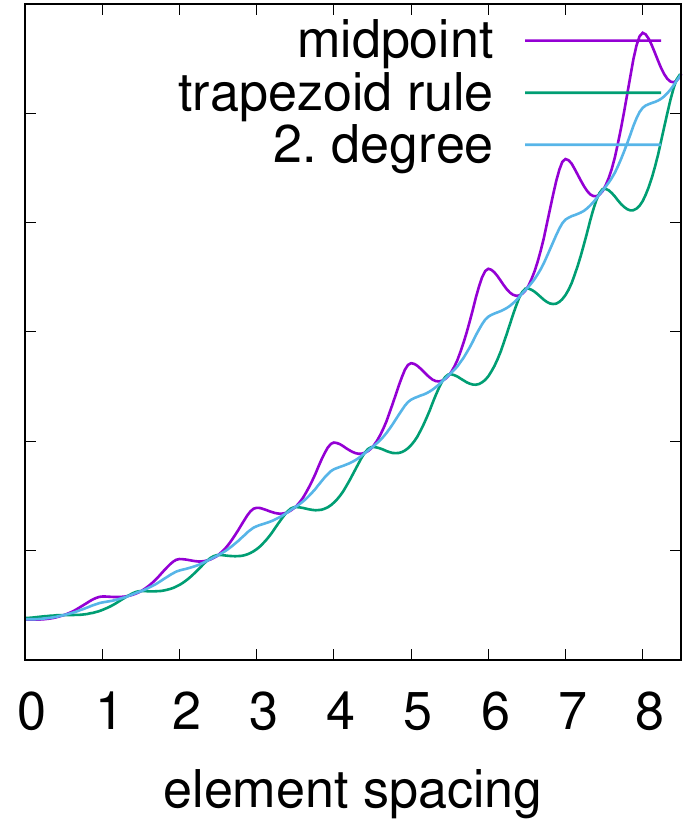}}
   \subfloat[$\text{tanh}$ \eqnref{eqn:tanh}]
   {\includegraphics[width=0.29\textwidth]{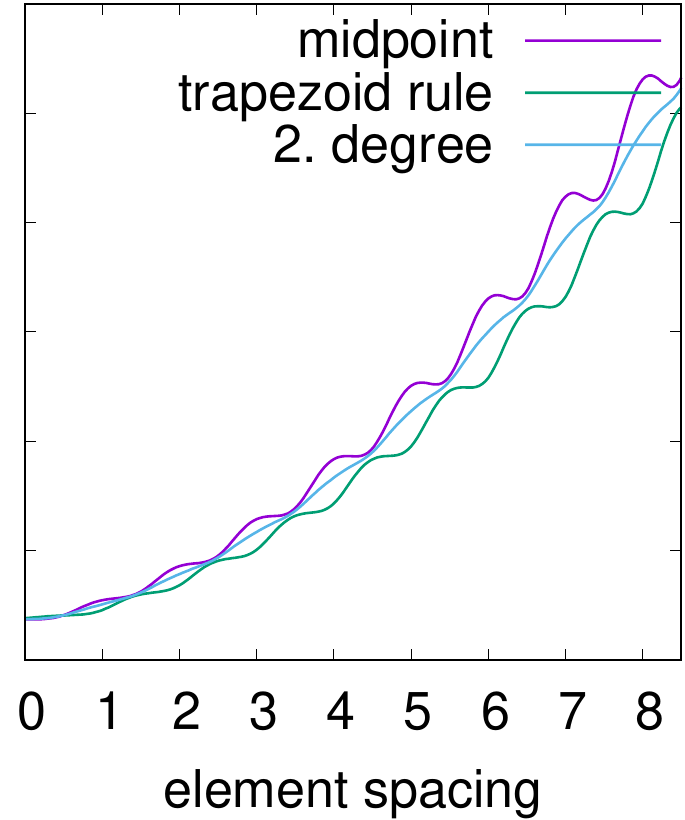}}
  \caption{We perform the similar experiment as in \figref{fig:shape_move} but this time with numerical integration via low order closed Newton-Cotes formulae. See \figref{fig:shapes}.}
  \label{fig:num_int_move}
\end{figure}

\figref{fig:num_int_move} clearly shows that, for the test problem, the smoothing effect shown with analytical integration in \figref{fig:shape_move} is lost when the number of sampling points in numerical integration is too low. Although, it should be noted that the test case is selected to reveal extreme response. 

\begin{figure}[ht!]
  \centering
  \includegraphics[width=0.33\textwidth]{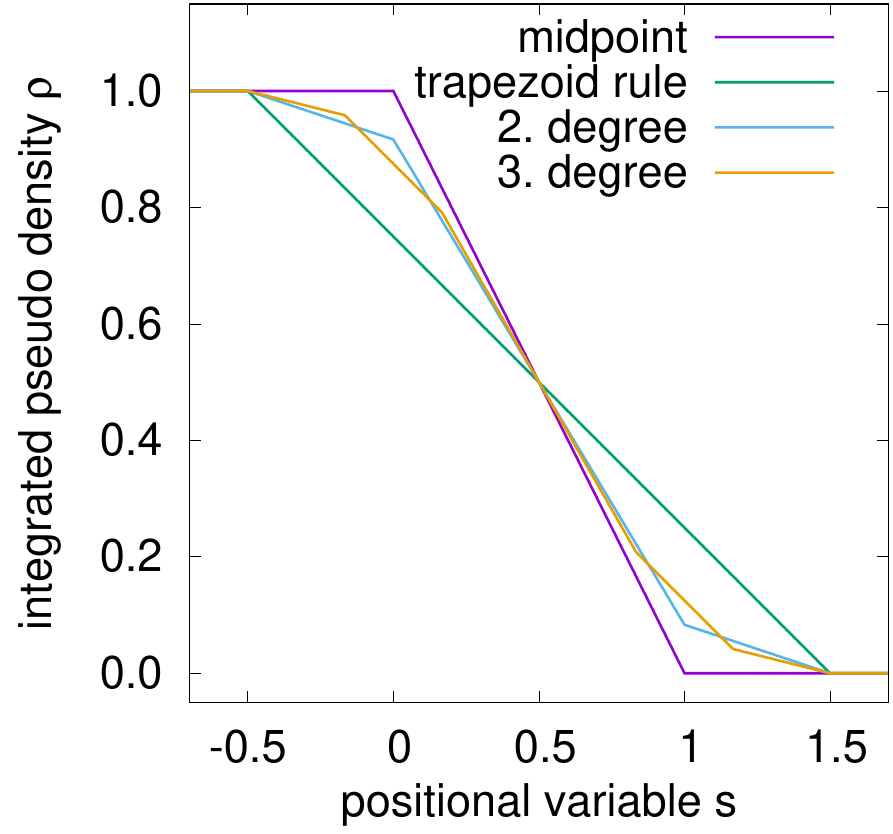}
  \includegraphics[width=0.33\textwidth]{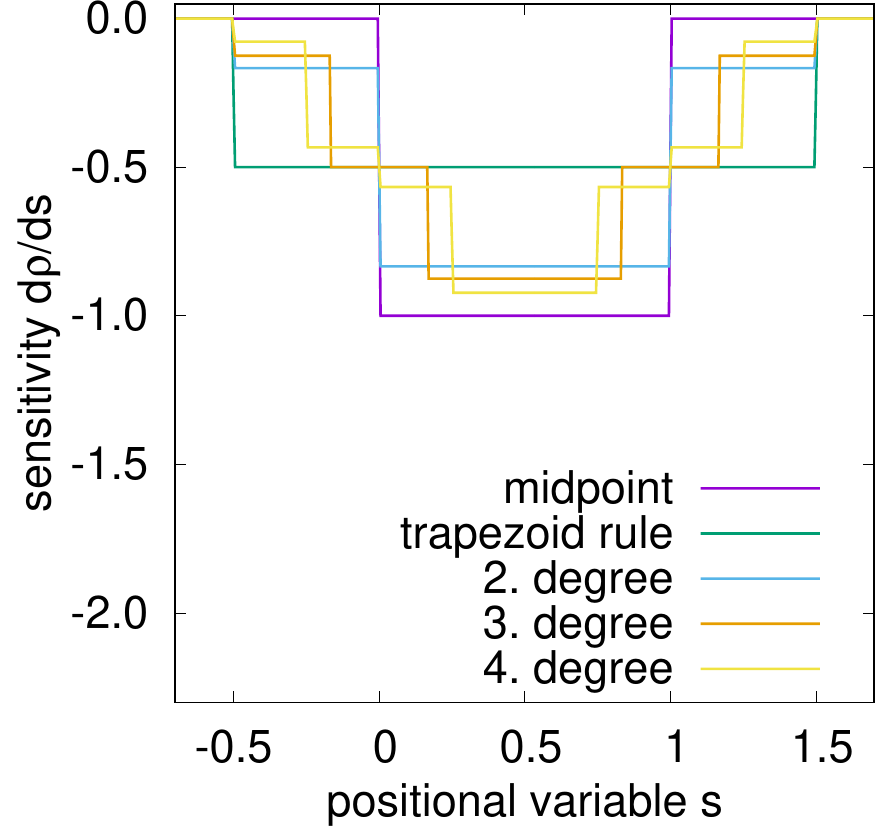}
  \\
  \includegraphics[width=0.33\textwidth]{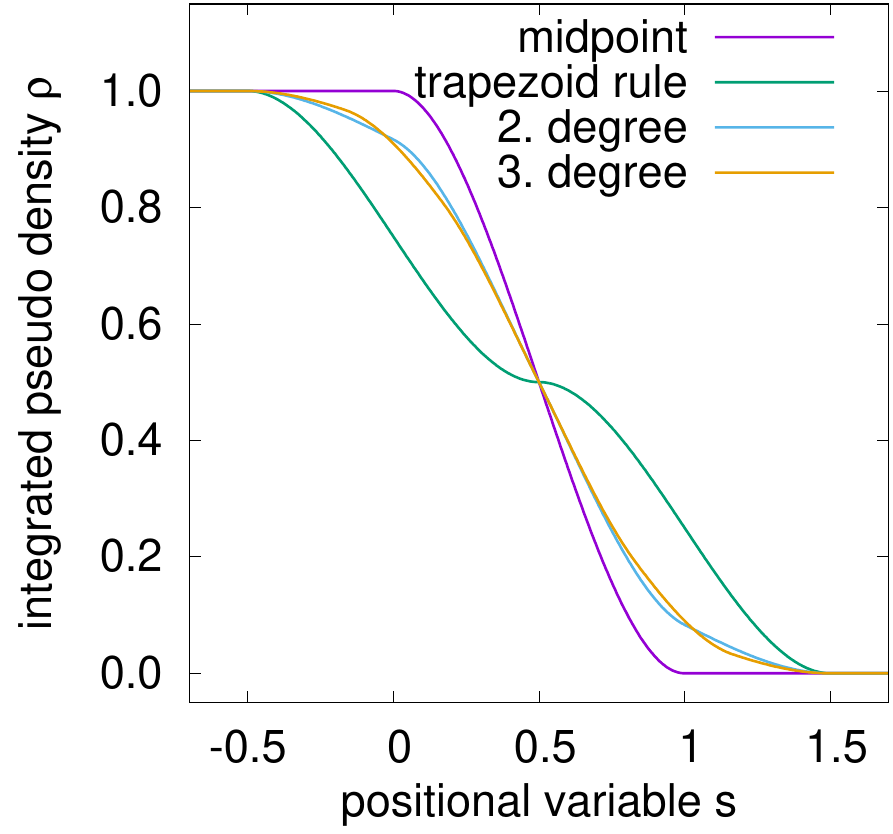}
  \includegraphics[width=0.33\textwidth]{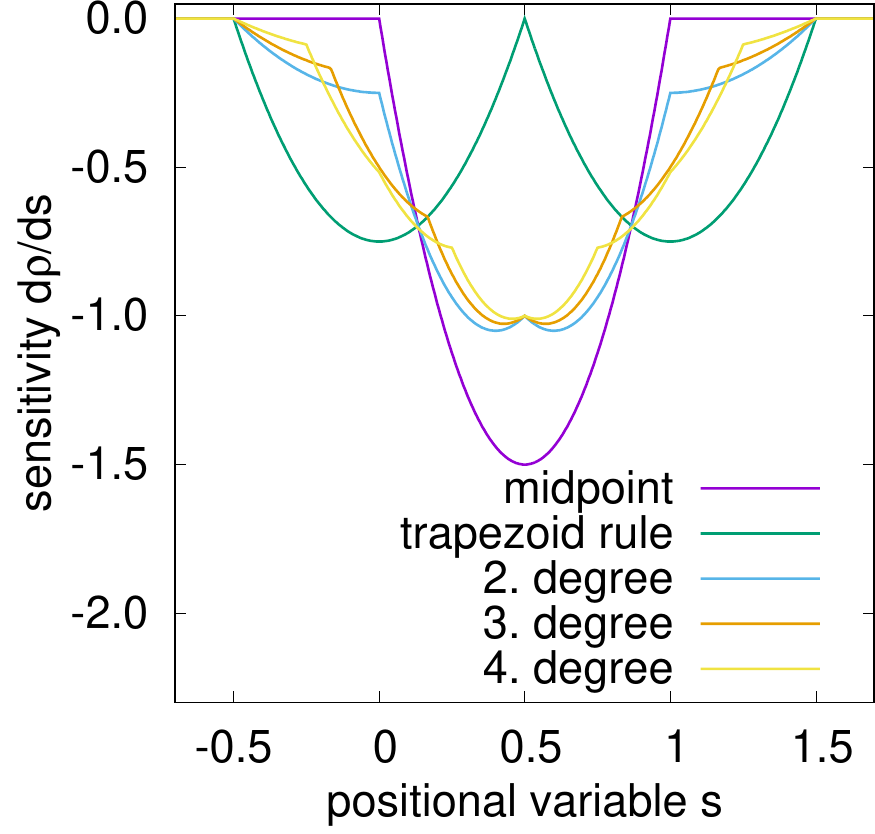}
  \\
  \includegraphics[width=0.33\textwidth]{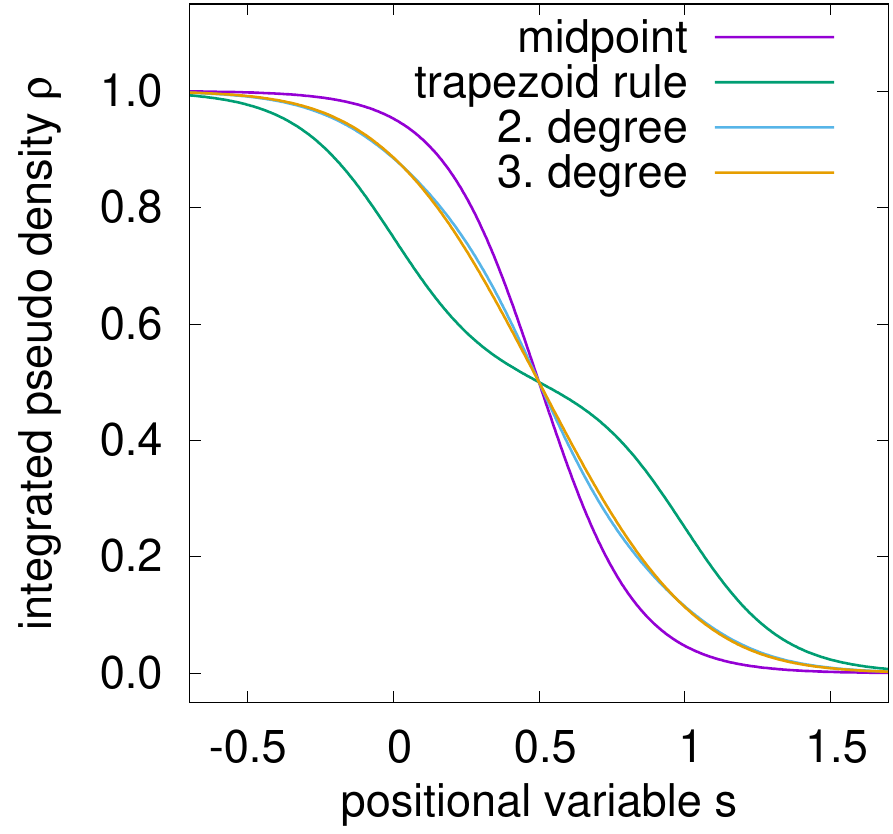} 
  \includegraphics[width=0.33\textwidth]{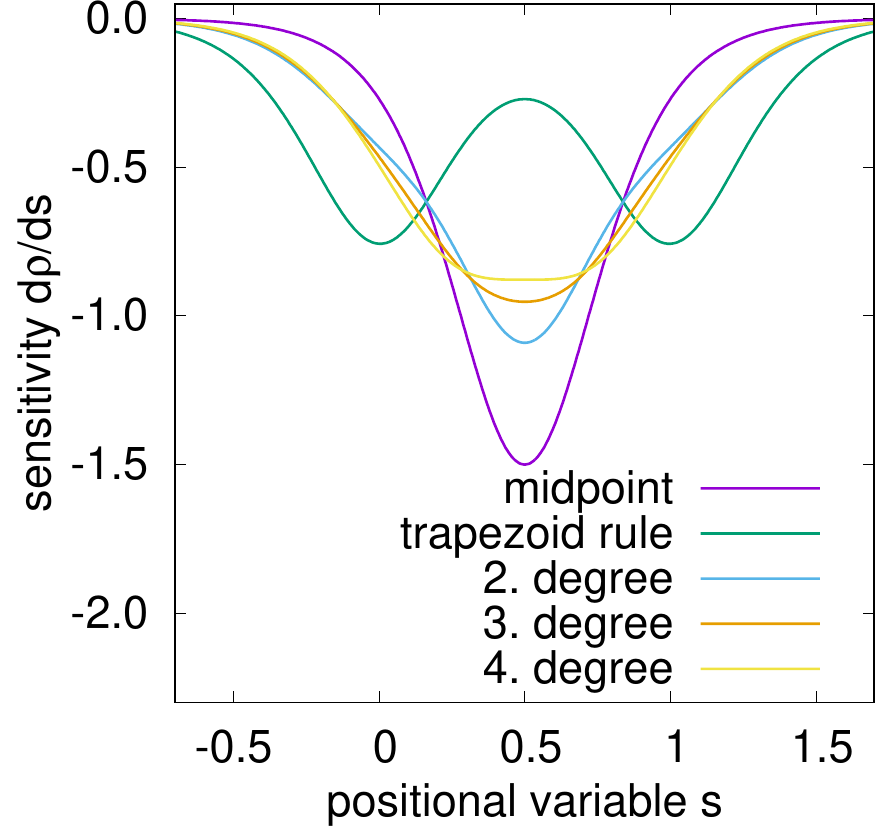} 
  \caption{For the first element $[0:1]$ we perform numerical integration \eqnref{eqn:rho_e_num_int} by closed Newton--Cotes in the left column and numerically evaluate the shape sensitivity in the right columns. The boundary modeling functions are linear (first row), polynomial (second row) and $\text{tanh}$ (last row).}
  \label{fig:num_int}
\end{figure}

When using numerical integration, the term $\sens{\rho_e}{s}$ in \eqnref{eqn:dJ_ds} is found by
\begin{equation}
   \sens{\rho_e}{s} = \sum_i^{N_\text{ip}} w_i \; \sens{\widetilde{H}(\Vec{x}_i)}{s}.
   \label{eqn:drho_ds_num}
\end{equation}
The test problem is now used to investigate the effect of the number of sampling points in numerical integration on \eqnref{eqn:drho_ds_num}. We consider the integral of the first element from the left and vary the bar position, see the upper row in \figref{fig:num_int}. We also introduce the element transition zone, $w_\rho$, which is defined as the number of elements across the boundary with intermediate density ($\rho_{\min} < \rho_e < 1$) multiplied by the element edge length, $l_\text{el}$.

For the linear boundary modeling function (top row in \figref{fig:num_int}) with midpoint integration at $\vec{x}_0$, $\rho_e(\vec{x}_0)$ varies from 1 to 0 with $w_\rho = w$. With trapezoidal rule, averaging the boundary modeling function values at the left and right node of element 0, $\rho_e(\vec{x}_0)$  `sees'  the boundary a half element earlier and a half element longer, and the element transition zone is: $w_\rho = w + l_\text{el}$.

The polynomial function \eqnref{eqn:transition_poly}, shown in \figref{fig:cut_smooth}, has zero slope at the end of its transition zone. Positioning the shape in the center of the element, a variation of the position has low impact when the boundary modeling function is only sampled at the ends of the transition zone by trapezoidal rule---see the center row in \figref{fig:num_int}. The $\text{tanh}$-like function \eqnref{eqn:tanh} shows similar behavior, but less pronounced.

The transition zone parameter $w$ for the smoothing function and number of sampling points in numerical integration are correlated. Enlarging the transition zone allows for a lower degree of numerical integration. Generally the density transition zone is $w \leq w_\rho \leq w + l_\text{el}$.

\begin{figure}[ht!]
  \centering
  \includegraphics[width=0.3\textwidth]{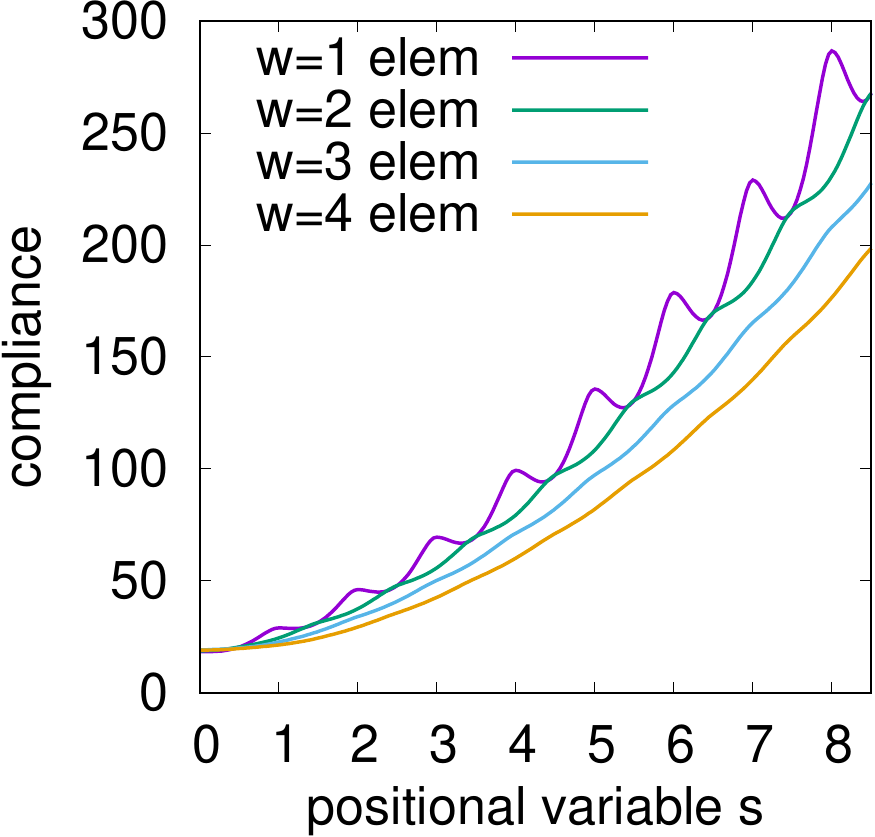}
  \includegraphics[width=0.3\textwidth]{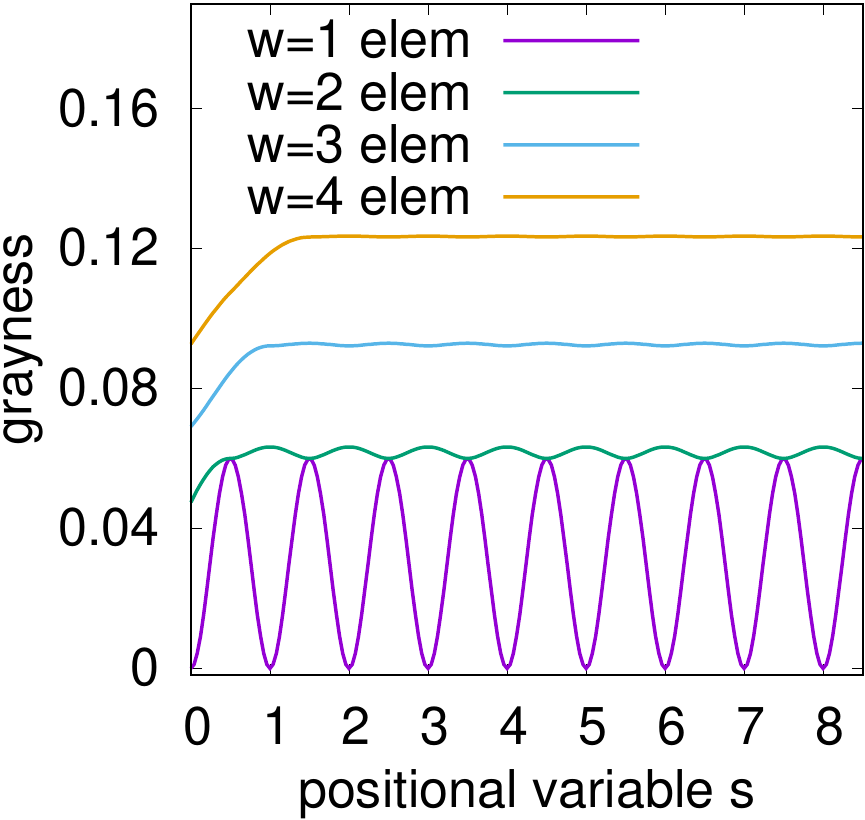}
  \includegraphics[width=0.3\textwidth]{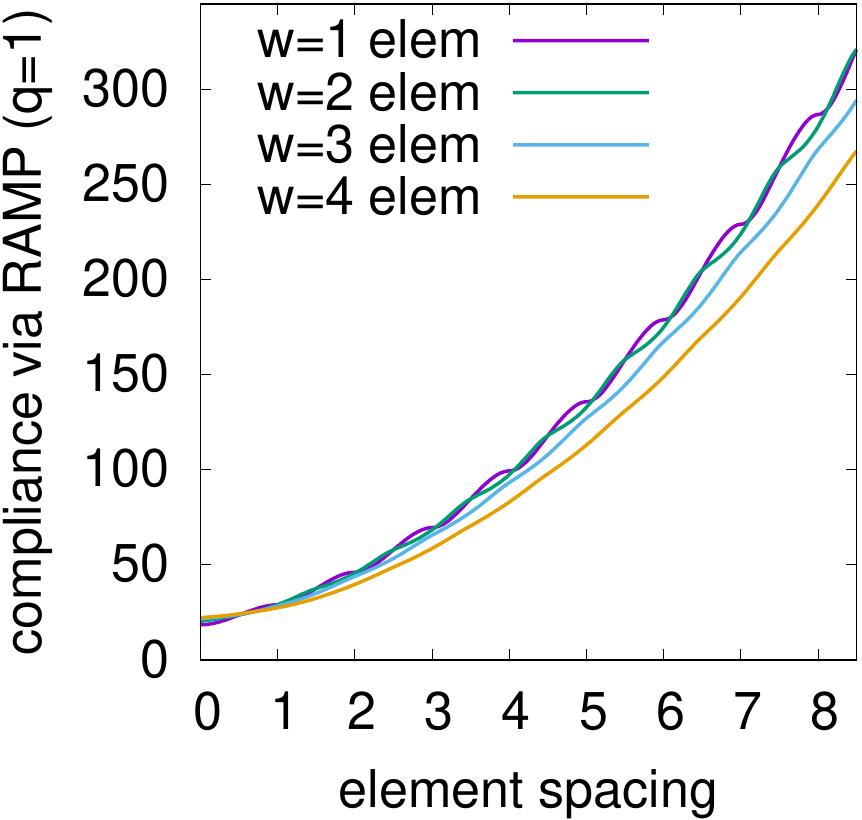} 
\caption{We vary the polynomial smoothing (center in \figref{fig:num_int_move}) by increasing the transition zone $w=2\,h$ from one to four elements, each time with midpoint integration. In the right figure RAMP interpolation \eqnref{eqn:ramp} is applied.}
  \label{fig:num_int_slope}
\end{figure}

In the following, we extend the transition zone $w$ from one element to up to four elements. We use polynomial smoothing with midpoint integration. The left figure in \figref{fig:num_int_slope} shows that already a doubled transition zone of two elements results in a significantly more monotonous compliance function over the parameter. However, due to the linear material interpolation, see \figref{fig:grayness}, the compliance value becomes artificially good due to a more blurred boundary. The grayness measured by \eqnref{eqn:grayness} is shown in the center image, revealing again that midpoint integration with transition zone of one element is not sufficient. Note that a wider grayness zone $w_\rho$ allows to include more information from the ersatz material sensitivity. Finally, we apply the RAMP material interpolation, discussed in \secref{sec:benchmark}; the right figure in \figref{fig:num_int_slope} shows that it helps compensate for the non-smoothness resulting from inaccurate numerical integration.

\subsubsection{Computing the signed distance}
\label{sec:signed_distance}
Feature-mapping methods often employ a signed-distance implicit function when using pseudo-density mapping, as this maintains the transition zone, as discussed in \secref{sec:boundary_smoothing}.

For some explicit geometry descriptions, the signed distance can be easily computed using an analytical expression. For example, the distance to the edge of a circular, or spherical feature can be directly computed from the feature parameters (center coordinates and radius). Features described by offset surfaces \citep{norato2015geometry, zhang2016geometry}, whose boundary is defined as the set of all points equidistant to a medial line segment or surface, readily provide direct expressions for the signed distance in terms of the design parameters. Also, for a design feature aligned to one principal axis, the signed distance is easily computed. We use this approach in the test problem in \secref{sec:benchmark}, where the design variable is simply the position of the left edge of the vertical feature. \cite{wein2018combined} use a simplified spline model of horizontal or vertical features, where the distance is obtained in the same way, see \figref{fig:spline}~(b), although the distance in their model is not necessary the exact distance to the feature.

To obtain the signed-distance for more complex explicit geometry descriptions, a popular method is to compute an equivalent implicit function, which is also a signed distance (i.e. a signed-distance level-set function). This can be achieved using schemes popular with level-set methods, such as the fast-marching method \citep{adalsteinsson1999fast}, or iteratively solving a Hamilton–-Jacobi equation in pseudo-time $t$
\begin{equation}
   \frac{\partial d(\vec{x})}{\partial t} + \text{sign}(\phi(\vec{x}))\left(||\nabla d(\vec{x})||-1 \right)=0.
   \label{eqn:signed.pde}
\end{equation}
The alternative is to directly compute the shortest distance from a point to the boundary of the feature, as done by \cite{norato2018topology} when using explicit features defined by supershapes.

These methods can also be used to compute the signed distance for other implicit geometry descriptions. However, these methods can be computationally expensive, especially if required each time the design changes. Thus, for implicit geometry representations, \cite{zhou2016feature} propose using a first-order Taylor approximation of the signed-distance function in the form of
\begin{equation}
   d(\vec{x}) \approx \frac{\phi(\vec{x})}{||\nabla \phi(\vec{x})||}.
   \label{eqn:signed.approx}
\end{equation}


\subsection{XFEM approaches}
\label{sec:xfem}
The alternative to pseudo-density mapping is to use an immersed boundary method.
The main challenge of mapping geometry onto a fixed grid is that the boundary does not align with the fixed-grid elements. Immersed boundary methods resolve this by introducing extra terms that model discontinuities within elements, while preserving the sharpness of the geometric interfaces. The eXtended Finite Element Method (XFEM) is a popular immersed boundary approach that has been utilized by several feature-mapping methods.

XFEM approaches model discontinuities by adding enrichment functions and additional degrees of freedom to nodes around the discontinuity. It was originally developed to model crack propagation without re-meshing  \citep{belytschko1999elastic,moes1999finite}. XFEM can also model discontinuities between different materials, or material and void, within an element. Thus, XFEM can be used to model the material discontinuities created by mapping features onto a fixed grid. The literature on XFEM is vast \citep{belytschko2009review,yazid2009state} and an in-depth review is not the focus of this paper. Instead we focus on relevant methods and issues encountered when using XFEM for feature-mapping methods.

There are two types of discontinuity that are considered in feature-mapping methods: material-void (or a strong discontinuity) and material-material (or a weak discontinuity). In general, three components are required to implement an XFEM scheme for material discontinuity: 1) enrichment strategy, 2) interface conditions, and 3) numerical integration.


\subsubsection{The simple scheme}
\label{sec:xfem_simple}
For the strong discontinuity case of a material-void interface, a simple scheme may be used, whereby a Heaviside enrichment is applied to the primary field (e.g. displacement or temperature) within the element
\begin{equation}
    \label{equn:HS_enrich}
    \vec{u}(\vec{x}) = \sum_{i=1}^n H(\phi(\vec{x})) N_i(\vec{x}) \vec{u}_i ,
\end{equation}
where $\vec{u}(\vec{x})$ is the physical field at point $\vec{x}$ within an element with $n$ nodes and shape functions $N_i(\vec{x})$, and $\phi(\vec{x})$ and $H(\vec{x})$ are implicit and Heaviside functions, as defined in \secref{sec:density_mapping}.

In this scenario, if the boundary is traction-free then there are no interface conditions and no additional degrees of freedom are required \citep{villanueva2014density}. This leads to a simple scheme, where element matrices are computed by numerical integration over the material domain. This is usually achieved by automatically sub-dividing the material domain into triangular sub-cells (e.g. using Delaunay triangulation) and using quadrature rules over each sub-cell. However, integration schemes without quadrature sub-cells have also been used \citep{li2012xfem}.


\subsubsection{Numerical aspects}
\label{sec:xfem_num}
The simple scheme for strong discontinuities has been utilized in several feature-mapping methods \citep{li2012xfem,zhou2013engineering,liu2014level}, its appeal being simplicity of implementation and ability to capture the sharp interface at the material-void boundary. However, there are several issues, or pitfalls, that can be encountered when using the simple scheme. These issues are discussed in the following along with potential solutions from the literature. In addition, the weak material-material discontinuity requires a more complex treatment.

During topology optimization, situations could occur where the design contains a material ``island", completely surrounded by void material and disconnected from the main structure. This causes the global system matrix to become singular, leading to numerical problems in solving the discretized governing equations. This can occur in feature-mapping methods if a solid component is mapped onto the fixed-grid, but does not overlap any other part of the solid region. A common remedy is to fill the void region with a fictitious weak material, which has properties several orders of magnitude lower than the real solid material \citep{wei2010study}. If the fictitious material is sufficiently weak, then the simple Heaviside enrichment scheme can still be used, as the error in ignoring the interface condition is small \citep{wei2010study}. An alternative was proposed by \cite{makhija2014numerical}, where each node is attached to a fictitious point in space by a soft spring. The advantages of this approach are that elements completely in the void phase are not assembled into the global matrix, reducing computational effort, and it avoids spurious load transfer through the void regions.

\begin{figure}[ht!]
  \centering
  \includegraphics[width=0.45\textwidth]{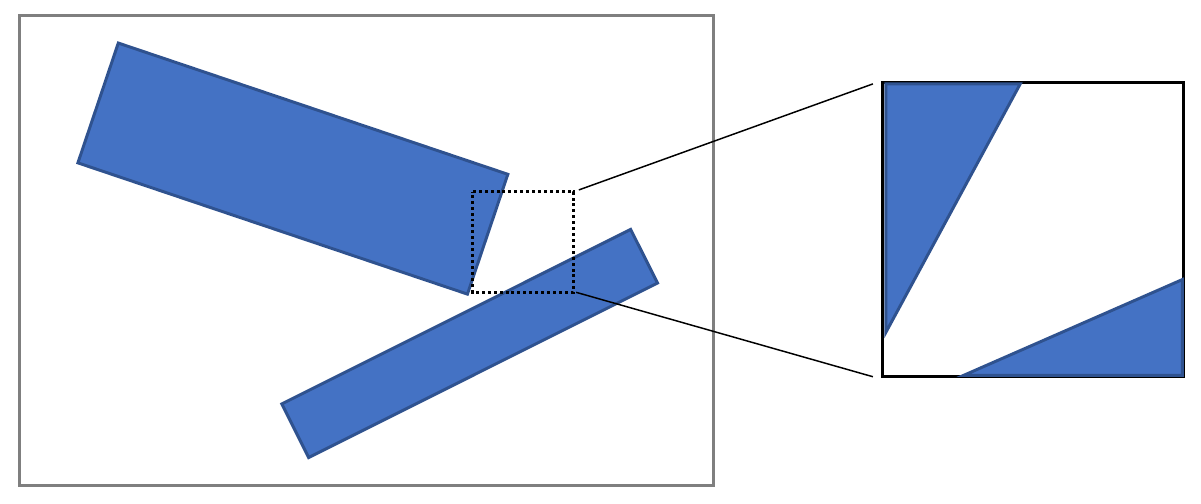}
  \caption{Potential non-physical coupling in XFEM.}
  \label{fig:nonphys}
\end{figure}

The simple scheme is only valid if the smallest geometric detail is larger than 2 elements \citep{villanueva2014density}. However, situations may occur during optimization when this is not true, potentially leading to interpolation error of the geometry and non-physical coupling between disconnected material phases (when the width of void feature is smaller than an element---see \figref{fig:nonphys}). This issue was demonstrated by \citet{makhija2014numerical} using a 1D bar example, where a non-zero reaction force was obtained when a gap in the bar was less than the element edge length. To address this issue, \cite{makhija2014numerical} proposed a generalized Heaviside enrichment strategy, based on the work of \cite{hansbo2004finite} and \cite{terada2003finite}, which captures the physical discontinuity by adding enrichment functions and additional degrees of freedom, depending on the order of discontinuity around a node.

It should be noted that the issue of non-physical coupling within an element is more common in free-form topology optimization methods, compared with feature-mapping methods, as these methods have some high-level control of the geometry. However, non-physical coupling can still occur in feature-mapping methods if solid components are close, such that the gap between them is less than the size of an element, as shown in \figref{fig:nonphys}. This type of non-physical coupling can also occur in pseudo-density mapping methods. However, a similar treatment has not been developed, possibly because pseudo-density methods do not aim to create a sharp interface in the analysis and therefore, this numerical artifact is not seen as an issue.

Modeling the weak discontinuity created by a material-material interface requires a more advanced XFEM scheme. Several authors have used an enrichment function proposed by \cite{moes2003computational}, which is a $C^0$-continuous enrichment function that inherently satisfies continuity in the primal solution (e.g. displacement, temperature) at an interface with a weak discontinuity in material properties. Alternatively, the generalized Heaviside enrichment strategy proposed by \cite{makhija2014numerical} can be used to model both strong and weak discontinuities. However, an additional constraint is required to enforce continuity across the interface for a weak discontinuity. This can be achieved using a scheme such as the stabilized Lagrange multiplier method, or Nitsche's method.

A further issue that affects both simple and more advanced XFEM schemes is the ill-conditioning of global system matrices due to very small integration regions, compared to the element size. This leads to convergence issues for nonlinear problems and iterative linear solvers \citep{lang2014simple}. A standard solution is to use a preconditioner to improve the condition number. For example, \cite{lang2014simple} introduced a simple and efficient geometric preconditioner for the generalized Heaviside enrichment scheme. It only requires knowledge of nodal basis functions and the interface geometry, so it can be computed before assembling the system matrix. This method proved effective at reducing the condition number, while maintaining accuracy.

Small integration regions can also affect the accuracy of results at the interface. For example, \cite{van2007stress} showed large errors in stress when the integration region of an element on the solid-void boundary was small. They discussed several possible remedies, including: removing elements with small solid parts, moving the closest mesh node or moving the boundary to eliminate the small area, post-processing to remove stresses in elements with small solid areas, or computing stresses using a smoothing scheme.

Finally, \cite{sharma2017shape} showed that the generalized Heaviside enrichment strategy proposed by \cite{makhija2014numerical} produces a smooth, non-oscillatory response function as the design changes. Thus, XFEM methods do not appear to produce a non-smooth response function when modeling the boundary using a sharp step function, in contrast to some of the pseudo-density material interpolation schemes (as shown in \secref{sec:density_note}). However, it was also shown that the shape sensitivity for the XFEM scheme can be oscillatory and that oscillations decreased with mesh refinement. Thus, it was concluded that the oscillations were mainly caused by the accuracy of mapping the geometry to the fixed-grid elements for numerical integration. Note that pseudo-density schemes can also produce a smooth, non-oscillatory response function, if implemented correctly (see discussion in \secref{sec:density_mapping}).


\subsubsection{Sensitivity analysis}
\label{sec:xfem_sens}
Sensitivity analysis for XFEM is generally more difficult compared with the pseudo-density approach. This is mainly because XFEM uses a more complex procedure to compute element matrices that involves numerical integration over material sub-domains. In the literature there are three main approaches to computing sensitivities when using XFEM in optimization.

The first approach is to differentiate, then discretize. This avoids computing derivatives of the change in integration regions as the interface moves, as sensitivities are derived from the continuum equations. A common example of this approach is to use shape sensitivities \citep{zhou2013engineering,liu2014level,wang2014topological}. However, it is well-known that convergence issues may occur due to the discretization error and often some form of regularization or smoothing is required \citep{VanDijk:2013:Level,liu2014level}.

The second approach is to discretize, then differentiate, but with a semi-analytical approach. The idea is to compute the derivative of the element matrices with respect to the design variables using the finite difference method. This derivative term is then inserted into the analytically derived sensitivity formula. Thus, sensitivities are consistent with the numerical discretization, but the semi-analytical approach avoids explicitly computing derivatives with respect to changes in the integration sub-domains. The finite difference approach is reasonably efficient, as it is only performed for elements that contain an interface and does not require assembling and solving a system of equations. This approach has proved effective and has been used in several feature-mapping methods \citep{van2007stress,sharma2017shape}.

However, the finite difference scheme should ensure that the design variable perturbation does not cause a change in element status, e.g. an element containing a material-void interface does not become either fully void, or fully solid. This causes problems in the derivative computation as it changes the number of degrees of freedom \citep{zhang2012integrated,noel2016analytical}. Several methods have been proposed to avoid this issue. One method is to  perform both forward and backward finite differences and check if either cause a status change. If neither cause a change, then the central difference is used. However, if the forward (or backward) difference causes a status change then only the backward (or forward) difference is used. Another method is to perturb the interface such that the finite difference perturbation cannot cause a status change \citep{sharma2017shape}. Alternatively, the finite difference perturbation step can be reduced to a magnitude that avoids a status change, although if the step magnitude is too small, numerical round-off errors can occur.

The third approach is to discretize first and then differentiate using a full analytical approach, without finite differencing.. The challenge is to compute the analytical derivative for the change in the integration sub-domains with respect to the design variables. \cite{zhang2012integrated} developed an analytical derivative for a material-material interface, when the geometry is represented by nodal implicit function values. \cite{noel2016analytical} and \cite{najafi2015gradient} proposed schemes utilizing a velocity field to efficiently compute the analytical derivatives. These fully analytical schemes are more complex and difficult to implement than the semi-analytical scheme, but are more efficient and avoid the status change issue.


\section{Combination of features}
\label{sec:combination}

The foregoing section describes the approaches that existing techniques use to map individual geometric features onto the fixed analysis mesh. To be able to modify the topology of the structure, it is also necessary to combine these features.  This is one of the key ingredients of performing topology optimization with high-level geometric features, and has received considerable attention in recent years. In this section we specifically consider the combination of closed regular sets (solids or holes). Unless otherwise stated and for brevity, whenever we refer to combination of solids we also refer to combination of holes.

The combination of solids in all approaches corresponds in effect to Boolean operations between solids. Just like other aspects in this review, it is possible to categorize approaches that combine solids in different ways. The main criterion we use to categorize combination methods is whether the combination occurs before or after mapping to the fixed analysis mesh.

\subsection{Smooth combination functions}
\label{subsec:smooth_comb}
Many feature-mapping methods utilize smooth combination functions so that derivatives with respect to the high-level geometric parameters are continuous. For example, the non-differentiable Boolean union of multiple solids represented by implicit functions $f_i$  corresponds to their maximum, i.e., $\bigcup_i f_i = \max_i(f_i(\vec{x}))$ \citep{Shapiro:2002:Solid}. Common choices for differentiable smooth approximations $\widetilde{\max}_i(f_i(\vec{x}))$ include the well-known Kreisselmeier–-Steinhauser (K-S) function
\begin{equation}
\label{eq:KS_func}
\widetilde{\max_i}^{KS}(f_i) := \tfrac{1}{p} \ln{\sum_i^N \left(\exp{(p f_i)} \right)}
\end{equation}
and the $p$-norm
\begin{equation}
\label{eq:p-norm}
\widetilde{\max_i}^{p}(f_i) := \left(\sum_i^N{f_i^p} \right)^ {\tfrac{1}{p}},
\end{equation}
where $N$ is the number of values and $p$ is a parameter that controls the sharpness and accuracy of the approximation (a larger $p$ results in a more accurate estimate of the true maximum and a sharper function).

Another type of functions used to perform Boolean operations is R-functions. The R-conjunction corresponds to the logical AND, whereas the R-disjunction corresponds to logical OR. Compositions of these two fundamental R-functions can be used to construct any Boolean expression. There are several forms of R-functions; for exampl,e \cite{chen2007shape} use the following definition
\begin{equation}
\label{eqn:r-func}
\begin{aligned}
f_1 \cap f_2 := f_1 + f_2 - \sqrt{f_1^2 + f_2^2},\\
f_1 \cup f_2 := f_1 + f_2 + \sqrt{f_1^2 + f_2^2},
\end{aligned}
\end{equation}
which is differentiable everywhere except at $f_1 = f_2 = 0$. It can be seen that $f_1 \cap f_2$ is positive if and only if both $f_1$ and $f_2$ are positive. Whereas, $f_1 \cup f_2$ is positive if either $f_1$ or $f_2$ is positive. For example, assume $\phi_{1}$ and $\phi_{2}$ are implicit functions for two solid features and $\phi_v$ the implicit function of a void feature. The combined implicit function can then be defined as
\begin{equation}
\label{eqn:phi_comb}
\phi(\vec{x}) := (\phi_{1} (\vec{x}) \cup \phi_{2} (\vec{x})) \cap \phi_v (\vec{x}).
\end{equation}

\subsection{Combine-then-map approaches}
\label{subsec:combine-then-map}

Since a geometric representation of the solids is available that is independent of the analysis mesh, a natural approach is to combine the solids directly using their geometric representation, and then map the combined solid onto the analysis mesh, as described in \secref{sec:geometry_mapping}.

\subsubsection{Implicit geometric representations}
\label{subsubsec:combine-then-map-implicit}

As mentioned in the previous section, the Boolean union or intersection of features represented by implicit functions can be attained by computing their maximum or minimum, respectively. This is a strategy that has been used to combine both solids \citep{cheng2006feature, zhou2013engineering, guo2014doing, zhang2016new} and holes \citep{cheng2006feature, chen2007shape, mei2008feature, wang2012high, zhang2017explicit}.

\begin{figure}[ht!]
 \centering
  \includegraphics[width=.7\columnwidth]{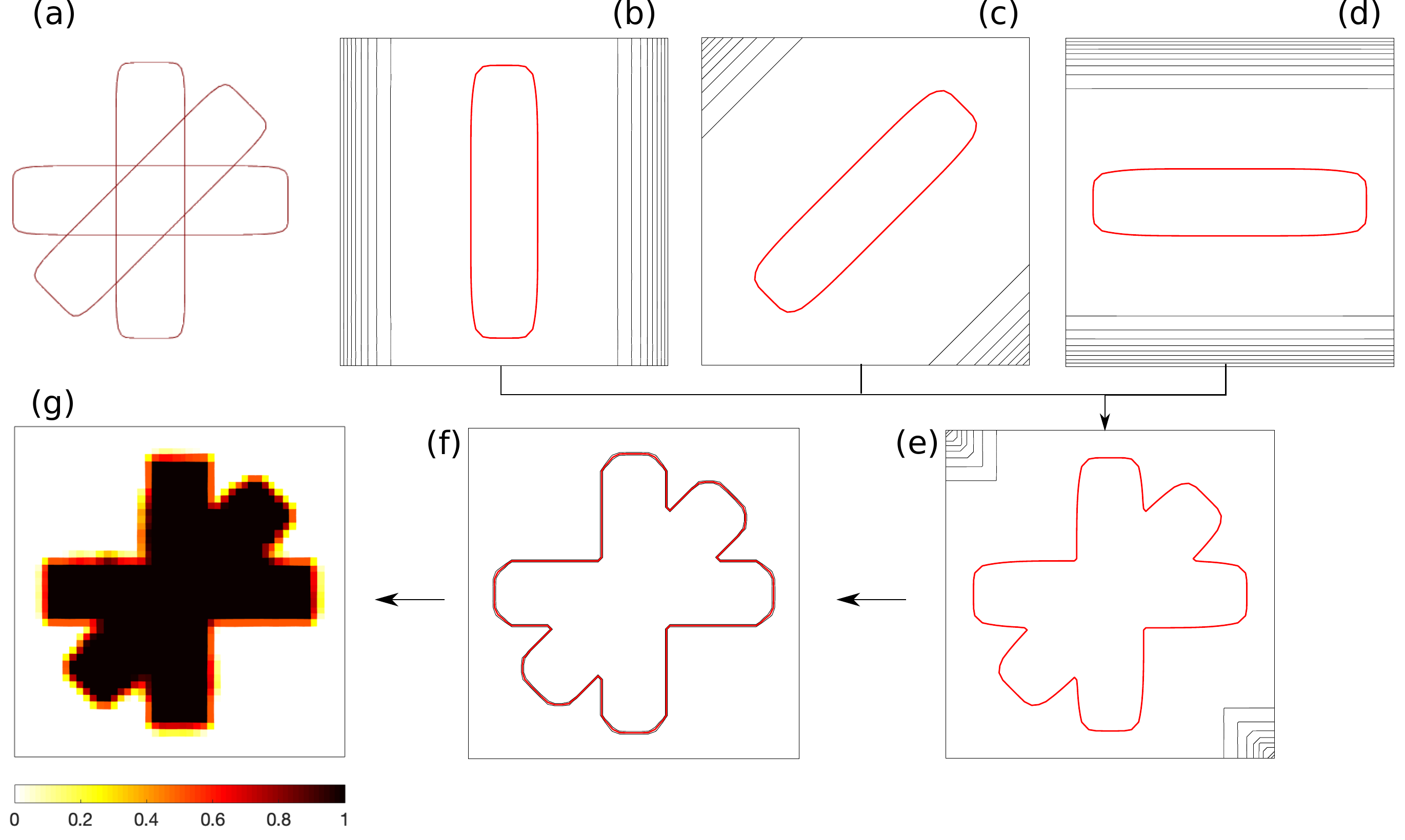} 
 \caption{\label{fig:mmc-combination} Combination by union of implicit functions-then-mapping: (a) rectangular bars modeled with hyperellipses; contour plots of (b)--(d)   implicit functions $\phi_i$, (e) union of implicit functions ($\Phi$) via maximum function, and (f) smoothed Heaviside of $\Phi$; and (g) pseudo-density used for analysis. Thick red line corresponds to zero level set for implicit functions, 0.5 level set for Heaviside (here and in following figures).}
\end{figure}

This combination approach is illustrated in \figref{fig:mmc-combination}, where the three rectangular bars are modeled using hyperellipse implicit functions (as used in several combine-then-map feature-mapping methods, e.g. \cite{guo2014doing,zhang2016new}).
 \figref{fig:mmc-combination}(b)--(d) show contour plots of $\phi_i$ for each of the three bars. All contour and fringe plots in this section are produced using a grid of $48 \times 48$ square elements. The Boolean union of the implicit functions for these three bars, as given by the true maximum function, is shown in \figref{fig:mmc-combination}(e). Note that smooth maximum functions could also be used.
 The combined implicit function is subsequently mapped onto the analysis mesh using a pseudo-density or immersed boundary approach, as discussed in \secref{sec:geometry_mapping}. A combined implicit function subject to a smooth Heaviside \eqnref{eq:smoothHS} is shown in \figref{fig:mmc-combination}(f) with the element constant pseudo-densities shown in \figref{fig:mmc-combination}(g) (which were obtained by the method used by \cite{zhang2016new}).

This approach is notable for the simplicity of the Boolean operations.  The simplicity is inherited from the implicit geometric representation, since in general it is much easier to perform Boolean operations with implicit, rather than with explicit geometric representations. Also, this combination approach is readily extended to 3-dimensional problems, e.g. \citep{liu20153D, zhang2016new}. 

Topological changes using this combination approach occur in one of three ways: (1) as solids move and overlap, the connectivity of the structure may change (holes may appear or disappear); (2) if a solid is engulfed inside another solid, it has no effect in the analysis due to the maximum operation, and thus the engulfed solid is effectively removed from the design; and (3) if one or more dimensions of a solid become sufficiently small, the effect of the solid on the analysis is negligible. 

\subsubsection{Explicit geometric representations}
\label{subsubsec:combine-then-map-explicit}

\begin{figure}[ht!]
 \centering
  \includegraphics[width=.7\columnwidth]{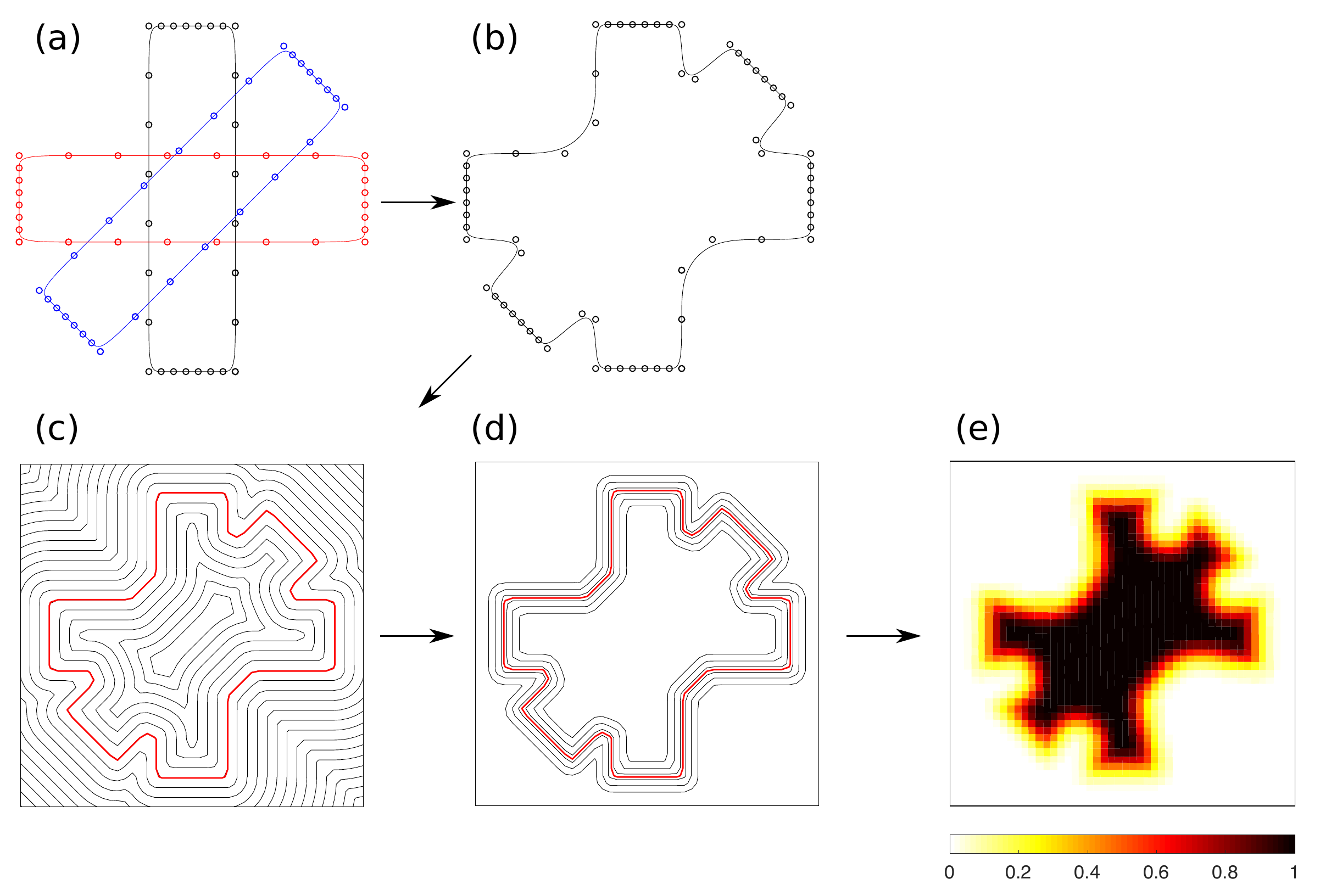} 
 \caption{\label{fig:bspline-point-deletion} Combination by union of explicit functions-then-mapping: (a) rectangular bars modeled with cubic B-splines, (b) union by deleting B-spline control points in the overlapping region, (c) signed-distance function of combined B-spline, (d) smoothed Heaviside of signed-distance function, and (e) pseudo-density used for analysis.}
\end{figure}

When the original geometric representation is explicit, two approaches have been employed to combine solids.  The first approach consists of performing the Boolean union directly on the explicit representation, and then converting the resulting design into an implicit geometric representation before mapping to the analysis mesh. This strategy is illustrated in \figref{fig:bspline-point-deletion}, where the three rectangles of the previous example are modeled using cubic B-splines with eight control points per side, see \figref{fig:bspline-point-deletion}(a). A combination technique used in this case \citep{lee2007smooth, seo2010isogeometric, zhang2017structural} consists of deleting from the current design those control points that lie in the overlapping region between the bars, so that the union of the primitives is given by a single B-spline made of the remaining control points, as shown in \figref{fig:bspline-point-deletion}(b). These works all consider B-spline-shaped holes; however, here we consider solid rectangles for consistency with the examples given for the other strategies.

After combining the solids, the explicit representation is transformed to an implicit representation, namely by computing the signed distance to the combined B-spline, as shown in \figref{fig:bspline-point-deletion}(c). An exact or smooth Heaviside approximation, such as the one presented in the preceding section, is then applied to the signed-distance function, at which point the mapping to the analysis can be completed in the different ways discussed in \secref{sec:geometry_mapping}, i.e. using pseudo-densities or an immersed boundary method.

This combination approach presents several challenges. In the particular case of B-splines, it is possible for the control points to be placed such that the B-spline can present self-intersections, which requires placing bounds on the positions of the control points \citep{lee2007smooth}, or employing special parameterizations of the B-spline \citep{zhang2017structural}. In both cases, it is necessary to determine the correct order of control points in the combined B-spline to generate the correct shape. Another, perhaps more pernicious challenge, is the potential lack of differentiability introduced by the control-point deletion approach. Suppose two solids overlap, and a control point of one of them lies exactly on the boundary of the other. Small positive and negative rotations of any of the primitives will cause that control point to be deleted or retained in the combined B-spline. Therefore, the combined B-spline may look appreciably different in both cases, which means the structural response will not be differentiable with respect to the orientation angle of either solid. To obtain an accurate union of the B-splines it is of course possible to introduce control points at the intersections with multiple knots to capture the sharp corners. However, this introduces additional challenges in the optimization, as the number of design variables (i.e., the positions of the control points) would increase with the additional control points. A third challenge lies in the computational cost incurred in translating the explicit representation to an implicit representation, e.g., the computation of the signed-distance field. Although there exist computational strategies to do this efficiently (see \secref{sec:signed_distance}), it still adds computational cost compared with directly using implicit geometric representations. Finally, while it is possible to perform Boolean operations of explicit representations of 3-dimensional solids, the aforementioned challenges are more difficult to solve for 3-dimensional problems. 

\begin{figure}[ht!]
 \centering
  \includegraphics[width=.7\textwidth]{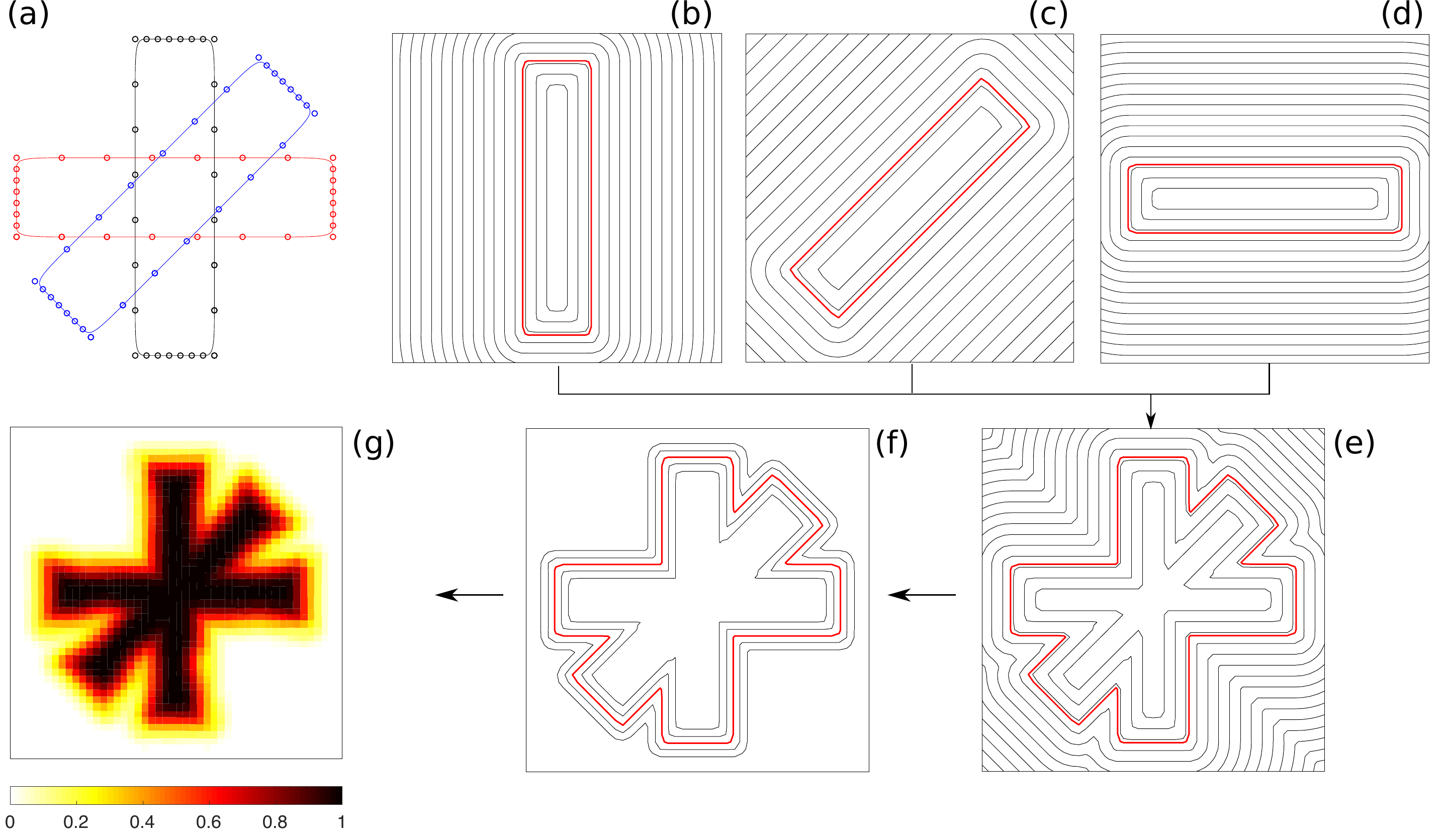} 
 \caption{\label{fig:bspline-implicit} Combination by conversion of explicit to implicit representation-then-mapping: (a) rectangular bars modeled with cubic B-splines, (b)--(d) signed-distance of individual bars, (e) union of implicit signed-distance functions via maximum function, (f) smoothed Heaviside of implicit union, and (e) pseudo-density used for analysis.}
\end{figure}

The second strategy to combine features with explicit representations is to first convert the explicit representation of each solid to an implicit representation, and then perform the combination of the individual implicit representations as in the previous section \citep{zhang2017topology}. This strategy is depicted in \figref{fig:bspline-implicit}, where each B-spline is first converted to an implicit representation (a signed-distance function), and then the combination of bars is achieved via the true maximum ---in the present case--- of the implicit functions.  This approach circumvents the problems arising from the deletion of control points and greatly facilitates the combination of primitives. However, there is additional computational cost, as a separate signed-distance function must be computed for each solid. This strategy is arguably similar to the map-then-combine approach described in the next section, because in this method the individual implicit functions are computed on the fixed grid prior to the combination. 

\subsection{Map-then-combine approaches}
\label{subsec:map-then-combine}

The alternative to the combination strategies described in the previous section, is to first map each individual solid to the analysis mesh and then combine the ensuing mapped variables, such as pseudo-densities, Heaviside function values, or even material property values. The combination could be done element-wise (e.g. element constant pseudo-densities), or at integration points (e.g. when using \eqnref{eqn:rho_e_num_int}).

All existing map-then-combine approaches require an implicit geometry description to achieve the mapping. Thus, map-then-combine approaches that start with an explicit geometry description first convert each feature to an implicit function, before mapping the individual implicit functions to the analysis grid. This is the same as the first step in \figref{fig:bspline-implicit}(a-d), where the geometry of each solid is described by a B-spline, which are then converted to implicit signed-distance functions. Therefore, the remainder of this section describes map-then-combine methods starting from an implicit geometry description.

\subsubsection{Property interpolation for hybrid approaches}
\label{subsubsec:map-then-combine-hybrid}

Some hybrid approaches described in \secref{sec:hybrid} combine a free pseudo-density field (as in density-based methods) with features using an extended material interpolation function, which interpolates between the solid-void pseudo-density field and solid features, e.g. \cite{qian2004optimal, wang2014topological}, or holes, e.g. \cite{kang2013integrated}. For example, if $E_s$ is the Young's modulus of the solid phase of the free pseudo-density field $\rho$, and $E_c$ modulus of the embedded solid component (assuming isotropic materials), an interpolation of the Young's modulus may be given by (cf., \cite{wang2014topological})
\begin{equation}
\label{eq:mod_Eint_solid}
    E(\vec{x}) = \widetilde{H}(\phi(\vec{x}))E_c + \left[1-\widetilde{H}(\phi(\vec{x}))\right] \left[ \rho^p(\vec{x}) \right]E_s,
\end{equation}
where $\phi$ is an implicit representation of the embedded solid primitive, $\widetilde{H}$ is a smooth approximation of the Heaviside function and $p$ is the SIMP penalization power. Note that in \eqnref{eq:mod_Eint_solid}, $\vec{x}$ is usually the element center, which is equivalent to mid-point integration (see \secref{sec:num_int}). Further approaches are discussed in \secref{subsec:combining-free-form-features}.

\subsubsection{Combining Heaviside functions}
\label{subsubsec:map-then-combine-HS}

Combination of multiple features can be achieved by combining Heaviside function values of each mapped feature using a maximum function. The single, combined Heaviside function can then be used to compute element pseudo-densities, or in immersed boundary methods (as discussed in \secref{sec:geometry_mapping}).
This approach was used by \cite{wein2018combined} to obtain element pseudo-densities by numerical integration via
\begin{equation}
\label{eq:comb_HS}
   \rho_e = \sum_j^{N_\text{ip}} w_j \; \widetilde{\max_i} \, \widetilde{H}_i(\Vec{x}),
\end{equation}
where the combination is done at the integration points. In \eqnref{eq:comb_HS} a smooth maximum and smooth Heaviside are used, thus making sensitivity analysis straight-forward.

It is interesting to note that in the case of midpoint integration and using \eqnref{eqn:rho_e_num_int}, \eqnref{eq:comb_HS} simplifies to
\begin{equation}
   \rho_e = \widetilde{\max_i} \, \rho_e^i,
\end{equation}
which is effectively the same as combining mapped element pseudo-density values (see \secref{subsubsec:map-then-combine-dens}). This highlights the close connection between the Heaviside function and pseudo-densities in feature-mapping.

\subsubsection{Combining pseudo-density values}
\label{subsubsec:map-then-combine-dens}

Another map-then-combine approach is to compute element pseudo-density values with respect to each feature, as described in \secref{sec:density_mapping}, and then perform the combination using a true or smooth maximum function. This approach also allows for an additional control of the combined features, by introducing variables that penalize the mapped densities of each solid feature separately. These are called size variables, which are penalized in the spirit of SIMP, so that a zero value indicates the solid has no effect on the analysis and thus can be removed from the design, whereas a value of unity indicates the solid must be retained \citep{norato2015geometry, zhang2016geometry}.    

\begin{figure}[ht!]
 \centering
  \includegraphics[width=.7\textwidth]{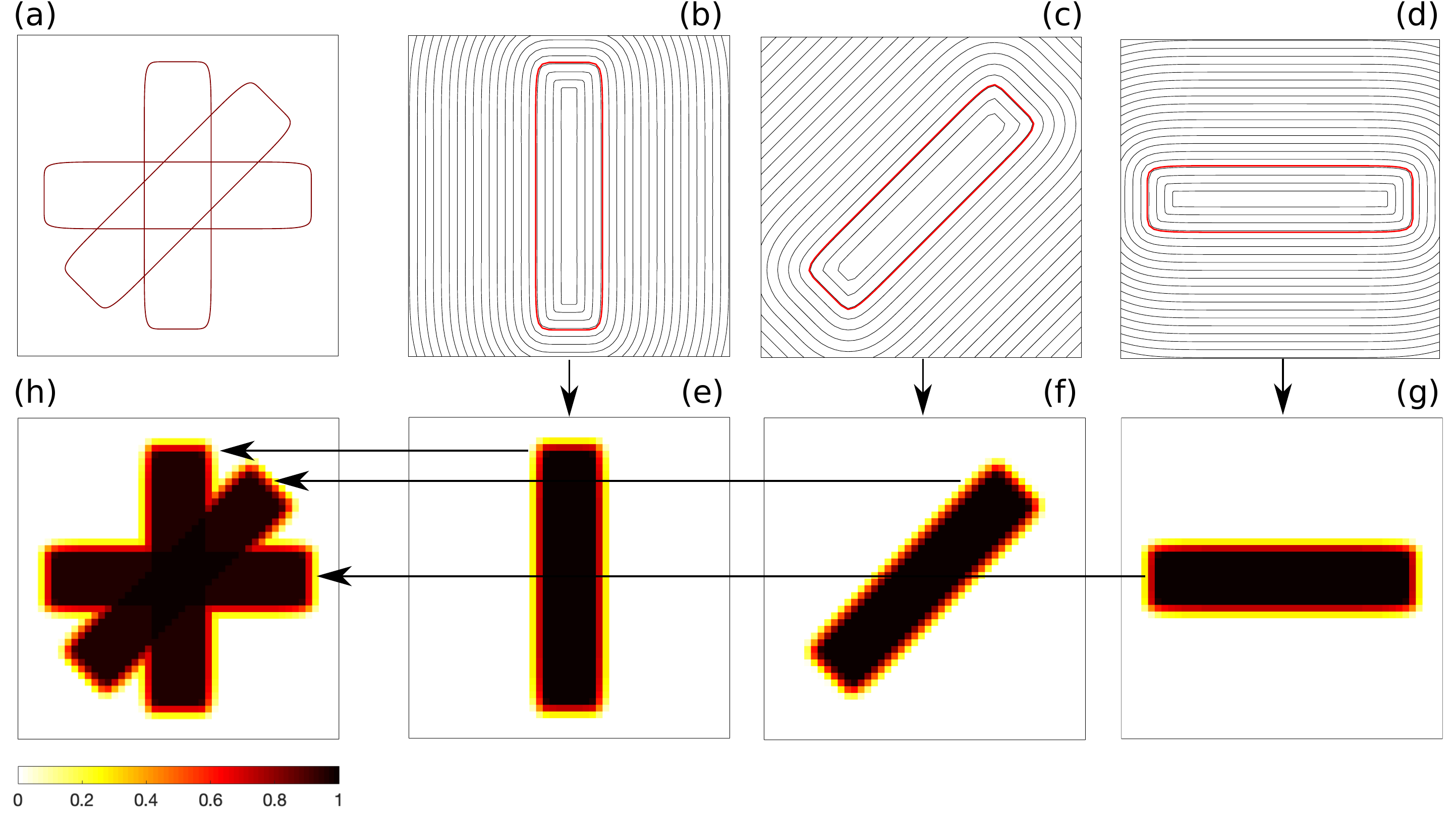} 
 \caption{\label{fig:gp-offset} Map-then-combine by union of offset surfaces: (a) bars modeled with offset surfaces (design variables are locations of endpoints of bars medial axes), (b)--(d) signed-distance field for each bar, computed directly from design parameters, (e)--(g) corresponding mapped pseudo-density for each bar, and (h) union of mapped densities using a smooth approximation of the maximum.}
\end{figure}

Without consideration for the aforementioned penalized size variable (that is, assuming all bars have a size variable of unity), the pseudo-density map-then-combine strategy is depicted in \figref{fig:gp-offset}.  \figref{fig:gp-offset}(a) shows the  three bars modeled as hyperellipses as before; \figref{fig:gp-offset}(b)--(d) are the signed-distance fields corresponding to these surfaces; \figref{fig:gp-offset}(e)--(g) show the mapped pseudo-densities for each bar, computed at each element of the mesh using \eqref{eq:geom_proj}; and \figref{fig:gp-offset}(h) shows the Boolean union of bars, obtained using a smooth approximation of the maximum function (here the $p$-norm \eqnref{eq:p-norm}).

\begin{figure}[ht!]
 \centering
  \includegraphics[width=.7\textwidth]{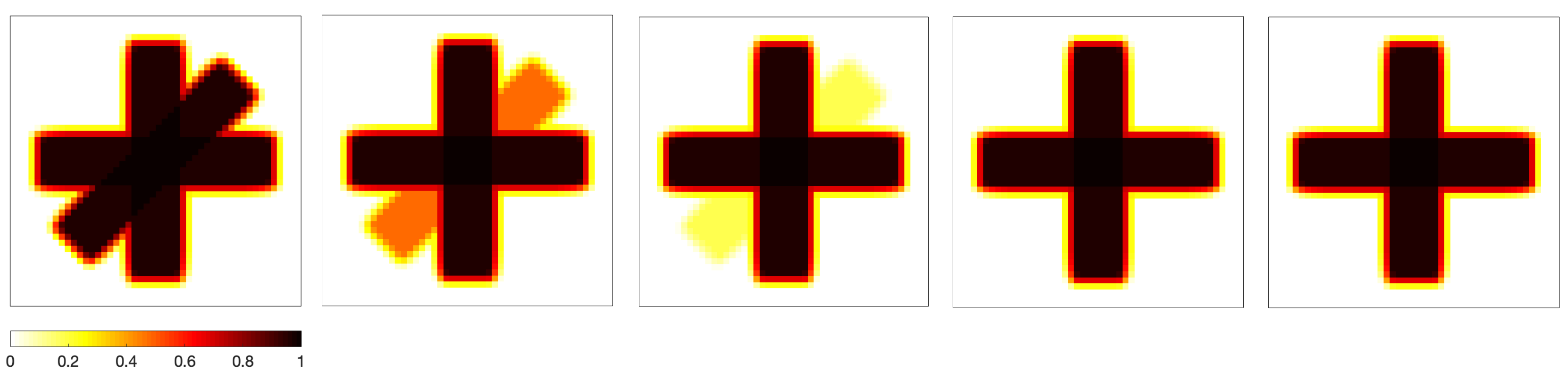} 
 \caption{\label{fig:gp-multiple-alphas} Three-bar example, where the diagonal bar has a size variable value $\alpha$ of (a) 1, (b) 0.8, (c) 0.6, (d) 0.1, and (e) 0. A power value of $p=3$ is used in \eqref{eq:effective-density} in all cases.}
\end{figure}

The combination of solids when a penalized size variable is used is as follows. First, element pseudo-densities are computed. The effective density at element $e$ for bar $i$ is subsequently computed as
\begin{equation}
\label{eq:effective-density}
    \widehat{\rho}_e^{i} = (\alpha^{i})^p \rho_e^i ,
\end{equation}
where $\alpha^{(i)}$ is the size variable corresponding to solid $i$ and $p$ is the penalization power. We note that if $\alpha^{(i)}=0$, then the effective density at element $e$ for bar $i$ is zero, hence this solid has no effect on the material properties at element $e$; this is true for every element for which $\rho_e \neq 0$, and so making the size variable zero effectively removes the solid from the design. The combination of the solids is subsequently obtained via, for example, a smooth maximum of the effective densities as
\begin{equation}
\label{eq:gp-composite-density}
    {\rho}_e = \widetilde{\max_i} \, \widehat{\rho}_e^{i}.
\end{equation}
\figref{fig:gp-multiple-alphas} shows the combined density after the union of all three bars when the diagonal bar has different values of its size variable $\alpha$ (while the other two bars have a size variable of unity). Clearly, as its size variable nears zero, the effect of the diagonal bar on the combined density vanishes. 

Topological changes using this approach can occur in different ways. As in the methods of \secref{subsec:combine-then-map}, when solids `move' in the optimization, holes can be created or disappear. We also note that some approaches (e.g., \cite{norato2015geometry,zhang2016geometry}) have the limitation that a solid cannot be removed by collapsing its dimensions, because for the sensitivities to be always well defined, it is necessary that the size of the sample window used to compute the volume fraction is smaller than the dimensions of the solid. The other removal mechanism, as aforementioned, is to make the size variable of the solid zero.

\subsection{Local minima}
\label{subsec:combination-local-minima}

\begin{figure}[ht!]
 \centering
  \subfloat[Setup]{   
  \includegraphics[width=.3375\columnwidth]{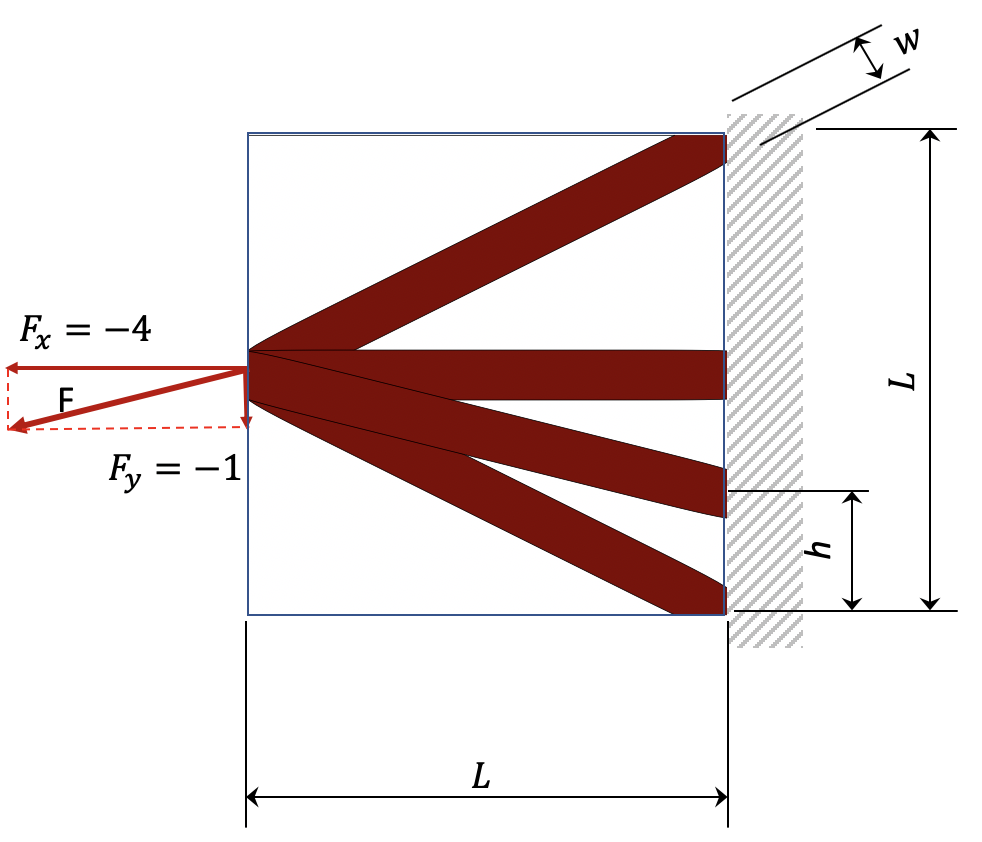}}  
  \subfloat[Compliance] {   
  \includegraphics[width=.4125\columnwidth]{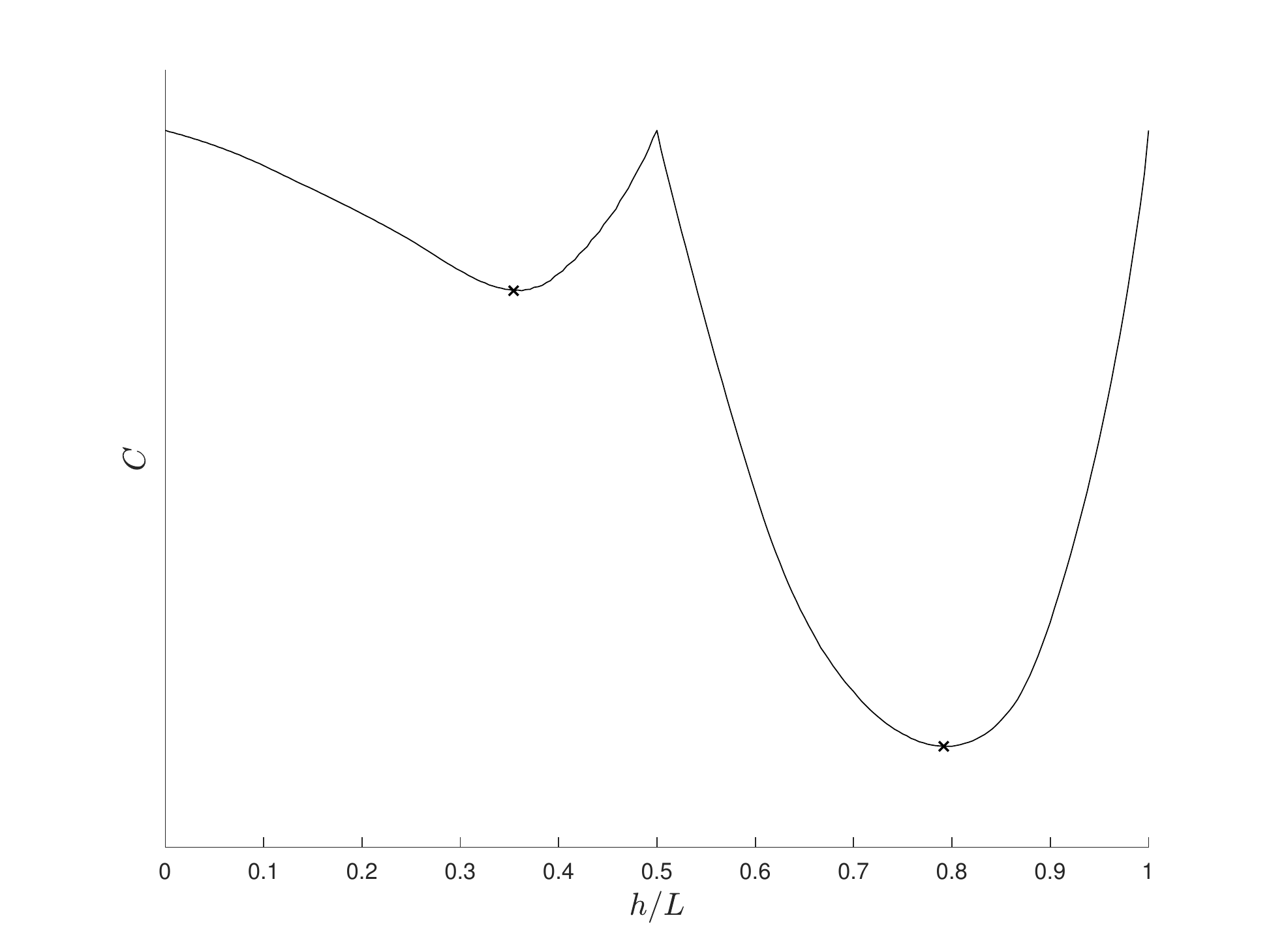}} 
 \caption{Effect of combination on local minima. The crosses mark the locations of the local minima.}
 \label{fig:combination-benchmark}
 \end{figure}

The combination of geometric features may lead to unfavorable local minima.  To illustrate this, we consider the example shown in \figref{fig:combination-benchmark}. Four bars are modelled with hyperellipses.  Three of the bars are fixed, and another one is moved by changing $h$. For $h \in \{0, L/2, L\}$ the moving bar entirely overlaps with one of the fixed bars. The design region is meshed with square bilinear elements with a relatively fine mesh. A binary pseudo-density mapping is used, where the element pseudo-density is either $\rhomin$ or $1$ depending on whether the element centroid is outside or inside of a bar, respectively. The combination is performed using a map-then-combine approach with a true maximum of the pseudo-density values.

\figref{fig:combination-benchmark} shows the compliance as a function of $h/L$. The actual magnitude of the compliance is not important; what is important is the presence of two distinct local minima, one of which ($h/L \approx 0.43$) is clearly worse than the other ($h/L \approx 0.79$). Therefore, if a gradient-based optimizer is used and the initial design has $h < L/2$, the optimizer will most likely converge to the poor local minimum. Thus, the more compact design representation used by feature-mapping methods (as opposed to the verbose representation used by density-based and level-set methods) is more prone to falling into unfavorable local minima depending on the initial design. Although all topology optimization techniques are dependent on the initial design (e.g., as shown by \cite{yan2018non} for density-based topology optimization of heat conduction structures), this dependency is more pronounced in feature-mapping techniques, as noted in \cite{norato2015geometry, zhang2016geometry}. We note that this has nothing to do with the particular feature-mapping technique used, but with the more restrictive geometric representation.

\section{Separation constraints}
\label{sec:separation}
Separation constraints are high-level geometric constraints that specify a minimum distance between solid components (or holes). When the minimum distance is set to zero, separation constraints are often called non-overlap constraints, as component overlap is prevented. These constraints can also be used to prevent components from leaving the design domain. Several techniques have been proposed to enforce this type of constraint in feature-mapping methods.

\begin{figure}[ht!]
  \centering
  \includegraphics[width=0.3\textwidth]{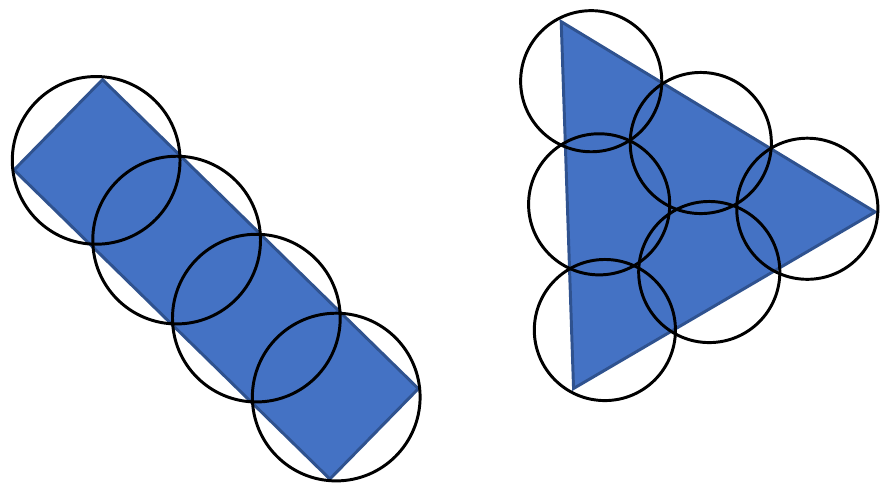}
  \caption{Finite circle method.}
  \label{fig:FCM}
\end{figure}

The earliest method is the finite circle method (FCM). The main idea is to approximate the shape of each component with a number of circles, \figref{fig:FCM}. Separation constraints can then be posed as simple geometric constraints on the minimum distance between circle centres. \cite{qian2004optimal} used a single circle for each component and \cite{zhang2006new} extended the idea to use multiple circles to approximate the shape of each component.

The main benefits of FCM are the simple definition of the constraints, which are continuous and differentiable. However, if there is a large number of components, then a large number of separation constraints is required, although most are usually inactive at the optimum. For example, even if only one circle is used to approximate each component, then $N(N-1)/2$ constraints are required for $N$ components. Also, component shapes are only approximated by circles, so separation constraints may not be able to reach their lower bound in some situations (due to circles covering a larger volume than the actual component). If high accuracy in the separation constraint is required, then more circles can be used for each component, which adds more constraints. A large number of constraints may affect the efficiency of the optimization \citep{zhang2011some}, although the number of constraints can be reduced by using different sized circles. \cite{xia2012superelement} also showed that to prevent components leaving the design domain using FCM, only a small circle at each corner of the component is required.

To use FCM more efficiently with a large number of components and constraints, \cite{gao2015improved} used constraint aggregation to combine all finite circle separation constraints into a single constraint function. An adaptive Kreisselmeier--Steinhauser (K-S) function is used, where the aggregation parameter is automatically determined to ensure the accuracy of the aggregation function. Another approach was proposed by \cite{zhu2017integrated} where the finite circle separation constraints are added to the objective using a combination of exterior and interior penalty methods.

A further limitation of the standard FCM is that it does not automatically adapt to components that are changing in size or shape. This is a challenging problem if the constraints are to remain continuous and differentiable throughout the optimization. However, \cite{zhang2012integrated} showed how this can be achieved for elliptically shaped components by linking the location and radius of each circle to the parameters of the elliptical shape.

\begin{figure}[ht!]
  \centering
  \includegraphics[width=0.32\textwidth]{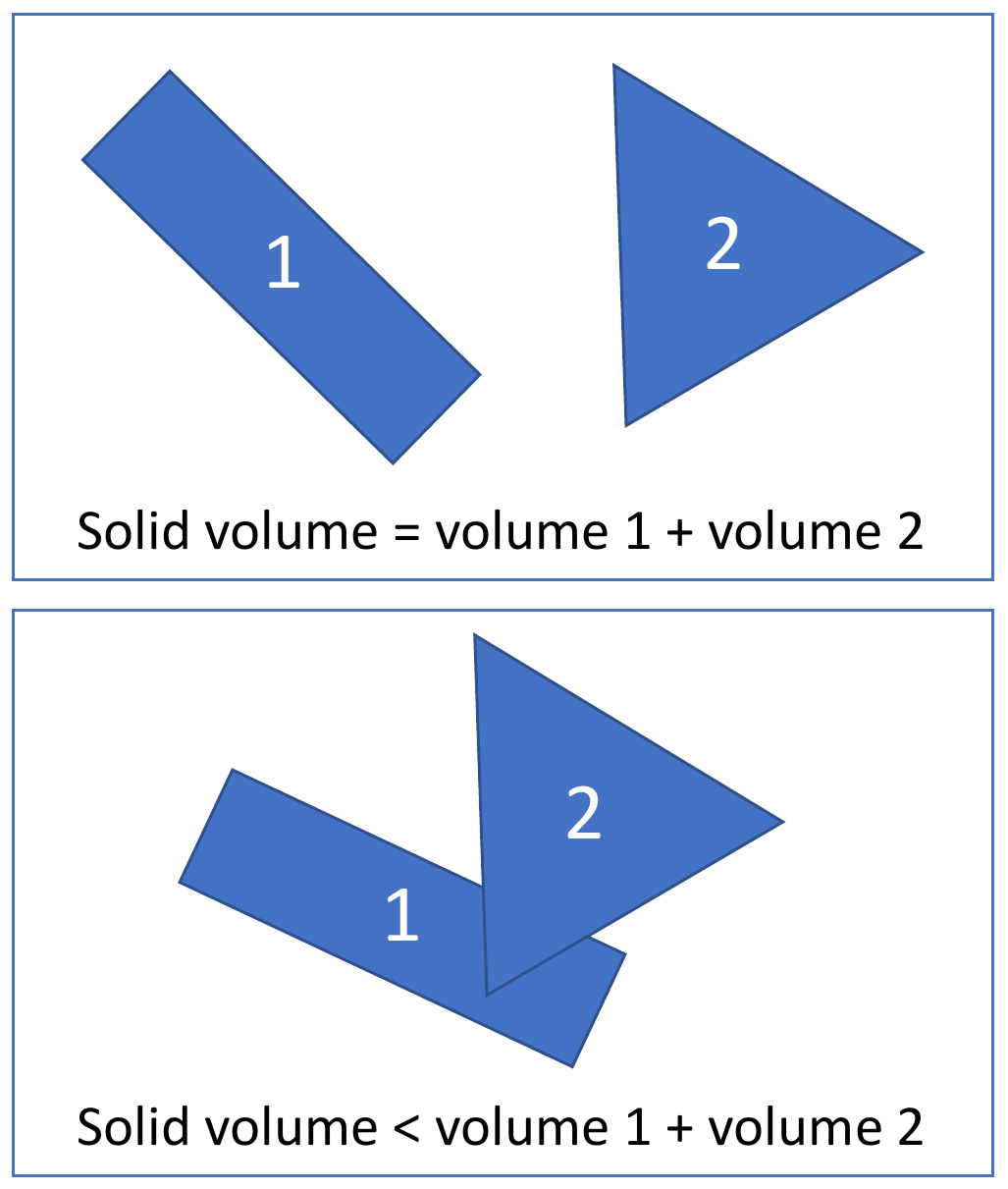}
  \caption{Integral method for preventing component overlap.}
  \label{fig:int_over}
\end{figure}

\cite{kang2013integrated} introduced an alternative to the FCM, suitable for feature-mapping methods where solid and/or void features are represented by an implicit function. The idea is to compare the integral of the solid region represented by the combination of implicit functions with the known volume of solid components. (The same idea also applies to void components). If the integral is less than the known volume, then there must be some overlap of components, or part of a component has left the design domain, \figref{fig:int_over}. This observation is used to formulate a single, differentiable constraint that prevents overlap for arbitrary shaped components and also prevents components leaving the design domain.

The integral method was extended by \cite{kang2016structural} for minimum distance separation constraints. To achieve this, each component is represented by a signed-distance implicit function. The signed-distance information is used to construct  ``virtual" components, whose boundaries are offset by half the minimum distance constraint value. The original separation constraint formulation is then applied to the ``virtual" components to ensure a minimum distance between components.

\cite{zhang2015explicit} introduced an approach using the structural skeleton, which is also suitable for methods where features are represented by implicit functions. The structural skeleton is defined as the set points that have at least two closest boundary points (this is also called the medial axis). To enforce a minimum distance separation constraint between two components, a signed-distance function is constructed for the combined implicit function of both components. A skeleton is then constructed to identify all points that are equidistant from both components. Finally, an explicit constraint is imposed on the minimum signed-distance value of all points belonging to the skeleton. This method can also be used to set maximum separation constraints. However, it requires construction of the signed-distance functions in each iteration and the number of constraints is linked to the number of components. Furthermore, the formulation is not differentiable and derivatives are approximated using finite differences and by assuming that the skeleton does not change when components move.

For methods based on pseudo-densities, \cite{zhang2017optimal} proposed a simple method using a map-then-combine approach. First the pseudo-densities for each component are mapped onto the fixed-grid. These are then combined using simple summation to create an auxiliary pseudo-density field. If any density value in the auxiliary field is greater than 1, then there must be some component overlap. Thus, an aggregated constraint function can be defined that ensures the maximum value of the auxiliary pseudo-density field is unity and hence prevent component overlap. This idea can easily be extended to provide a minimum separation constraint, by uniformly enlarging the size of the components by half the minimum separation distance before computing the auxiliary pseudo-density field.
 
The FCM- and integral-type methods for preventing component overlap can also be used to prevent components leaving the design domain. The integral approach achieves this without any modification to the original method, as it automatically detects when any part of a component has left the design domain. The FCM approach needs additional distance constraints to prevent components leaving the domain. This is straightforward for convex design domains \citep{zhu2008simultaneous}, but non-convex domains present difficulties in defining continuous differentiable distance constraints.
An approach for pseudo-density methods is to use a layer of ghost points that lie a short distance outside the design domain \citep{zhang2018geometry}. The pseudo-densities at ghost points are then computed and if any value is non-zero, then a component must be outside the domain. This idea is used to create a aggregated constraint function that ensures the maximum pseudo-density value at all ghost points is zero. This approach can be used for both convex and non-convex design domains without modification.


\section{Feature-mapping methods for shape optimization}
\label{sec:shape_opt}
In this section we discuss the application of feature-mapping to solve shape optimization problems, which is essentially a classical shape optimization approach with the design mapped to a fixed-grid. To this end we start with a brief description of what we consider as \textit{classical shape optimization}.

\subsection{Classical shape optimization}
\label{sec:classical_shape_opt}
In classical shape optimization, only the structural part is discretized using finite elements. The structural interface is exactly modeled, which is for some applications an essential feature.  Void regions are not part of the finite element analysis, which can significantly reduce the computational effort.

As a consequence, the boundary mesh nodes are moved during optimization. To maintain mesh quality for accurate analysis, mesh smoothing and/or re-meshing is necessary. However, re-meshing can become rather involved and if the quality of the finite element approximation is insufficient, there is the potential for the mesh to be optimized for numerical artifacts, similar to the checkerboard effect in density-based topology optimization. 

Classical shape optimization has been  successfully applied over several decades. However, the mathematical and technical realization is rather involved, especially compared to density-based topology optimization. The mathematical approach is often formulated in an infinite dimensional setting, see e.g. \cite{Sokolowski:1992:Book} or \cite{Haslinger:2003:Book} and differentiable mesh generation must be provided \citep{Haslinger:2003:Book}. The gradient information is based on the shape gradient.

Shape optimization is performed with a wide range of different parameterizations. These can be categorized as either boundary-node-based parameterization, or higher-order forms of design parameterization. We begin with the first variant where each surface node is a design variable. This is called the independent node movement approach \citep{Imam:1982:3DShape} or parameter free shape optimization. This provides a large space of admissible shapes, but it comes with its own challenges in terms of regularization and feature size control, see e.g. \cite{Le:2011:ParameterFreeShape}. Closing and creation of holes are generally difficult to achieve, or even impossible. During the optimization process, insertion or deletion of boundary nodes may be necessary, as well as re-meshing. This generally prevents the use of first-order mathematical programming algorithms. As a consequence, constraint functions need to be handled indirectly. Furthermore, no rigorous convergence criteria are available.   

Aside from the parameterization of boundary nodes there are many  variants of higher-order parameterizations established in shape optimization, see \cite{Haftka:1986:Survey} for an early survey. Conveniently, this corresponds to the construction of geometries by spline functions in computer aided design (CAD). Here, the mapping from the design parameters onto the boundary nodes is differentiable and thus allows gradient-based optimization, see \cite{Braibant:1984:B-splines}. The separate meshing of the geometry can be alternatively handled by isogeometric shape optimization, see \cite{Wall:2008:Isogeometric}. Provided a differentiable parameter-to-boundary mapping, the shape sensitivity of an arbitrary parameterization can be obtained from the nodal shape gradient by the chain rule.

\subsection{Using pseudo-density feature-mapping}
\label{sec:shape_dens}
Some major advantages of shape optimization are the crisp boundary description and the wide and versatile range of design parameterizations. However, the modeling and technical realization is often quite involved. Density-based topology optimization is known for its elegant and easy modeling, standard approaches for sensitivity analysis and straight-forward technical realization, e.g. no re-meshing is necessary. On the other side, the extremely rich design space allows only an implicit design description and is, for some applications, difficult to control. Enforced by standard regularization approaches, one has to deal with a more or less blurred interface description which is anyway rasterized by the fixed analysis mesh.

It is now a natural approach to combine features of both worlds, density-based topology optimization and parameterized shape optimization. This can be seen from two points of view. One can perform standard parameterized shape optimization, but skip the re-meshing by using a fixed analysis grid, leading to an inexact boundary modeling. The same picture is given from the point of view of density-based topology optimization, where a higher level design parameterization is mapped to the pseudo-density field, allowing much closer control of the design and, as such, intentionally sacrificing the rich design space of density-based topology optimization. Under certain circumstances, this approach can be seen as a simplified shape optimization where, aside from the design to density mapping, the technical implementation and sensitivity follows that of density-based topology optimization.  

In principle, any higher-order classical shape parameterization can be combined with feature-mapping, by simply mapping the boundary transition to a fixed-grid using pseudo-densities (see \secref{sec:density_mapping}).

\begin{figure}[ht!]
  \centering
  \subfloat[B-spline]{   
  \includegraphics[width=0.22\textwidth]{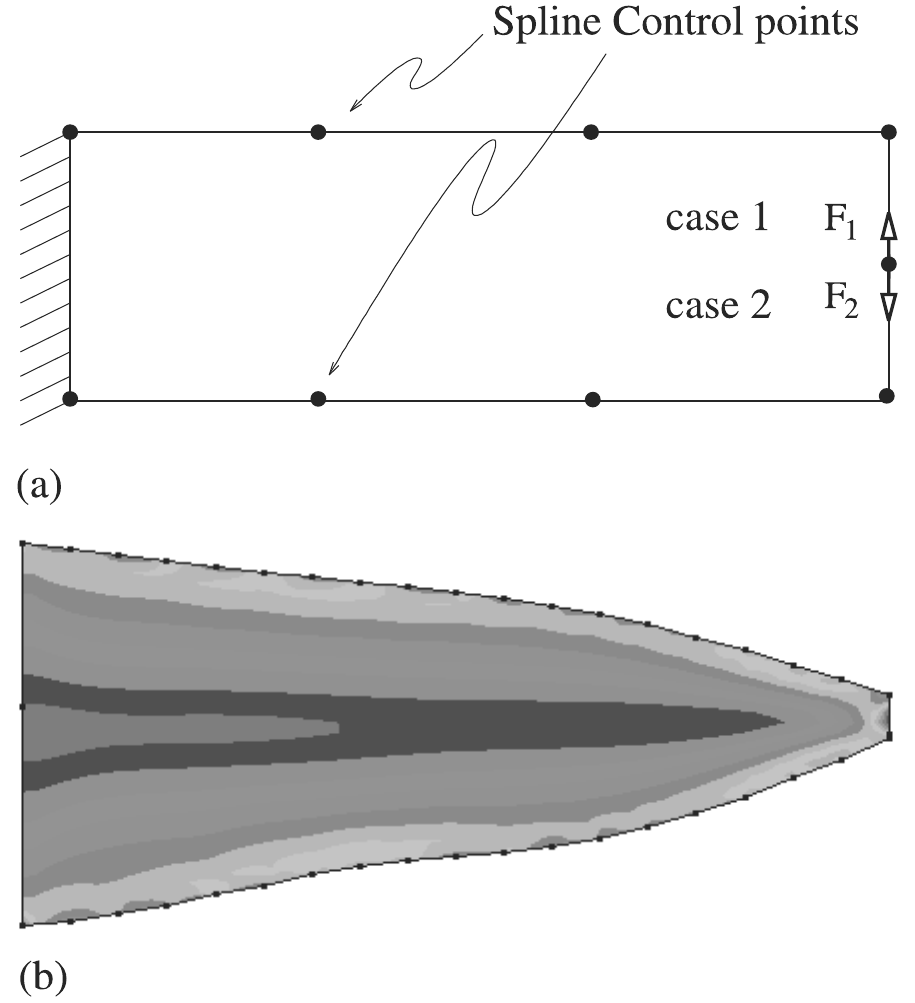}} \,
  \subfloat[first order spline]{   
  \includegraphics[width=0.22\textwidth]{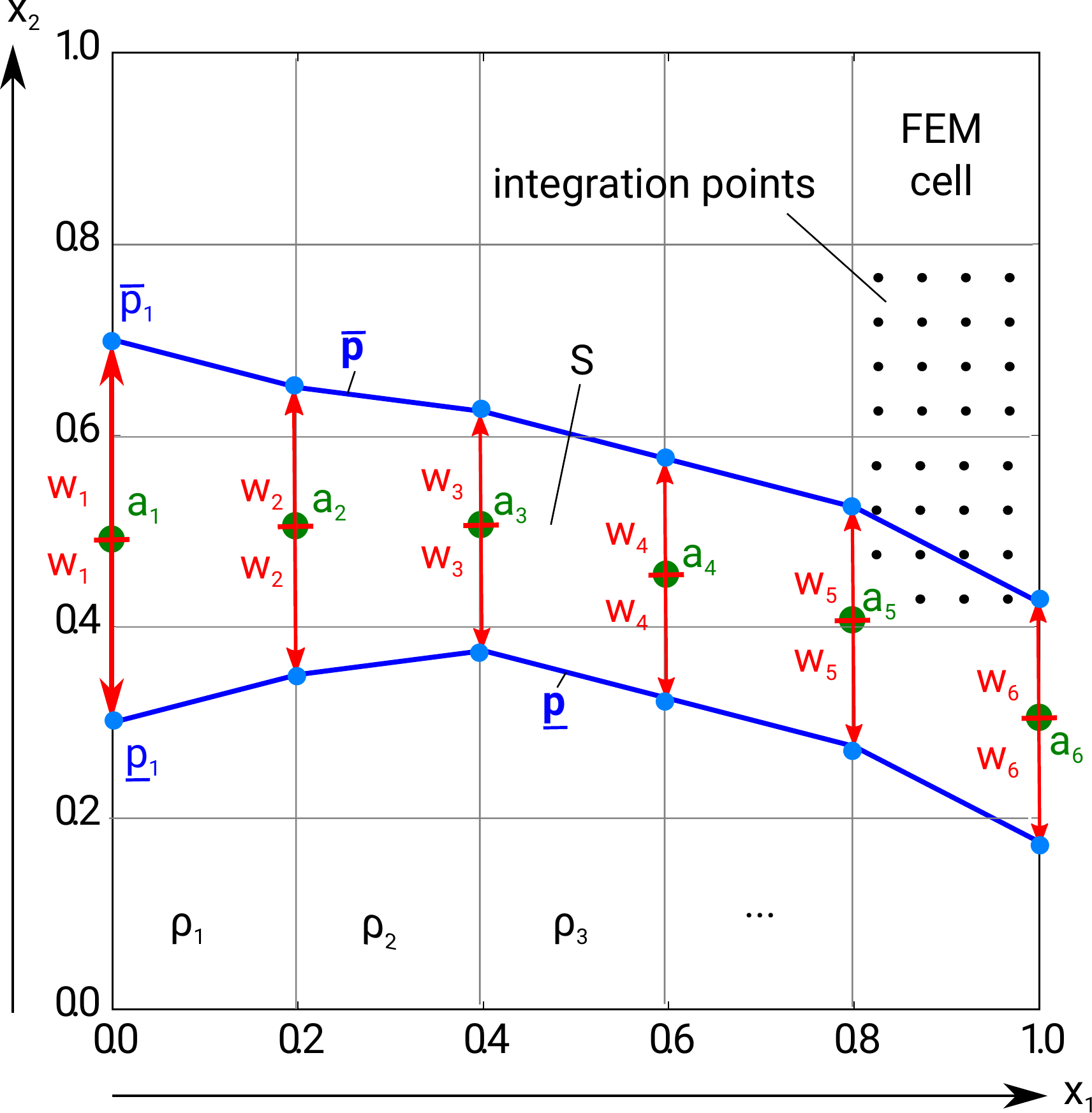}} 
  \caption{(a) B-spline parameterization in \cite{Garcia:2004:EvolShape}; (b) piecewise linear spline with mesh aligned control points in \cite{wein2018combined}.}
  \label{fig:spline}
\end{figure}

In \cite{Garcia:2004:EvolShape} B-splines are optimized with a small number of control points, see \figref{fig:spline}. The design parameters are the location of the control points. The optimization in \cite{Garcia:2004:EvolShape} is performed without sensitivity analysis using a gradient-free method. 

The first rigorous high-level geometry-to-density mapping was demonstrated in \cite{Norato:2004:Shape}, where boundary smoothing is discussed and the convergence of the rasterization of an explicitly given primitive is shown theoretically and numerically. Optimization results based on implicit shapes (via radial basis functions) are also shown, but the approach is not discussed in detail. 

In \cite{wein2018combined}, splines of order 1 (piecewise linear functions) are used, where the control points are aligned to the analysis mesh. The design space is rather restricted, as either the $y-$ or $x-$component of the control points are design variables for horizontal or vertical structures, respectively. Additionally, the width is parameterized by thickness variables. This facilitates, but also requires, the control of first and second spatial derivatives of the design variables in the form of slope constraints and curvature constraints, which are both directly defined from the design variables. Also, as the strips are primarily aligned with one of the axes, see \figref{fig:spline}, the computation of an approximate signed-distance function is easily obtained and is differentiable with respect to the control points.

\begin{figure}[ht!]
  \centering
  \includegraphics[width=0.26\textwidth]{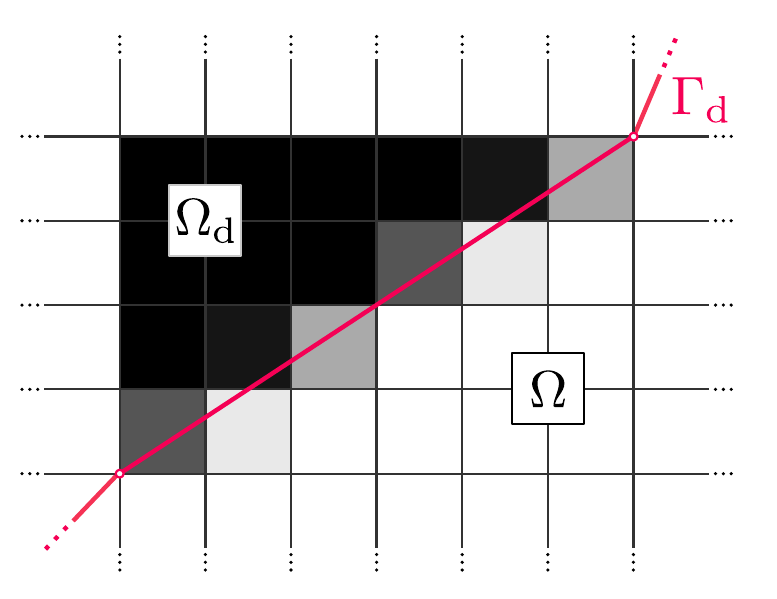} \,
  \includegraphics[width=0.18\textwidth]{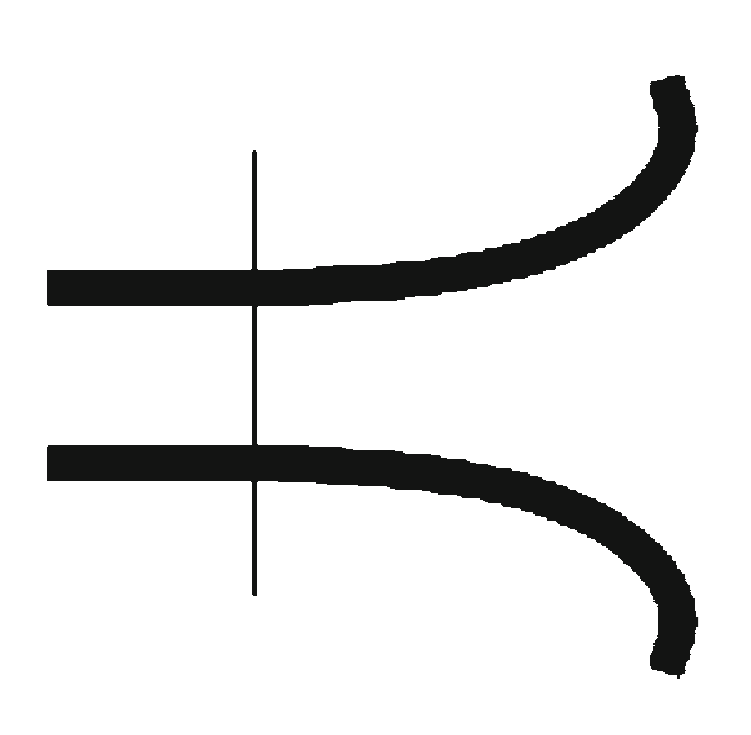} 
  \caption{Optimization of the angle between line segments in \cite{Kasolis:2012:Horn}. The right figure shows an optimized acoustic horn in 2D.}
  \label{fig:horn}
\end{figure}

An interesting variant of the first-order spline representation is used in \cite{Kasolis:2012:Horn}. Here the angle $\vartheta_i/l$ between line segments of length $l$ expressing the boundary $\Gamma_\text{d}$ is the design variable, see \figref{fig:horn}. The application is the optimization of an acoustic horn. The line segments are larger than the mesh element size and no regularization is necessary. The structure is of constant thickness and fixed on the left side. A change in the angle $\vartheta_i$ leads to a rigid body movement of all line segments to the right. The design to pseudo-density mapping is differentiable. 

\subsection{Using immersed-boundary feature-mapping}
\label{sec:shape_immersed}
Combining immersed-boundary feature-mapping with higher-order classical shape parameterization has also been demonstrated. Again, in principle, any higher-order parameterization can be combined with feature-mapping using an immersed boundary method, such as XFEM.

\cite{van2007stress} use XFEM and a fixed-grid to optimize the shape of a 2D fillet for stress minimization. The fillet is parameterized by an implicit geometry description, such as a super-circle, or super-ellipse. \cite{noel2016analytical} use XFEM and feature-mapping for shape optimization of bimaterial structures, such as those with elliptical stiff and soft inclusions. The same framework is also used to optimize the shape of inclusions when designing micro-structural material layout to minimize the maximum stress \citep{noel2017shape}.

The interface-enriched generalized FEM (IGFEM) is used by \cite{najafi2015gradient} for shape optimization of bimaterial structural and thermal problems, where an explicit geometry parameterization is used, e.g. circular and ellipitcal shaped inclusions. XFEM and IGFEM both use a fixed-grid and enrichment functions. The main difference is that the additional degrees of freedom in XFEM enrichment are added to the fixed-grid nodes, whereas in IGFEM  they are added to points where the interface intersects a fixed-grid element edge.

\section{Hybrid feature-mapping / free-form methods for topology optimization}
\label{sec:hybrid}
Hybrid methods are defined as those that combine feature optimization with free-form topology optimization, where free-form means topology optimization without any high-level feature control (such as conventional density-based and level-set methods). Features can be solid or void and feature design variables can include position, orientation, size, shape and/ or number. Historically, the development of hybrid methods began before pure feature-mapping methods for topology optimization (`pure' meaning without any free-form optimization). Some techniques developed for hybrid methods were later exploited in the development of pure feature-mapping methods.

\subsection{Combining free-form with features}
\label{subsec:combining-free-form-features}
A key technical challenge of hybrid methods is to combine the properties of features with the properties of the free-form structure in a single analysis model. Ideally, the method should enable gradients of all design variables to be computed, so that an efficient gradient-based optimization approach can be used to solve the problem. This is achieved by combining feature and free-form geometry using either a map-then-combine, or combine-then-map approach (as discussed in \secref{sec:combination}).

\cite{qian2004optimal} used an implicit Gaussian, or peak function to map the stiffness of rectangular solid features onto the analysis mesh via
\begin{equation}
\label{eqn:peak_func}
    \hat{E}_i(x,y) = E_i \exp\left ( -\frac{(\Delta \Tilde{x}_i)^{\eta}}{\sigma^2_{xi}} -\frac{(\Delta \Tilde{y}_i)^{\eta}}{\sigma^2_{yi}}\right)
\end{equation}
with
\begin{equation}
\label{eq:peak_rotate}
    \left( 
  \begin{array}{l}
  \Delta \Tilde{x}_i \\
  \Delta \Tilde{y}_i
  \end{array} \right) = 
   \left(
  \begin{array}{ll}
  \cos{\theta_i} & \sin{\theta_i} \\
  -\sin{\theta_i} & \cos{\theta_i}
  \end{array} \right) 
  \left( 
  \begin{array}{l}
  x - x_i \\
  y - y_i
  \end{array} \right),
\end{equation}
where $\hat{E}_i(x,y)$ is the stiffness of the feature at location $(x, y)$, $E_i$ is the stiffness of the feature material, $x_i$, $y_i$ are the location design variables and $\theta_i$ is the orientation. If the exponent $\eta >= 4$, then the feature is approximately rectangular and the $\sigma$ values control the dimensions.

The peak function was originally developed for multi-material topology optimization. Here, it is used to combine the stiffness of the features with the stiffness of free-form structure using a simple summation as
\begin{equation}
\label{eq:hy_comb_sum}
    E(x,y) = \hat{E}_0(x,y) + \sum_i^N{\hat{E}_i(x,y)},
\end{equation}
where $N$ is the number of components and $\hat{E}_0$ is the stiffness of the free-form structure, which is defined by element-wise density design variables and penalized using a similar peak function. This map-then-combine approach allows straight-forward computation of sensitivities, as the effective stiffness of an element in the fixed mesh is a linear combination of free-form and feature design variable values. The multi-material formulation also allows for features and free-form structure to have different properties. However, the implicit geometry representation in \eqnref{eqn:peak_func} is limited by the shapes it can represent.

Several methods define features using implicit functions with more geometric freedom. A popular choice is to define the boundary of a feature as the zero level-set of an implicit scalar function \eqnref{eqn:imp_func}.

\cite{chen2007shape} represent both the features and free-form structure using implicit functions. The benefit of this approach is that the free-form structure can be combined with arbitrary features using a combine-then-map approach with Boolean operations, resulting in a single implicit function that describes the overall structure. However, Boolean operations cannot be differentiated. Thus, smooth approximations of Boolean operations are used instead, such as smooth R-functions (see \secref{subsec:smooth_comb}). However, this combine-then-map approach makes solid features become indistinguishable from the free-form structure in the analysis model. Thus, solid features are limited to have the same material properties as the free-form structure.  R-functions have also been used in other hybrid methods to create complex features by combining primitives represented by implicit functions \citep{xia2013implicit}.

If features are defined using implicit functions then a smoothed Heaviside function,,e.g. \eqref{eq:smoothHS}, can be used to map feature properties onto a fixed-grid using pseudo-densities \citep{xia2013implicit, zhang2015explicit}. This method is discussed in more detail in \secref{subsubsec:combine-then-map-implicit}. Using the smoothed Heaviside also allows for features with different stiffness values to be combined in the analysis model. For example, \cite{xia2013implicit} combined implicit feature representations with an element density-based free-form structure representation and used the following map-then-combine approach to define the Young's modulus within each element of the fixed mesh
\begin{equation}
\label{eq:E-comb_Hs}
    E_e = \rho^p \left(E_0 + \sum_{i=1}^{N} \widetilde{H}(\phi_i)(E_i - E_0) \right).
\end{equation}
This formulation works for void features where $E_i \approx 0$, although it does not inherently prevent element density variables beneath solid features from becoming zero, potentially leading to changes in the shape or topology of a feature. However, this can be avoided in the optimum for the minimization of compliance problem (subject to a volume constraint on the free-form structure) by excluding elements associated with solid features from the volume calculation.

It is also possible to use an exact Heaviside function when mapping implicit void features onto the fixed mesh \citep{kang2013integrated}. The free-form structure is parameterized using element-wise density variables and the pseudo-density within an element, centred at $x$, is simply computed by
\begin{equation}
\label{eq:hy_dens_Hs}
    {\rho}(\vec{x}) = \tilde{\rho}(\vec{x}) H(\phi_\omega(\vec{x})),
\end{equation}
where $\rho$ is the physical density used to formulate the finite element stiffness matrices and $\tilde{\rho}$ is obtained using a density filter of the free-form density design variables. The inherent discontinuity of the exact Heaviside function does not allow explicit derivative calculation for feature location and orientation. Thus, shape derivatives (obtained from the analytical pre-discretized governing equations) are used to update the location and orientation variables. These shape derivatives are then used to obtain velocity components for feature location and orientation, which are then used with the Hamilton-Jacobi equation to update the implicit functions that define the void features. However, using this approach does not allow for feature design variables (location and orientation) to be combined with free-form structure density variables in the same nonlinear programming optimization problem. Thus, to achieve simultaneous optimization of free-form structure and feature location / orientation, the feature velocity components are included as design variables when using a gradient-based optimizer.

\cite{li2017integrated} introduced the stiffness spreading method to combine properties of explicit geometry features and the free-form structure. The method starts with a fitted finite element mesh of the component whose shape does not change, but it only moves and rotates. To combine this with the free-form structure, defined on a fixed-mesh, the component stiffness matrix is transformed into an  equivalent matrix that has the same dimensions and degrees of freedom as the fixed-mesh.
The free-form structure and transformed feature matrices are then simply summed for analysis. Thus, this can be considered a map-then-combine approach. The transform essentially provides a linear map between a solid feature and fixed-mesh degrees of freedom. This could be achieved by a local transform, where feature nodes are simply mapped to the closest fixed-mesh nodes (i.e. those associated with a single element). However, this creates a discontinuity in the derivatives of the feature location and orientation. Thus, a non-local transform is used, where the connection is smoothed, or spread over several nodes in the fixed-mesh within a predefined radius.

Multi-point constraints (MPCs) have also been used to connect explicit solid  features to the fixed-mesh used to parameterize the free-form structure \citep{zhu2015multi,gao2015improved}. The solid features are meshed independently and connected to the free-form structure fixed-mesh at a number of predefined locations that represent bolts or rivets. This approach is in contrast to other methods that assume a perfect bonding at the interface between solid features and free-form structure. The use of MPC connections has also been extended to a multi-frame problem, where component locations and free-form structure for a number of frames are simultaneously optimized with the location of the frames in a larger design space \citep{zhu2017integrated}. MPCs are also used to connect frames to the free-form structure in the larger design space.

The above methods all aim to map feature properties onto a fixed analysis mesh in the spirit of feature-mapping methods as defined in this review. The alternative approach of re-meshing has also been used in hybrid methods. Solid features of fixed size and shape are discretized using a conforming mesh, whereas a fixed regular mesh is used to parameterize the free-form structure. The two meshes are combined into a single conforming mesh for analysis. This is achieved by first placing feature meshes over the fixed mesh. Elements in the fixed mesh that are partly covered by a feature are then divided, or remeshed, to create a single conforming mesh for analysis \citep{zhu2008simultaneous}.

This local re-meshing at the boundary leads to changes in the free-form mesh whenever a feature moves. Thus, the usual element-wise density parameterization for the free-form structure cannot be used. Instead, the concept of density points is introduced, where the density at a point influences several elements within a local domain \citep{zhu2008simultaneous}.

Re-meshing approaches add computational complexity and cost (especially in 3D) compared with the fixed-grid mapping methods. Also, semi-analytical derivatives of feature design variables are required, where the derivative of the global stiffness matrix with respect to a small perturbation in the variables is computed using finite differences \citep{zhang2011some}. This process also adds significant computational cost, as the mesh is perturbed and stiffness matrices recomputed for each feature design variable.

Techniques have been proposed to reduce the computational cost of semi-analytical derivatives. \cite{xia2012superelement} used a superelement formulation to reduce the size of feature stiffness matrices. The same authors also proposed using a smoothed Heaviside approximation of the boundary to facilitate full analytical derivatives \citep{xia2012sensitivity}. The idea is to compute stiffness matrix derivatives based on the change in material properties as a feature moves. The smoothed Heaviside approximation is required to smooth the sharp discontinuity in material properties. However, re-meshing is still performed when features move.

One benefit of re-meshing is that the interface between solid features and free-form support structure is explicitly modeled. This could be important if the bonding  between solid features and free-form support structure is included in the analysis as a nonlinear interface behavior. A recent example of this was presented by \cite{liu2018integrated}, where a cohesive zone interface model was used and the dissipated interface energy added to a compliance objective function. It was observed that features move to positions where the free-form support structure is under compression, as this helps minimize the dissipated interface energy.

XFEM is an alternative to re-meshing that can also capture the explicit interface (see \secref{sec:xfem}). Hybrid methods that use XFEM also use implicit functions to describe features, as this provides a natural way to identify elements in the fixed-mesh that contain an interface between the free-form structure and a feature. Furthermore, analytical derivatives of feature location / orientation design variables can be obtained using the chain rule \citep{zhang2012integrated}. XFEM has been combined with density-based \citep{zhang2012integrated,wang2014topological}, level-set based \citep{zhou2013engineering} and parameterized level-set \citep{liu2014level} methods for simultaneous optimization of features and free-form structure.

In summary, various methods have been proposed to combine properties of features and free-form structure in hybrid methods. Implicit and explicit features have been combined with density and level-set-based free-form structure. A key difference in the methods is whether they follow a combine-then-map or map-then-combine approach, as only the latter allows different material properties for solid features and free-form structure.

\subsection{Feature design variables}
\label{sec:hybrid_vars}
The original motivation for developing hybrid methods was to simultaneously optimize the location and orientation of a number of solid features (with predefined shapes) and a free-form structure to support the features \citep{qian2004optimal}. The solid features represent components to be embedded within a free-form structure and their properties (such as stiffness) contribute to the overall performance of the structure. Thus, by simultaneously optimizing the free-form structure and position / orientation of solid features, an overall optimum can be found. 

This idea also applies to void features of fixed size, which could, for example, represent necessary cutouts for systems to pass through the structure. Again, by simultaneously optimizing the position / orientation of void features and the free-form structure, an overall optimum can be found.

Beyond features of fixed shape and size, hybrid methods have been proposed that can also simultaneously optimize the size and/or shape of features. \cite{cai2015stress} demonstrated the simultaneous optimization of circular hole radius and support structure. The hole had a fixed location and the solution was trivial, as the optimal radius is simply the lower bound. The method presented by \cite{zhou2013engineering} can simultaneously optimize the location, orientation and shape of features. Furthermore, the method can be used to optimize only the shape of selected feature edges, thus providing several levels of control over feature geometry. This is achieved by modeling the free-form and feature geometries using implicit functions. These are then combined by multiplying Heaviside functions to give an indicator function
\begin{equation}
\label{eq:prodH}
    \chi_\omega(\vec{x}) = \prod_i^M H(\phi_i(\vec{x})),
\end{equation}
where $M$ is the total number of implicit functions (i.e. one for the free-form structure and $M-1$ features). The indicator function is then used to map the combined structure onto the fixed mesh for analysis (i.e. a combine-then-map approach). Note that this formulation assumes solid features and free-form structure have the same material, although the framework does allow for regions with different material properties.

The method presented by \cite{liu2014level} also uses implicit functions to represent both free-form and feature geometries, which are combined together using R-functions (also a combine-then-map approach). The features considered are primitives that can also change shape and size. For example, the  semi-major and semi-minor axes of an elliptical void feature were simultaneously optimized with the free-form structure. In addition, \cite{liu2014level} introduced the idea of a tolerance zone, which is a fixed-width region of material around a void feature. This was achieved by ensuring the presence of free-form structure within a certain region around void feature edges.

\cite{lin2015isocomp} proposed a hybrid method that can optimize the size, shape and number of features. Features are defined explicitly using splines and the free-form topology is parameterized using a pseudo-density type approach. A criterion based on the topological derivative is used to create new features, such as holes and inclusions.

\subsection{Optimization strategy}
\label{sec:hybrid_opt}
Most hybrid methods aim to simultaneously optimize the features and free-form structure. However, it has long been recognized that the hybrid problem formulation with movable features is non-convex and the solution is often highly dependent on the starting position of the features \citep{qian2004optimal,zhu2015multi,cheng2018coupling,liu2018integrated}.

A simple method to overcome this is to solve the problem several times with different starting positions of the features, although this adds computational cost. Another approach is to initially optimize just feature positions/ orientations. Then, after a set number of iterations, or after the feature-only problem converges, both features and free-form structure are optimized simultaneously \citep{kang2013integrated,wang2014topological,li2017integrated}. This strategy may provide a better design, as features can get stuck near their starting positions if the simultaneous approach is used from the start.

A variant of this approach is to fix feature designs and positions after the initial optimization \citep{bakhtiarinejad2017component}. This results in a sequential optimization approach that could be repeated several times, i.e., optimize the feature design and then the free-form structure alternately, until convergence \citep{xia2012superelement}. An advantage of the sequential approach, is that different optimization strategies can be applied to feature design and free-form structure.

\subsection{Applications}
\label{sec:hybrid_apps}
Most hybrid methods have been applied to maximize the stiffness (minimize compliance) of the combined feature and free-form structure. Thus, solid features are placed to best utilize their stiffness, whereas void features are placed to least disrupt the free-form structure from obtaining the stiffest design.

 Developments for other applications include the design of compliant smart structures with embedded movable actuators \citep{wang2014topological}, where the actuator features  provide an input force and the free-form structure is used to design the compliant part of the mechanism. By simultaneously optimizing embedded piezoelectric actuator placement and the compliant mechanism design, an overall optimum can be found. \cite{cai2015stress} used features and free-form structure defined by implicit functions and combined by R-functions to optimize structures with stress constraints. Finally, there have also been applications for heat transfer problems \citep{li2017integrated, cheng2018coupling}, where the position of solid features are optimized to take advantage of their thermal conductivity and/or heat flux boundary conditions.

Currently, there are few applications of hybrid methods beyond designing stiff structures. Thus, there is potential for future research into other application areas, such as electromagnetic and electrical-mechanical problems, fluid flow and fluid-structure interactions problems and problems involving nonlinear mechanics.

\section{Feature-mapping methods for topology optimization}
\label{sec:TO-feat-map-methods}

In this section we discuss feature-mapping methods for topology optimization, where the structure is exclusively defined by the combination of high-level parametric descriptions of voids in a solid design region, or of solids in a void design region.  For brevity in the discussion, we will refer to the former as holes and the latter as components.  

As a historical note, a mention is due to the bubble method \citep{eschenauer1994bubble}, which describes holes in a solid region using B-splines and is, in fact, one of the first topology optimization techniques.  However, this method does not fit into the family of methods covered in this review, as it employs a conforming mesh and thus requires re-meshing upon design changes, and also lacks a mechanism to merge holes.  

\subsection{Combine-then-map methods}
\label{subsec:TO-combine-then-map-review}

\begin{figure}[ht!]
	\centering
	\subfloat[]   
	{\includegraphics[width=.225\columnwidth,valign=c]{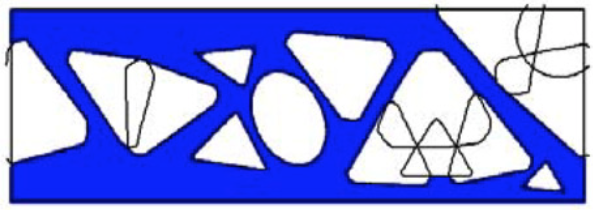} %
	\vphantom{\includegraphics[width=.225\columnwidth,valign=c]{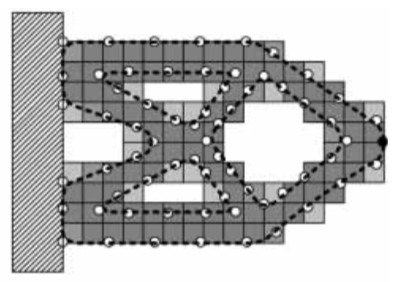}}} \,
	\subfloat[]
	{\includegraphics[width=.225\columnwidth,valign=c]{sect8_fig1_2008Lee_Kwak.png}} \\
	\subfloat[]   
	{\includegraphics[width=.225\columnwidth,valign=c]{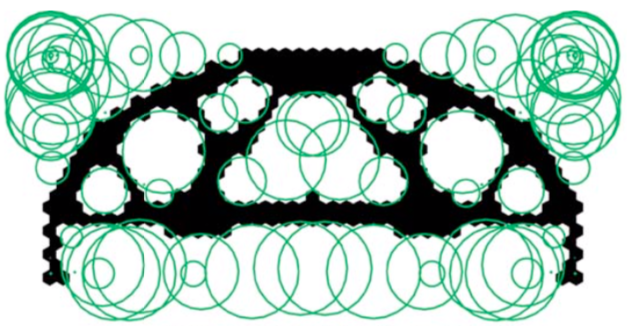}}  \,
	\subfloat[]
	{\includegraphics[width=.225\columnwidth,valign=c]{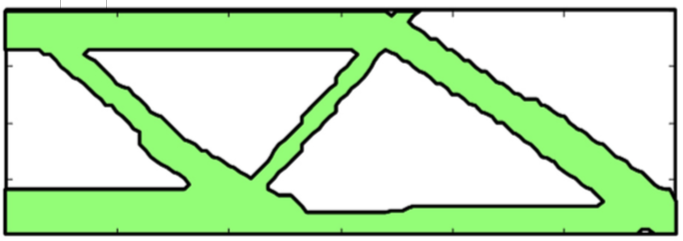}
	\vphantom{\includegraphics[width=.225\columnwidth,valign=c]{sect8_fig1_2011_Saxena_MMOS.png}}} \\	
	\subfloat[]   
	{\includegraphics[width=.225\columnwidth,valign=c]{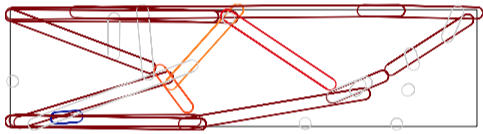}
	\vphantom{\includegraphics[width=.225\columnwidth,valign=c]{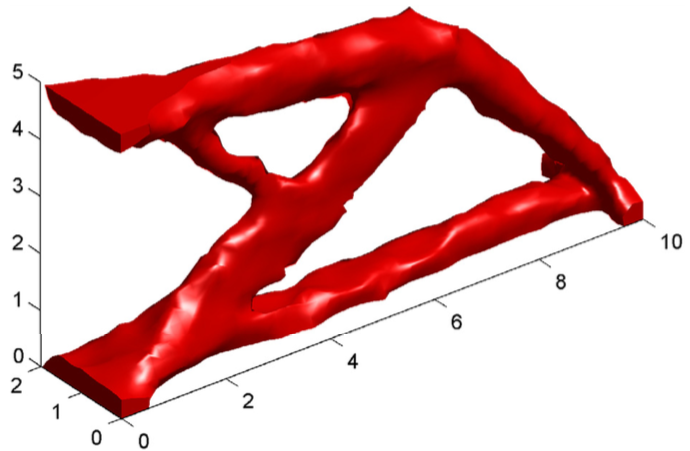}}}  \,
	\subfloat[]
	{\includegraphics[width=.225\columnwidth,valign=c]{sect8_fig1_2017_Zhang_et_al_MMV.png}} 	
	
	\caption{Examples of different feature-mapping methods for topology optimization: a) \cite{mei2008feature}, b) \cite{lee2008smooth}, c) \cite{saxena2011circular}, d) \cite{guo2014doing}, e) \cite{norato2015geometry} and f) \cite{zhang2017explicit}.}
	\label{fig:sample-feature-mapping}
\end{figure}

As described in \secref{sec:combination}, feature-mapping methods for topology optimization render structures exclusively made of the combination of geometric features. Examples of these methods are shown in \figref{fig:sample-feature-mapping}.
The pioneering work in combine-then-map methods was introduced in \cite{cheng2006feature, mei2008feature}.\footnote{The latter work seems to be a more detailed journal version of the former, hence from hereon we only reference the latter.}  This method formulates and successfully demonstrates all of the hallmark features of methods in this category, as well as other features that have not been explored elsewhere at the time of writing this review, such as the primitives interpolation strategy described in \secref{subsec:TO-geometric-reps}. To our surprise, this contribution went unnoticed in many subsequent works; only until recently have some works started citing it. We hope our review aids in restoring credit to this foundational contribution. In this method, holes are combined using (non-smooth) R-functions. 

The first combine-then-map method to employ primitive-shaped solid components using gradient-based optimization is the method of moving morphable components (MMC) of \cite{guo2014doing}. This work has many derivatives and, at the time of writing this review, is the most cited feature-mapping topology optimization approach.  It uses a true maximum function to combine components into a single implicit function.  Virtually all combine-then-map methods with implicit representations use a non-differentiable maximum.  Exceptions are the works by \cite{zhou2016feature, sharma2017advances,du2019moving}, which use a smooth Kreisselmeier-Steinhauser (KS) approximation to combine the signed-distance representations, and the works by \cite{zhang2016minimum, xie2018new}, which use smooth R-functions.  The work by \cite{sharma2017advances} smooths out the combined level-set function with an anisotropic filter to produce smooth boundaries, although it indicates that, for the examples presented, such filtering is not necessary to produce good results.

The works by \cite{lee2007smooth,lee2008smooth} and \cite{kim2008smooth} are the first to employ an explicit representation, namely B-splines. Holes are merged by deleting control points, as outlined in \secref{subsubsec:combine-then-map-explicit}.  \cite{wang2009structural} is perhaps the first work to design structures made exclusively of solid components, consisting of wide curves modeled using B{\'e}zier curves.  The pseudo-density mapping in this method is not differentiable, as described in \secref{subsec:TO-analysis-approaches}, thus this method uses a genetic algorithm for the optimization. 

In the material mask overlay strategy (MMOS) of \cite{saxena2011circular}, holes are represented as the union of primitive-shaped regions (the `masks'), such as circles, rectangles and ellipses. The mapping to the analysis mesh uses a small or unity pseudo-density depending on whether an element's centroid is inside or outside any of the masks, respectively. As this mapping is not differentiable, a gradient-free hill-climbing method is used for the optimization. \cite{wang2012high} introduced a differentiable version of MMOS for the design of photonic waveguides that combines a free-form density field with circular holes. The boundary of the holes is smoothed as described in \secref{sec:boundary_smoothing} and holes are combined using a product of the smooth Heavisides, similarly to \eqnref{eq:prodH}. In the examples provided in this work, however, the circular holes do not overlap.  \cite{hoang2017topology} present a modification of this differentiable MMOS approach for designs with solid components.

%
%

\subsection{Map-then-combine methods}
\label{subsec:TO-map-then-combine-review}

The most common method that embodies the map-then-combine strategy is the geometry projection method of \cite{bell2012geometry, norato2015geometry}. This method employs a volume fraction calculation using a circular (in 2D) or spherical (in 3D) sample window (cf.~\eqref{eq:geom_proj}) to obtain a pseudo-density for each of the components. To combine features, the preliminary work by  \cite{bell2012geometry} used the minimum signed distance to any of the components, which is equivalent to a true maximum. The work by \cite{norato2015geometry} introduced the size variable that allows the complete removal of a component via penalization (cf.~\eqref{eq:effective-density}).  It also uses a $p$-norm (cf.~\eqref{eq:p-norm}) as a smooth approximation of the maximum to perform the combination of the penalized pseudo-densities. 

A mention is due to the works of \cite{guest2012casting} and \cite{ha2014optimizing}, which remove holes or produce components of a fixed shape and size (e.g., circular), respectively, by using projection filters with the given shape on a grid of density variables.  In these methods, the design is described by a density field instead of a high-level parametric description of the holes or components, hence they do not fall under our categorization of feature-mapping techniques.  The shape of the hole or component determines the neighborhood of elements that contribute to the filtered density in each element. 

\subsection{Geometric representations}
\label{subsec:TO-geometric-reps}

Most of the combine-then-map methods that use upfront implicit representations use hyperellipses to represent various primitives. MMC and most of its derivatives employ a direct hyperellipse representation to model 2D bars  \citep{guo2014doing}, and 3D bars and plates \citep{zhang2017new, sharma2017advances}.  \citep{zhang2016new} and \cite{guo2016explicit} modify the hyperellipse equation to produce curved bars with linearly, quadratically or sinusoidally varying width.  Other methods, notably \cite{mei2008feature} and \cite{zhou2016feature}, use a hyperellipse representation as well, but they convert it to a signed-distance representation, which, as mentioned in \secref{sec:boundary_smoothing}, can lead to better convergence.  The geometry projection method in \cite{norato2018topology} models 2D solid components using supershapes, which correspond to a generalization of the superellipse formula that has variable symmetry, and thus can approximate a wide range of primitive shapes with a single equation. As in all geometry projection techniques, a signed distance to each supershape is computed to obtain the mapped pseudo-density.

 \cite{lee2007smooth,lee2008smooth} and \cite{kim2008smooth} modeled holes with B-splines.  The 2D solid components in  \cite{wang2009structural} are modeled using wide B{\'e}zier curves, with design variables being the control point locations. The curves connect support and loading points and a constraint in the optimization is required to prevent curve self-intersections. 2D \citep{zhang2017structural,du2019moving} and  3D \citep{zhang2017explicit} B-splines are used to model holes in the Moving Morphable Voids (MMV) method, which is a variant of the MMC method.  Control points inside the intersection are made inactive, and the resulting combined B-spline is converted to a signed-distance function prior to mapping to the analysis. 
 
An implicit geometry representation common in geometry projection works is offset surfaces (cf.\ \cite{bloomenthal1990interactive}), in which the component boundary is given by the set of points equidistant to a medial axis or surface.  For instance, a 2D bar is given by a rectangle with semicircular ends \citep{norato2015geometry}, a 3D bar by a cylinder with semispherical ends \citep{watts2017geometric,kazemitopology}, and a 3D plate as a cuboid with semicylindrical edges and quarter-sphere corners \citep{zhang2016geometry}.  Offset bars have also been used in conjunction with the MMC method in \cite{deng2016design}  to produce designs where bars are connected at all times by sharing endpoints of their medial axes, effectively rendering a ground structure approach; and in \cite{hoang2017topology} to impose a minimum thickness in the structure (see \secref{subsec:TO-complexity}). Offset plates with a single curvature radius are modeled in \cite{zhang2018geometry}.  The offset surface representation is combined with a free-form density field in \cite{zhang2016geometry} to produce designs made of plates with free-form holes.

An interesting geometric representation is employed in  \cite{mei2008feature}, whereby each hole is represented as a weighted sum of prescribed geometric primitives (e.g., circles and triangles) that are implicitly represented using signed-distance functions.  A constraint is added to the optimization to penalize the primitive weights so that the holes converge to being, for example, either a pure circle or a pure triangle. A similar constraint is used in \cite{norato2018topology} to penalize supershapes parameters so that the optimal design is exclusively made of, e.g., rectangles and ellipses.  \cite{mei2008feature} also use shape matching techniques to take the resulting designs and consolidate them into fewer primitives, after which an additional optimization stage is performed. 

\subsection{Analysis aproaches}
\label{subsec:TO-analysis-approaches}

As mentioned in \secref{subsec:combine-then-map}, combine-then-map approaches can perform the analysis using either immersed boundary or pseudo-density techniques. Map-then-combine approaches, on the other hand, only use pseudo-density techniques, although they could in principle use immersed boundary approaches, as discussed in \secref{subsubsec:map-then-combine-HS}. 

Combine-then-map approaches that utilize an implicit representation and pseudo-densities for the analysis (cf.\ \cite{mei2008feature,zhang2016new}) compute a smoothed Heaviside of the combined implicit function (see \secref{sec:boundary_smoothing}). Those that use an immersed boundary approach, have used the XFEM method (cf.\ \cite{guo2014doing, sharma2017advances}).  Other approaches use pseudo-densities and isogeometric analysis (IGA)  \citep{hou2017explicit,xie2018new}. 

Different analysis techniques are used in combine-then-map approaches that use a direct explicit representation (i.e., without conversion to an implicit representation such as a signed distance).  In \cite{lee2007smooth,lee2008smooth} and \cite{kim2008smooth}, the analysis uses a pseudo-density that is computed by determining the intersections of the B-spline holes with the element boundaries, replacing the intersected boundary with a straight line, and using a volume fraction approach, as in \cite{Garcia:1999:FixedGridFE}. The work by \cite{seo2010isogeometric} is similar to these methods, however it uses an immersed boundary mapping approach via IGA; their technique decomposes elements that intersect the structural boundary into multiple cells to perform the integration.  The work by \cite{zhang2017comprehensive} uses the weighted B-spline finite cell method for the analysis. The technique with wide curves of \cite{wang2009structural} uses an element pseudo-density, which is calculated by successively subdividing elements that intersect the structural boundary, and computing the area ratio of the subcells that lie inside the curve (or of overlapping curves). In this approach, a non-gradient-based optimizer is used. 

An interesting analysis technique is presented in the MMV method of \cite{zhang2017explicit}, whereby elements in the void regions are removed from the analysis at each iteration to decrease the size of the analysis problem. These elements may be reintroduced in subsequent iterations if they become non-void.

%

\subsection{Complexity and minimum size control}
\label{subsec:TO-complexity}

We start this section with a brief discussion on complexity control.  By complexity, we here mean the topological genus of the structure, i.e., its number of holes---the higher the genus, the more complex the structure. One of the well known characteristics of most density topology optimization techniques is that, in the absence of a control mechanism (e.g., filtering, or a perimeter or slope constraint), the design is usually mesh-dependent.  This is a byproduct of the fact that in these methods the design representation is tied to the analysis mesh.  An alternative in these techniques is to use a grid of design variables that is independent of the mesh, cf.\ \cite{nguyen2010computational}.  

Feature-mapping methods in general do not suffer from this mesh-dependency because the representation of the components or holes is entirely independent of the mesh.  Moreover, if lower bounds are imposed on the dimensions of the components, it is of course not possible for the optimization to obtain smaller feature sizes. Therefore, as long as the total number of components is kept to a maximum, the designs produced by these methods exhibit similar complexity and member size.  This does not entirely mean, however, that different mesh sizes produce the exact same design due to several reasons.  First, the mesh resolution affects the accuracy of the geometry mapping and the analysis solution and consequently the sensitivities, hence different mesh resolutions can lead to different local minima. Furthermore, in the case of pseudo-density approaches, mesh alignment (in conjunction with the interpolation approach) may also cause small differences, as noted in \secref{sec:benchmark}. Second, even if a lower bound on the component dimensions guarantees a minimum feature size, when components intersect or `touch' they can produce regions whose size is smaller (and possibly zero, i.e., a point intersection) than the desired minimum size.  

In addition to imposing bounds on the geometric parameters to ensure a minimum component size, the work in \cite{sharma2017advances} adds a penalty to the objective function to prevent the onset of plates whose medial surface area is smaller than a specified threshold. The motivation to impose this penalty is that small components with zero shape sensitivities may disconnect from the structure and cause oscillations in the optimization.

Several approaches have been proposed to enforce a minimum size for the combined structure in feature-mapping methods for topology optimization. In \cite{zhang2016minimum}, a constraint is added to impose a minimum distance between any two members that intersect. Since the constraint is applied only to components that overlap, it appears to be non-differentiable. \cite{hoang2017topology} employ bars modeled with offset surfaces (i.e., rectangles with semi-circular ends) and effectively control the minimum thickness by introducing two constraints: one that ensures a minimum volume in a mask that covers the rectangular portion of the bar, and another that limits the amount of intermediate material in the semicircular ends of the bars. The work by \cite{niu2018equal}, which uses rectangular bars and the MMC method, adds a penalty term to the objective function to enforce a minimum size. This term consists of a smoothed Heaviside of the minimum distance between the centerlines of two bars for all pairs of bars; this term is multiplied by another that consists of the angle difference between the two bars (raised to a power), so that the orientation between bars that are closer to being parallel is penalized, but those that are closer to being perpendicular are not. While the foregoing methods are effective, they are associated with specific geometric representations (i.e., bars), and a minimum size control that works for combined structures with components of any shape is currently missing.  The recent work by \cite{wang2019imposing} constrains the minimum distance between the medial axes of nearby or intersecting components to be less than a prescribed value, together with a constraint on the minimum thickness of the components.

A note is worth making with regards to intersections between components, which has an effect on the minimum size at intersections.  For pseudo-density approaches, it is possible that components that should intersect in the optimal design  are not fully connected and there is a small gap between them.  As discussed in \cite{norato2018topology}, this is due to the fact that some smooth approximations of the maximum, such as the $p$-norm of \eqref{eq:p-norm}, approximate the true maximum from above, and thus they may render artificially high pseudo-densities, which increases the stiffness of small gaps between components. This issue could also be associated with, e.g., linear volume-to-stiffness interpolations. However, a tighter convergence tolerance in the optimization can aid in resolving this issue. 

\subsection{Design space modification}
\label{subsec:TO-design-space}

As discussed in \secref{subsec:combination-local-minima}, feature-mapping methods for topology optimization are more prone to falling into poor local minima than free-form methods.  Therefore, the choice of initial design, as well as any mechanisms to adaptively modify it, are important aspects in these methods. 

One strategy employed by several of the techniques that design holes in a solid region is to adaptively introduce holes at locations that attain the most negative value of the topological derivative of the compliance.  This strategy, which has been used by level-set techniques (cf.\ \cite{sigmund2013topology} and references therein), was introduced in feature-mapping methods by \cite{mei2008feature} and later used by other methods (e.g., \cite{lee2007smooth,lee2008smooth, kim2008smooth}). To our knowledge, this strategy has not been used yet introduce solid components into a void domain.  In addition to adaptively introducing holes, the work of  \cite{mei2008feature} also removes holes that do not intersect any solid region to simplify the design. The generative design method of \cite{li2019generating} adaptively adds components to generate tree-like structures for area-to-point conduction problems. It performs a sequence of MMC optimizations; at the end of each optimization run, it adds or removes components at the ends of existing branches depending on whether the leaf-side thickness of the component exceeds or falls short of prescribed thresholds, respectively. 

To prevent entrapment in poor local minima, the MMC method of \cite{zhang2017structuralcomplexity} introduces a function that smoothly approximates the number of `effective' bars in the design. Two or more bars are considered to form a single effective component if they overlap (near) collinearly .  The function examines the intersection between the bars: if the intersected area is greater than some specified threshold, and if the angle between the bars does not nearly equal zero or $\pi$, they are considered separate components.  A constraint is then imposed on the number of effective components.

A `bootstrapping' strategy to produce good initial designs for 2D problems with bars is presented in \cite{weiss2018data}, which uses density-based topology optimization to produce an initial design. This result is first converted to a 0-1 design by thresholding, and then it is skeletonized to produce the medial axis of the structure.  The nodes of this skeleton are connected with straight bars to produce an initial design for the MMC method. 

In \cite{zhang2018finding}, the gradient-based tunneling method is used in conjunction with the geometry projection method to go from one local minimum to a better one.  In essence, after converging to a local minimum (the optimization phase), this technique adds a term to the objective function that makes the current minimum a pole of the modified objective, and uses this function to find another design with an equal or lower objective (the tunneling phase).  If successful, a new optimization phase is started with the original function.  

Finally, taking advantage of the reduced number of design variables in feature-mapping methods, and to attempt to prevent entrapment in poor local minima, recent works have used statistical techniques such as: an evolutionary strategy \citep{bujny2018identification}, Bayesian optimization \citep{sharpe2018design} and machine learning methods such as support vector regression and K-nearest neighbors \citep{lei2019machine} to perform the optimization. \cite{raponi2019kriging} also used a Kriging-based surrogate model to improve the convergence rate and optimal designs when using a statistical optimization technique.

\subsection{Geometric constraints}
\label{subsec:TO-geometric-constraints}

One of the most appealing aspects of feature-mapping methods is the ability to impose geometric constraints. These may be motivated by, for example, wanting to make the structure out of stock material (such as bars or plates), or facilitating its manufacturing. The simplest of these constraints is to enforce bounds on the component geometric parameters.  For instance, it is possible to obtain a structure that is made of bars with the same cross section (cf.\ \cite{norato2015geometry}), or to impose upper bounds on the dimensions of rectangular plates to account for commercial availability \citep{zhang2016geometry}.  Although conceptually simple, enforcing these requirements in free-form methods is difficult. 

Surprisingly, however, there have only been a few works to date that incorporate other geometric constraints. A common requirement, particularly for structures made of stock material, is to ensure that the components lie entirely within the design region, as otherwise a component that is only partially inside may require cuts (for example, through the thickness of a plate) that are impractical to manufacture.  For rectangular and cuboid design regions in 2D and 3D, respectively, this requirement is easily satisfied by imposing bounds on the design variables that determine the positions of the components.  However, for other design region shapes, particularly those that are not convex, using bounds or simple constraints on the design parameters does not work. The geometry projection method in \cite{zhang2018geometry} proposes a way to address this by creating a layer of points slightly outside of the boundary of the design region, and imposing a constraint in the optimization that the maximum projected pseudo-density in any of these points is zero.  By imposing the constraint on the pseudo-densities, this technique works for components of any shape.  An alternative to address this issue is to employ the no-overlap techniques discussed in \secref{sec:separation}, and require that the geometric components do not overlap the exterior of the design region.

Another type of geometric constraint aims to control the orientation of components.  In the MMC technique of \cite{guo2017self}, a constraint is imposed on the orientation of bars in structures fabricated via additive manufacturing to ensure their angle with respect to the print direction is smaller than the overhang angle.  In the same work, holes on a solid structure are modeled using B-splines, and a constraint is imposed on the control points positions to ensure the boundaries of the holes do not exceed the overhang angle. In the latter approach, however, another constraint is also imposed to prevent holes from merging in order to completely avoid violations of the overhang angle requirement, thus this method does not change the topology of the design. In the work by \cite{wein2018combined}, described in \secref{sec:shape_dens} and shown in \figref{fig:spline}(b), a maximum overhang on the strip boundaries is readily enforced by imposing constraints on the relative position of boundary points (which are themselves a function of the control points and the strip semi-widths). Similarly, using constraints on the positions of the control points and on the semi-widths, this method enforces slope and curvature constraints on the overall shape of the strip.

In the geometry projection technique of \cite{smith2019geometric}, a constraint is imposed on the minimum angle between any two bars (modeled as offset surfaces) to ease manufacturing.  The form of the angle constraint is somewhat similar to that of \cite{niu2018equal} in that it multiplies a term that penalizes the angles between two bars with another that measures the distance between their medial axes.  However, the angle term penalizes angles that are smaller than a prescribed value (instead of the angle difference), and the distance term is used to impose the angle constraint only between bars that are closer than a prescribed value.  Moreover, the penalty term is multiplied by the size variables of both bars, so that no angle constraint is imposed when one or both of the bars have a zero size variable.  

Feature-mapping methods often produce designs where components ---particularly bars--- intersect such that they have a `long' overlap, i.e., they are close to each other and they are near parallel. If the structure is fabricated with stock material, these overlaps make fabrication difficult. The aforementioned techniques to impose an angle constraint between bars are very effective in preventing this situation.  \cite{smith2019geometric} define a no-overlap region in the bar, corresponding to the bar minus its circular ends.  A constraint is then imposed that the sum of the projected pseudo-densities for each individual bar on the overlapping region is at most unity. In effect, these techniques ensure that a bar intersects another only at its ends. However, it may be intersected by other bars anywhere along its length.  Another technique to prevent long overlaps is the one proposed by \cite{watts2017geometric} for design of periodic lattices, where a constraint is imposed that ensures the largest distance from the endpoint of a bar's axis to the closest endpoint of another bar is at most the bar's half width. This ensures a bar only intersects another at its ends but not along the bars, thus it is more restrictive than the foregoing techniques.

Another type of geometric constraint is symmetry.  Due to the more restrictive design representation, feature-mapping methods may render designs that are not exactly symmetric even if the design region shape and the boundary conditions are such that a symmetric design is expected with free-form methods.  In the case of compliance minimization with a volume constraint, this behavior, as noted in \cite{norato2015geometry}, owes to the fact that a symmetric design for the number of available components may have a volume lower than the constraint; since adding material to this design may decrease the compliance, it may be possible to find an asymmetric design with an active volume constraint that is better. This occurs particularly for components of fixed width.  This behavior is consistent with the known fact that optimal truss designs can also be asymmetric \citep{stolpe2016truss}. To ensure a symmetric design, the geometry projection techniques in \cite{watts2017geometric,kazemitopology} use a simple strategy in which the point where the pseudo density is computed is reflected with respect to the symmetry plane before computing the signed distance to the components.  This strategy can be applied to any number of symmetry planes, as demonstrated in \cite{watts2017geometric} to design periodic truss lattices whose homogenized properties exhibit desired symmetries. It can also be applied to any component shape. In \cite{wein2018combined}, the control points that define the strips are reflected with respect to symmetry planes to obtain, e.g., square symmetric structures. 

\subsection{Applications}
\label{subsec:TO-applications}

\begin{figure}[ht!]
	\centering
	\subfloat[]   
	{\includegraphics[width=.225\columnwidth,valign=c]{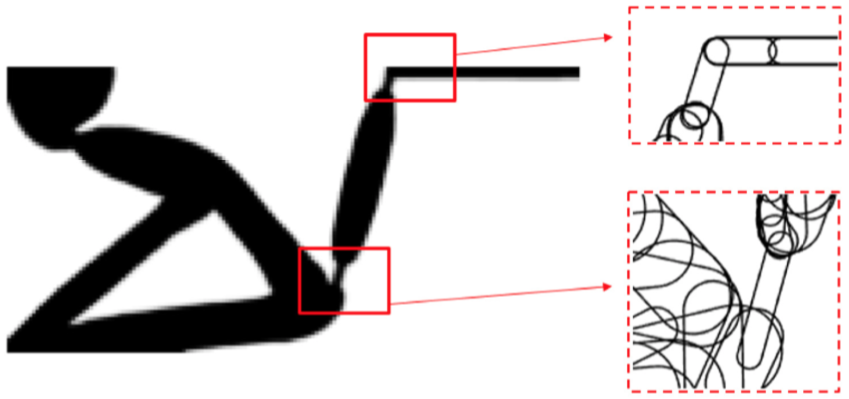} %
		\vphantom{\includegraphics[width=.225\columnwidth,valign=c]{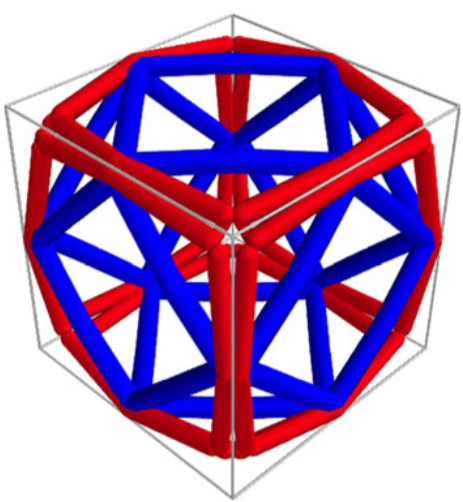}}
	} \,
	\subfloat[]
	{\includegraphics[width=.225\columnwidth,valign=c]{sect8_fig2_2017_Watts_et_al.png}} \\
	\subfloat[]   
	{\includegraphics[width=.225\columnwidth,valign=c]{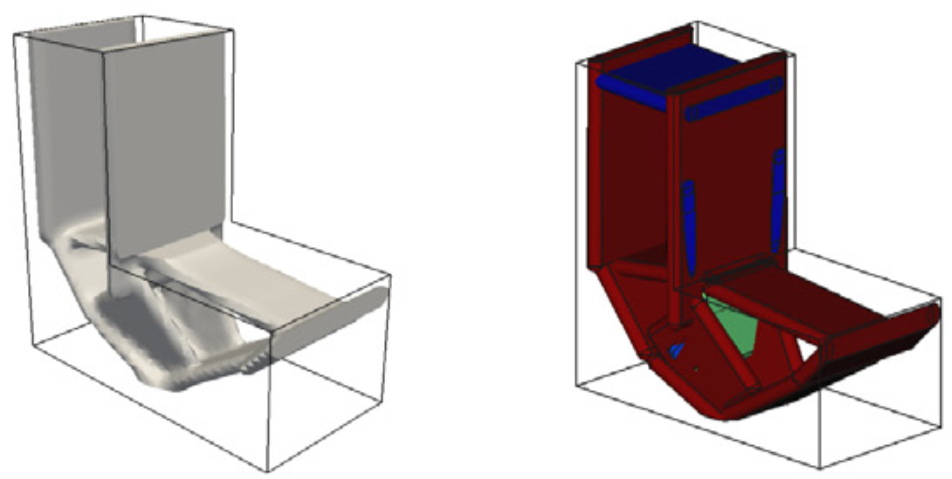}
	\vphantom{\includegraphics[width=.225\columnwidth,valign=c]{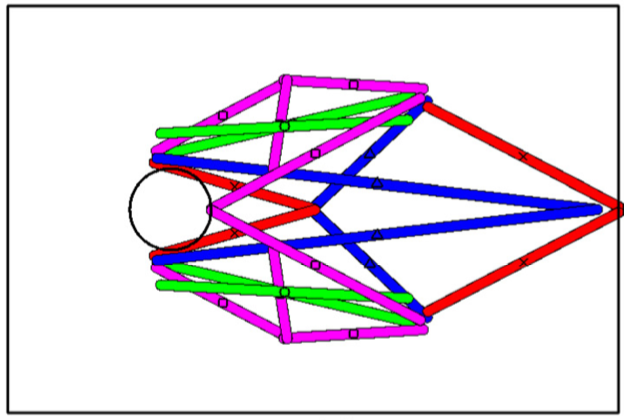}}
	}  \,
	\subfloat[]
	{\includegraphics[width=.225\columnwidth,valign=c]{sect8_fig2_2018_Kazemi_et_al.png}
	} \\	
	\subfloat[]   
	{\includegraphics[width=.225\columnwidth,valign=c]{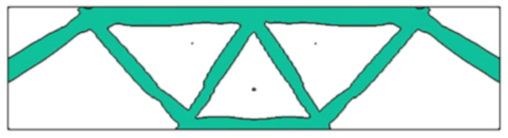}
		\vphantom{\includegraphics[width=.225\columnwidth,valign=c]{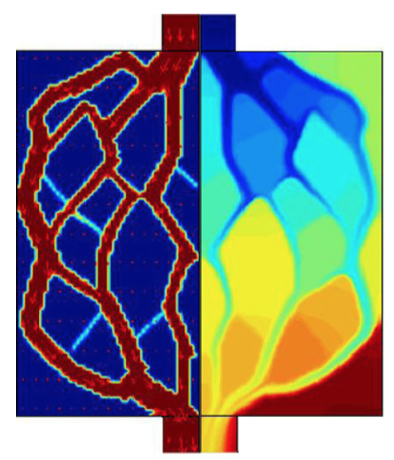}}
	}  \,
	\subfloat[]
	{\includegraphics[width=.225\columnwidth,valign=c]{sect8_fig2_2019_Yu_et_al.png}} 	
	
	\caption{Examples of different feature-mapping methods for topology optimization: a) compliant mechanism design \citep{hoang2017topology}, b) material design \citep{watts2017geometric}, c) stress constraints \citep{zhang2017stress}, d) multi-material structures \citep{kazemitopology}, e) geometric nonlinearities \citep{zhu2018structural} and f) thermal-fluid problems \citep{yu2019topology}.}
	\label{fig:sample-applications}
\end{figure}

Not surprisingly, most works in feature-mapping methods for topology optimization consider minimization of compliance with a volume constraint.  Recently, some works have incorporated other structural responses and physical regimes.  Some of the applications mentioned in this section are shown in \figref{fig:sample-applications}.

One problem which feature-mapping methods have been used for is the design of linear compliant mechanisms \citep{deng2016design,guo2016explicit,zhang2016minimum,hoang2017topology}, whereby the objective function is typically related to maximizing displacements at certain locations. The mechanisms to impose a minimum size in \citep{zhang2016minimum,hoang2017topology} are effective in preventing single-node hinges often seen in compliant mechanism design using free-form methods. Another application in elastostatic structures is an MMV method for design-dependent loading \citep{zhou2019concurrent}.  Stress constraints are incorporated in the geometry projection technique of \cite{zhang2017stress} and in the MMV method of \cite{zhang2018moving}; and stresses were minimized in the former work and in the MMC method of \cite{takalloozadeh2017implementation}, which uses the topological derivative to compute the design sensitivities of the stress function.

In the area of material design, the geometry projection method of \cite{watts2017geometric} is applied to design periodic truss latices with desired material symmetries, to maximize the effective bulk modulus or to obtain negative Poisson's ratio. 

The design of multi-material structures, where each component is made from one material from a set of available materials, is another recent application of these methods. The first work to accomplish this is the aforementioned work of \cite{watts2017geometric}.  By adapting the density-based multi-material technique of \cite{sigmund1997design} to the geometry projection framework, their method simultaneously designs the layout of the bars within the unit cell and selects the best material for each component.  The MMC method of \cite{zhang2018topology} designs structures with components made of multiple materials. However, the material for each component is assigned a priori and does not change during the optimization, and only the component layout is optimized.  \cite{kazemitopology} simultaneously design the topology and material choice of multi-material structures by adapting the discrete material optimization (DMO) technique of \cite{stegmann2005discrete} to the geometry projection method.

The MMC method of \cite{zhu2018structural} minimizes the compliance of geometrically nonlinear structures.  It uses a neo-Hookean material along with the material interpolation of \cite{wang2014interpolation}, which models low-density elements with a linear material and thus stabilizes the analysis.  The asymptotes in the optimizer ---namely, the method of moving asymptotes (MMA) of \cite{Svanberg:87:MMA}--- are adaptively updated to enforce a more conservative approximation of the optimization functions when the nonlinear finite element analysis has difficulty converging. The MMV approach in \cite{xue2019explicit} designs minimal-compliance structures and compliant mechanisms with finite deformations. To circumvent the convergence issues brought upon by low-density elements, this method adaptively removes elements inside the holes whose pseudo-density is lower than a prescribed threshold. 

An interesting application of feature-mapping methods is the design of particular types of structures that are typically manufactured using stock material, such as plates or rods.  These include the layout design of reinforcing ribs made of plates, including the geometry projection technique of \cite{zhang2017optimal} and the MMC method of \cite{zhang2018ribs}.  The latter work incorporates buckling constraints.  A subset of these works is devoted to the design of aircraft wingbox design: the method of \cite{smith2019geometry}, which aims to find an optimal layout of ribs for the wingbox, and the method of \cite{li2018multidisciplinary}, in which the layout of ribs is fixed, and MMV is used to find the optimal design of holes in the ribs to minimize fuel sloshing.  The recent work in \cite{chu2019design} uses MMC to design periodic truss cores of sandwich panels with minimal compliance.

Another application area with growing interest is dynamics.  For instance, several recent works have been devoted to design of flexible multibody systems, such as the MMC techniques of \cite{sun2018simultaneous, sun2018topologyND} to design a variable-length structure, and \cite{sun2018topology} to design a component with large motion and large deformation. \cite{sun2019topology} use MMC to design a rotating thin plate to maximize its first eigenfrequency, or maximize the gap between two consecutive eigenfrequencies. In \cite{xie2019explicit}, MMC is used to design the layout of damping patches on a vibrating plate to minimize its average kinetic energy over a frequency range. \cite{wormser2017design} design periodic lattice structures that act as phononic band gaps---a problem that is very hard to solve with free-form topology optimization techniques, which render disconnected structures that cannot be manufactured.
The MMC method has also been combined with gradient-free optimization methods to optimize structures for crash-worthiness \citep{bujny2018identification,raponi2019kriging}.

In the realm of heat transfer, combine-then-map methods have been applied to the design of thermo-elastic structures.  \cite{takalloozadeh2017implementation} use the MMC method to design minimum compliance structures subject to a uniform temperature change. \cite{sharma2017advances} minimizes the mismatch between the resulting displacement and a prescribed one for both steady-state and unsteady heat conduction (in the latter case, the prescribed displacement is specified at a given time).  The method of \cite{lohan2017topology} designs 2D structures for heat conduction in what amounts to an MMOS method.  The structure, which has higher conductivity than its surrounding medium,  is defined as the union of circles, and an element pseudo-density of 0 or 1 is assigned depending on whether the centroid is outside or inside of the structure, respectively. A gradient-free space colonization algorithm is used for the optimization. The work in \cite{yu2019topology} addresses thermal-fluid problems with the MMC method to design the pipes in a cooling device, with the objective function being a weighted sum of the thermal compliance and the power dissipation of the system. 

Finally, in the area of electromagnetics, the work of \cite{li2019generating} designs a tree-like power distribution network on an integrated circuit by minimizing the power mean voltage subject to a volume constraint. \cite{liu2019moving} use MMC to solve an inverse problem, namely to reconstruct the shapes of objects inside a body from electrical impedance measurements made on the surface of the body. 

\section{Discussion}
\label{sec:discuss}

Feature-mapping is a powerful and promising approach in structural optimization. As detailed in this review, there is a wide range of choices to accomplish the two main ingredients of these methods, namely the mapping of high-level parameters onto the fixed grid for analysis and the combination of features. 

We believe there is still sufficient room for further exploring the aspects of feature mapping and to advance novel techniques, which will likely happen based on applications of feature mapping to solve problems not (or not as efficiently) solvable by established methods, notably density-based and level-set topology optimization.


One area that we identified for potential further work is to consider more problems where high-level geometric constraints are essential for the application. For example, this could include manufacturing constraints and problems where connectivity between two or more points is essential.

An often discussed benefit of feature-mapping approaches is the potentially straightforward transfer of optimized designs to the CAD environment (due to use of a high-level geometric parameterization). However, there is little evidence to demonstrate this benefit in the literature.

The issue of local minima and initial design-dependent solutions is currently one of the key challenges of feature-mapping methods, particularly for hybrid methods and topology optimization. Some remedies to this problem are suggested in \secref{subsec:TO-design-space}. However, we consider there is scope for new methods and ideas to solve (or at least mitigate) this issue.

Feature-mapping methods for topology optimization have mostly considered Boolean unions of components or holes. While the use of the Boolean union in conjunction with other operations (e.g., intersection and subtraction) was already conceived in the pioneering work of \cite{mei2008feature}, this has not been demonstrated in methods for topology optimization.  Hybrid methods, on the other hand, have demonstrated this possibility. This is an important capability because it allows to obtain more complex shapes (albeit still with high-level parameters), and because it mirrors the way in which CAD systems construct complex geometries. We thus believe there is room for topology optimization methods to incorporate and demonstrate this capability. 

Finally, in preparing this review we have noticed that there is a level of reinvention of methods and techniques. We hope this review helps researchers understand the different aspects of feature-mapping and the various techniques proposed to-date in order to avoid duplication and to duly credit existing work.

\section*{Acknowledgement}
We thank Dr. Lukas Pflug from the Department of Mathematics at the  Friedrich-Alexander-Universit\"at Erlangen-N\"urnberg (FAU), Germany, for fruitful discussion and support. 

The initiative for this review goes back to critical yet constructive comments by Prof. Kurt Maute, from the University of Colorado Boulder, USA.

We also thank Prof. Horea Ilies from the University of Connecticut, USA, for guidance and insight into some of the geometric aspects of this work.

\bibliographystyle{abbrvnat}     
\small
\setlength{\bibsep}{0pt plus 0.3ex}
\bibliography{Review_arXiv}

\end{document}

%% file: fabian_macros.tex
\usepackage{amsmath,amssymb,amsfonts,amsbsy}
\usepackage{mathptmx}
\usepackage{setspace}
\usepackage{bm}

\newcommand{\bs}{\boldsymbol}

\newcommand{\intd}{\mathrm{d}\,}

\newcommand{\sens}[2]{\frac{\partial{#1}}{\partial{#2}}}
\newcommand{\drho}[1]{\frac{\partial{#1}}{\partial \rho_e}}
\newcommand{\total}[2]{\frac{\text{d}\,{#1}}{\text{d}\,{#2}}}

\newcommand{\brho}{{\bs \rho}}

\newcommand{\blmbd}{{\bs \lambda}}

\newcommand{\rhomin}{\rho_{\min}}

\newcommand{\figref}[1]{Fig.~\ref{#1}}
\newcommand{\secref}[1]{Sec.~\ref{#1}}

\newcommand{\eqnref}[1]{(\ref{#1})}